\documentclass[3p,number,sort&compress,times,lefttitle]{elsarticle2}

\journal{Computer Methods in Applied Mechanics and Engineering}

\usepackage[utf8]{inputenc}
\usepackage{hyphenat}
\usepackage{subcaption}
\usepackage{multirow}
\usepackage{graphicx}
\usepackage{bm}
\usepackage{amsmath,amssymb,amsfonts}
\usepackage{setspace}
\usepackage{tikz}
\usetikzlibrary{arrows, arrows.meta}
\usepackage[colorlinks=true]{hyperref}
\usepackage{xcolor}
\usepackage{cancel}

\renewcommand{\Re}{\mathbb{R}}

\newcommand{\sref}[1]{Section~\ref{#1}}

\newcommand{\fref}[1]{Fig.~\ref{#1}}
\newcommand{\eref}[1]{(\ref{#1})}

\newcommand{\vm}[1]{{\bm{#1}}}
\newcommand{\vx}{\vm{x}}

\DeclareMathOperator*{\argmin}{argmin}
\DeclareMathOperator*{\argmax}{argmax}
\DeclareMathOperator*{\argsup}{argsup}

\newcommand{\revised}[1]{{#1}} 

\graphicspath{ {./figs/}
             }

\begin{document}

\doublespacing

\title{Variational formulation based on duality to solve partial 
       differential equations: Use of 
       B-splines and machine learning approximants}

\author[1]{N.\ Sukumar\corref{cor}}
\ead{nsukumar@ucdavis.edu}

\author[2]{Amit Acharya}
\ead{acharyaamit@cmu.edu}

\cortext[cor]{Corresponding author}

\address[1]{Department of Civil and Environmental Engineering,  One Shields Avenue,
University of California, Davis, CA 95616, USA}

\address[2]{Department of Civil and Environmental Engineering, and
            Center for Nonlinear Analysis,
            Carnegie Mellon University, Pittsburgh, PA 15213, USA}

\begin{abstract} 
Many partial differential equations (PDEs) such as Navier--Stokes equations
in fluid mechanics, 
inelastic deformation in solids, and transient parabolic and hyperbolic 
equations do not have an exact, primal variational structure. Recently, a variational principle based on
the dual (Lagrange multiplier) field was proposed.
The essential idea in this approach is to treat the given
\revised{PDEs} as constraints, and to invoke an arbitrarily chosen 
auxiliary potential with strong convexity properties to be optimized.  
\revised{On requiring the vanishing of the gradient of the Lagrangian with respect to
the primal variables, a mapping from the dual to the
primal fields is obtained.}
This leads to requiring a convex dual functional to be minimized subject to 
Dirichlet boundary conditions on dual variables, with the guarantee that even 
PDEs that do not possess a variational structure in primal form can be solved 
via a variational principle. The vanishing of the
first variation of the dual functional is, up 
to Dirichlet boundary conditions on dual fields, the weak form of the primal PDE problem with the dual-to-primal change of variables incorporated. 
We derive the dual weak form for the linear, one-dimensional,  
transient convection-diffusion equation. A Galerkin
discretization is used to obtain the discrete equations, with
the trial and test functions \revised{in one dimension}
chosen as linear 
combination of either 
\revised{shallow neural networks with Rectified Power Unit (RePU)
activation functions} or B-spline basis functions; the corresponding stiffness matrix is symmetric. 
For transient problems, a 
space-time Galerkin implementation is used with tensor-product
B-splines as approximating functions. 
Numerical results are presented for 
the steady-state and transient convection-diffusion equations and transient heat conduction. The proposed method
delivers sound accuracy for ODEs and PDEs and 
rates of convergence are established in the $L^2$ norm and $H^1$ seminorm for the steady-state 
convection-diffusion problem.
\end{abstract}

% limit to six

\maketitle

\section{Introduction}\label{sec:intro}
Many classes of ordinary differential equations (ODEs)
and partial differential equations (PDEs) do not possess a 
natural, exact, variational structure:\footnote{In the sense of 
being derived, along with all associated boundary and initial 
conditions, as Euler--Lagrange equations of some functional of 
fields defined on space-(time) without approximation.} convection-diffusion and 
Navier--Stokes equations in fluid mechanics (cf.~\cite{ghoussoub2009anti,ortiz2018variational}), inelastic deformation of solids (cf.~\cite{ortiz1999nonconvex,ortiz1999variational,petryk2020quasi,carstensen2002non}), 
and time-dependent parabolic (heat) and 
hyperbolic (wave)~\cite{gurtin1964variational} problems, to name just a few; an insightful discussion of the issues involved in the context of the equations of continuum mechanics and classical field theories can be found in \cite[Sec.~1]{seliger1968variational}. 
For such problems,
it is desirable to provide a variational
principle that can lead to alternate solution strategies with
the finite element method, meshfree methods, virtual element methods, B-spline
approximations, and neural networks.
Recently, 
a scheme for generating dual variational principles for such problems was proposed in~\cite{Acharya:2022:VPN} which is
cast in terms of dual (Lagrange multiplier) fields. In~\cite{Acharya:2022:VPN, Acharya:2023:DVP, Acharya:2024:HCC,Singh:2024:HCN,Acharya:2024:VDS,Acharya:2024:APD}, that scheme has been placed within the broader context of the prior efforts for developing the variational principles mentioned above,
including the method of least squares. 
There exists a superficial similarity between the dual formulation 
to solve differential-algebraic systems and 
the adjoint method~\cite{Errico:1997:WAM,Cao:2003:ASA}, 
the latter used in constrained optimization
to compute the derivatives of a function with respect to
unknown parameters. Their similarities and differences 
are noted in~\ref{appendix:A}. 

The methodology developed and demonstrated in~\cite{Acharya:2022:VPN,Acharya:2023:DVP,Acharya:2024:VDS,Acharya:2024:APD,
Acharya:2024:HCC, Singh:2024:HCN,Kouskiya:2024:HCH,Kouskiya:2024:IBD} applies to  quasistatic and dynamic, conservative and dissipative 
PDEs that arise in continuum mechanics. 
In a notable prior work~\cite{seliger1968variational},
a strategy for deriving variational principles for a limited set\footnote{In the words of Seliger and Whitham describing their work~\cite{seliger1968variational}: ``We are still short, however, of a general theorem stating the conditions under which a variational principle can be
found for any given system of equations and of an automatic fool-proof method of
producing it.''} of PDE systems from continuum mechanics was 
devised; the approach was successfully applied to nonlinear elastodynamics in the spatial setting, as well as the compressible Euler equations, Maxwell's equations, and those of the collisionless plasma.
The essential idea, as introduced in~\cite{Acharya:2022:VPN}, is to treat the given \revised{PDEs} as constraints for an 
arbitrarily chosen auxiliary potential with strong convexity properties that 
is optimized. The Lagrange multipliers associated with the constraints are 
referred to as the dual fields.
On requiring the vanishing of the gradient of the Lagrangian with respect to
the primal variables, a dual-to-primal (DtP) mapping from the dual to the
primal fields is obtained. This leads to requiring 
a convex dual functional to be minimized subject to 
Dirichlet boundary conditions on dual variables, with the guarantee that even 
PDEs that do not possess a variational structure can be solved as 
Euler--Lagrange equations of the dual functional. The first variation 
(set to zero) of the dual functional is, up to Dirichlet boundary conditions on dual fields, the weak form of the primal PDE problem with the dual-to-primal change of variables incorporated.
The Euler--Lagrange equations of the dual functional are shown to be
locally degenerate elliptic, regardless of the properties of the primal 
system~\cite{Acharya:2024:HCC} (for elliptic primal problems, the dual problem is elliptic). Our ideas are related to those of Brenier~\cite{Brenier:2010:HCN,Brenier:2018:IVP}, and 
generalizes~\cite{Brenier:2018:IVP} in various 
ways~\cite{Acharya:2024:VDS}.

When posed in primal form, the weak formulation for convection-diffusion 
problems and incompressible Navier--Stokes equations require special 
treatment via the use of upwinding and/or stabilization schemes. 
The streamline upwind/Petrov--Galerkin (SUPG) finite
element method produces a
variationally consistent scheme that delivers stable and accurate
solutions~\cite{Brooks:1982:SUP}. In the Galerkin/least squares (GLS)
method~\cite{Hughes:1989:FEF,tezduyar1991stabilized}, 
a given PDE (linearized, when not linear) is squared and the first variation of this functional, after scaling with a tuning parameter $\tau$, is added to the standard Galerkin weak form. This enhances the stability of the standard Galerkin
method. However, the weak formulation 
in GLS recovers the strong form only when $\tau \to 0$. The
derivation of variational methods for multiscale phenomena and the connection
of subgrid models therein to enriched bubble functions and to 
stabilized methods is presented in~\cite{Hughes:1995:MPG}. 

An alternative approach is taken in the least-squares finite element
method (LSFEM)~\cite{Bochev:2009:LFE}: if $Au = 0$ is a PDE system ($A$ can be a nonlinear operator),
then the LSFEM functional $||Au||^2$ in an appropriate norm (e.g.,
$L^2$) is minimized.
On setting the first variation of this functional to zero leads
to the variational equations. 
A positive attribute of LSFEM is that it results in a symmetric
stiffness matrix.
However, though the solution of the strong form minimizes the LSFEM functional, notably for nonlinear problems
it is not guaranteed that all solutions to the
Euler--Lagrange equations of the LSFEM functional are solutions of the original PDE.\footnote{This was understood, and explained in~\cite[Sec.~1, p.~2]{seliger1968variational}.} In addition, there is no 
systematic procedure to incorporate boundary conditions 
\revised{in the variational statement}. 

From a computational viewpoint, the availability of a variational principle is beneficial; for instance, 
regardless of nonlinearities, discretizations based on a variational principle result in symmetric stiffness/Jacobian matrices. To solve nonlinear PDEs (e.g., Allen--Cahn and Burgers' equations), 
transient wave phenomena as well as chaotic dynamical systems, the
strong form is used in physics-informed neural 
networks (PINNs)~\cite{Raissi:2019:PIN,Karniadakis:2021:PIM,Cuomo:2022:SML},
Kolomogorov--Arnold networks (KANs)~\cite{Liu:2024:KAN,Koenig:2024:KAN}, and neural networks 
based on the
Kolmogorov superposition theorem~\cite{Guilhoto:2024:DLA}. 
Use of duality principles can provide alternate routes to solutions via functional minimization with
neural networks to better capture the underlying physics. For instance, deep Ritz~\cite{E:2018:DRM} provide certain
advantages to solving variational formulations with neural 
networks: only Dirichlet boundary conditions are required to be imposed (advances in~\cite{Sukumar:2022:EIB} can be leveraged), and as a consequence, admissible neural network approximations ensure that the potential energy functional
for the Poisson equation or in nonlinear elasticity is bounded, and
since only first-order derivatives are computed, the computational overhead is
less than when using the strong form for second-order PDEs.
With the dual approach, smooth approximants for the dual
fields are desirable, which are readily constructed using
neural network approximants. 
In addition, since imposing homogeneous Dirichlet boundary conditions on dual fields 
suffices, it is easier to  
construct (a priori) ansatz
that are kinematically admissible.
Furthermore, a transient problem
in the dual approach is posed as a boundary-value
problem with a terminal
Dirichlet boundary condition that can be solved
using a space-time neural network discretization in finite time 
slices (see~\cite{Kouskiya:2024:HCH} 
for time-slicing in the dual context).

In previous studies, linear $C^0$ finite elements have been used in the
dual variational formulation to solve linear and nonlinear PDEs 
and an ODE: transient heat conduction, linear
transport and Euler's nonlinear ODE system for the motion of a rigid body~\cite{Kouskiya:2024:HCH}; nonconvex one-dimensional
elastostatics and elastodynamics~\cite{Singh:2024:HCN}; and (inviscid) Burgers 
equation~\cite{Kouskiya:2024:IBD}. 
Since in general 
the primal fields depend on the derivative of dual fields, an $L^2$ 
projection of the DtP generated primal fields onto
$C^0$ finite element spaces was used to reconstruct continuous
primal fields.  In this paper, we adopt smooth 
neural network approximations and high-order 
B-spline approximating functions in the dual
variational formulation.
In the numerical implementation, we show that 
it is desirable to use smooth, higher order
approximations for the 
dual fields if the primal field is $C^0$. However, for the
pure transport problem, a dual field that is $C^0$, obtained by adding $C^0$ basis functions to the $C^k$ ($k \ge 1$) approximants already employed may be 
beneficial since it accommodates discontinuous solutions in the primal field.
In one dimension, shallow neural networks with 
Rectified Power Unit (RePU) activation 
function~\cite{He:2024:EAP,He:2024:DNN} and univariate
B-splines~\cite{deBoor:2001:PGS} are used, and bivariate tensor-product B-splines are adopted in two dimensions.  

The remainder of this 
paper is structured as follows. In~\sref{sec:formalism},
we introduce the general formalism of the dual variational principle and 
derive the dual function for a system of nonlinear equations.
The solution procedure via the duality approach is presented for least 
squares minimization, solving a 
system of two quadratic equations, and to obtain nonnegative
solutions for an underdetermined system of linear equations 
via entropy maximization. In~\sref{sec:dual}, we first
derive the dual formulation for an
initial-value problem; the weak form of the associated dual second-order
boundary-value
problem is 
also derived, with the latter further elaborated in~\ref{appendix:B}.
Then we present the dual formulation
for the transient convection-diffusion equation. Numerical
discretization of the variational form, with expressions for the 
stiffness matrix and force vector, are provided 
in~\sref{sec:varform_cd}.
In~\sref{sec:numerical_examples}, numerical results are
presented for the Laplace equation, one-dimensional 
steady-state convection-diffusion, transient convection-diffusion, and
transient heat equations. Accurate results are obtained with sound convergence,
and we summarize the main findings from this work in~\sref{sec:conclusions}.
            
\section{General formalism for the finite-dimensional 
         case}\label{sec:formalism}
Following \cite[Sec.~1]{Acharya:2024:HCC}, let 
$\vm{G}(\vm{U}) = \vm{0}$ represent a system of linear equations, 
nonlinear equations, ODEs or PDEs, where 
$\vm{G} : \Re^n \to \Re^m$. To fix
ideas, we first present
the dual formalism to solve a system of nonlinear equations.  To formulate
an optimization problem that is preferably convex, consider the objective 
function {$\vm{U} \mapsto H(\vm{U})$}, where {$H(\vm{U})$} is a convex 
(choice is flexible) auxiliary potential that is tailored to dominate the 
nonlinearities in the primal problem.

Let $\vm{\lambda} \in \Re^m$ be the Lagrange multiplier vector that is 
associated with the constraints. We now optimize (minimize) $H$
subject to the constraint $\vm{G}(\vm{U}) = \vm{0}$, and can write the 
Lagrangian and the stationary conditions as:
\begin{subequations}\label{eq:minH}
   \begin{align}
     {L}_H (\vm{U},\vm{\lambda}) &:=  
              H(\vm{U}) + \vm{\lambda} \cdot \vm{G}(\vm{U}),\label{eq:minH_a} \\
             \nabla_{\vm{U}} {L}_H (\vm{U},\vm{\lambda}) 
          &=  \nabla H(\vm{U}) + \vm{\lambda} \cdot \nabla \vm{G}(\vm{U}) 
              = \vm{0}, \label{eq:minH_b} \\
     \nabla_{\vm{\lambda}} {L}_H(\vm{U},\vm{\lambda}) &=  
              \vm{G}(\vm{U}) = \vm{0} . \label{eq:minH_c}
   \end{align}
\end{subequations}
In~\eref{eq:minH}, $\vm{U} \in \Re^n$ are the primal variables and 
$\vm{\lambda} \in \Re^m$ are the dual variables. In the critical point formulation, which shows the consistency of our scheme as generating solutions to the primal problem, one solves~\eqref{eq:minH_b} for $\vm{U}$ as a function of $\vm{\lambda}$, the dual-to-primal 
mapping $\vm{U}_H(\vm{\lambda})$, and then look for $\vm{\lambda}^\ast$ that solves~\eqref{eq:minH_c} in the form 
$\vm{G} \bigl( \vm{U}_H(\vm{\lambda}^\ast) \bigr) = \vm{0}$. Defining the 
function $\vm{U}_H$ is facilitated 
by the choice of $H$. In the related convex optimization formulation, 
one solves the maximization problem
\begin{equation}
   \sup_{\vm{\lambda}} \inf_{\vm{U}} 
  {L}_H(\vm{U},\vm{\lambda}) = \sup_{\vm{\lambda}} S(\vm{\lambda}),
\end{equation}
where $S(\vm{\lambda})$ is the dual function, which is 
concave, since it is the pointwise infimum of a family of affine functions, regardless of the nonlinearity of $\vm{G}$. A maximizer (and a critical point) $\vm{\lambda}^\ast$ of the dual function
$S(\vm{\lambda})$ is defined as
\begin{equation}\label{eq:lambda_ast}
   \vm{\lambda}^\ast = \argsup \limits_{\vm{\lambda}} S(\vm{\lambda}).
\end{equation}
Note that the convex optimization problem does not need a notion of a DtP mapping (in fact, a concave dual problem can be obtained even with $H = 0$). However, for a meaningful scheme that addresses the question of finding solutions to the primal problem, a choice of a strictly convex function $H$ facilitates the definition of a DtP map. The maximizer of the dual problem in this case, with some consideration on the % of 
smoothness of $S$ at the maximizer $\vm{\lambda}^\ast$ as well as the smoothness of the DtP map, defines a primal solution as 
\begin{equation*}
\vm{U}^\ast = \vm{U}_H(\vm{\lambda}^\ast).
\end{equation*}
In this connection, we note that while the hypotheses of the Implicit Function theorem, primarily related to the invertibility of the Hessian $\nabla^2_U L_H(\cdot, \vm{\lambda})$  in $\mathcal{O}_n \times \mathcal{O}_m \subset \mathbb{R}^{n} \times \mathbb{R}^m$, ensures the existence of the function 
$\vm{U}_H: \mathcal{O}^m \to \mathbb{R}^n$ (at least locally),
the existence of such a function does not imply
$\inf \limits_{\vm{U}} L_H(\vm{U}, 
\vm{\lambda}) = L_H(\vm{U}_H(\vm{\lambda}), \vm{\lambda})$.
Strict convexity of $L_H(\cdot,\vm{\lambda})$ in $\mathcal{O}_n$ for each 
$\vm{\lambda} \in \mathcal{O}_m$ is a guarantee of this fact. 

To verify that $\vm{U}^\ast$ satisfies~\eref{eq:minH_c}, that is
$\vm{G}(\vm{U}^\ast) = \vm{0}$, assume that in a neighborhood of 
$\vm{\lambda}^\ast$, $S(\vm{\lambda}) = \inf \limits_{\vm{U}} L_H(\vm{U}, 
\vm{\lambda}) = L_H \bigl( \vm{U}_H(\vm{\lambda}), \vm{\lambda}
\bigr)$ holds, and observe that solving the optimization
problem posed in~\eref{eq:lambda_ast} by
 setting the gradient of $S$ to vanish and using~\eref{eq:minH_b}
 yields
 \revised{
\begin{equation*}
\left. \vm{0} = \dfrac{\partial S}{\partial \vm{\lambda}} 
             \right|_{ ( \vm{\lambda}\,=\,\vm{\lambda}^\ast ) } 
             = 
             \left.
             \cancelto{\vm{0}}{
             \dfrac{\partial L_H}{\partial \vm{U}}
             }
            \cdot
            \dfrac{\partial \vm{U}}{\partial \vm{\lambda}}
            \right|_{( \vm{U}\,=\,\vm{U}_H(\vm{\lambda}^\ast), \,
            \vm{\lambda} \,=\, \vm{\lambda}^\ast ) }
            + 
            \left.
            \dfrac{\partial {L}_H}{\partial \vm{\lambda}} 
            \right|_{( \vm{U}\,=\,\vm{U}_H(\vm{\lambda}^\ast), \,
            \vm{\lambda} \,=\, \vm{\lambda}^\ast ) }
            = 
            \left.
            \dfrac{\partial {L}_H}{\partial \vm{\lambda}} 
            \right|_{( \vm{U}\,=\,\vm{U}_H(\vm{\lambda}^\ast), \,
            \vm{\lambda} \,=\, \vm{\lambda}^\ast ) }
            = \vm{G} \bigl( \vm{U}_H(\vm{\lambda}^\ast) \bigr) 
            = \vm{G}(\vm{U}^\ast), 
\end{equation*}
}
which establishes that the primal problem is solved via dual maximization.

We now apply the duality principle to three test cases: system of linear 
equations, solution of two quadratic equations and nonnegative solutions 
for an underdetermined system of linear equations.

\subsection{Linear system of equations}\label{subsec:LS}
Following \cite[Sec.~2]{Acharya:2023:DVP}, we seek $\vx \in \Re^n$ that solves the system of linear 
equations:
\begin{equation}\label{eq:Ax=b}
\vm{A} \vx = \vm{b}
\end{equation}
by the dual methodology, where $\vm{A} \in \Re^{m \times n}$ and $\vm{b} \in \Re^m$. If
$m = n$, then~\eref{eq:Ax=b} has a unique solution if $\vm{A}$ has full
rank of $n$. If $n > m$ ($n < m$), the linear system is underdetermined (overdetermined). 

\revised{The simplest choice of a convex potential (density) is a 
quadratic function, since its differential is linear and will lead to an explicit DtP map for algebraic problems up to degree 2 and for the 
linear PDEs that are considered in this work. To 
solve~\eqref{eq:Ax=b} with the dual approach, 
we choose $H(\vx) = \frac{1}{2} \vx^\top \vx$.}
Let $\vm{\lambda} \in \Re^m$ be the Lagrange
multiplier vector. We form the Lagrangian ${L}(\cdot,\cdot)$,
and pose the optimization problem as:
\begin{equation}\label{eq:H_Ax}
\max_{\vm{\lambda}} \min_{\vx}
\left[ {L}(\vx,\vm{\lambda}) :=  
H(\vx) + \vm{\lambda}^\top (\vm{b} - \vm{A}\vx ) \right].
\end{equation}
On setting {$\frac{\partial {L}}{\partial \vx} = \vm{0}$, we obtain
the DtP map: 
\begin{equation}\label{eq:dtp_Ax}
      \vx := \vx_H = \vm{A}^\top \vm{\lambda},
\end{equation}
which on substituting in~\eref{eq:H_Ax} yields the dual function:
\begin{equation}
S(\vm{\lambda}) = \min_{\vx} {L}(\vx,\vm{\lambda}) 
                = {L}(\vx_H,\vm{\lambda}) 
                = -\frac{1}{2} \vm{\lambda}^\top \vm{A} \vm{A}^\top \vm{\lambda} + \vm{\lambda}^\top \vm{b} ,
\end{equation}
which is a quadratic form in $\vm{\lambda}$. Now, the dual maximization
problem is:
\begin{equation}
\vm{\lambda}^\ast = \argmax \limits_{\vm{\lambda}} S(\vm{\lambda}),
\end{equation}
which is solved using the necessary (stationary)
condition $\frac{\partial S}{\partial \vm{\lambda}} = \vm{0}$
to yield
\revised{
\begin{equation}\label{eq:lambdadual_Ax}
\big(\vm{A} \vm{A}^\top \bigr) \vm{\lambda} = 
\vm{b},
\end{equation}
whose solution is $\vm{\lambda}^\ast$. For a primal system~\eqref{eq:Ax=b} that is consistent, any dual solution satisfying the above equation (regardless of the method used to generate it, e.g., using an inverse or a generalized inverse in some specific sense if $\vm{A} \vm{A}^\top$ is not invertible) suffices to guarantee a primal solution.
Hence, on using~\eref{eq:dtp_Ax}, the solution for $\vx$ is:
\begin{equation}\label{eq:xdual_Ax}
\vx_H = \vm{A}^\top \vm{\lambda}^\ast .
\end{equation}
}

\revised{
Note that the least squares (minimum $L^2$ norm) solution of~\eref{eq:Ax=b} is:
\begin{equation}\label{eq:xLS_Ax}
\vx_{\textrm{LS}} 
= \vm{A}^{\dagger} \vm{b},
\end{equation}
where $\vm{A}^\dagger$ is the Moore--Penrose pseudoinverse of $\vm{A}$. 
When~\eqref{eq:Ax=b} is consistent,
the least squares solution 
is identical to the dual solution in~\eqref{eq:xdual_Ax} for the
considered $H$.
However, in the context of obtaining a solution to $\vm{A}\vx = \vm{b}$,
it provides a spurious `solution' when the linear system 
in~\eref{eq:Ax=b} is  
not consistent, i.e.,
$\vm{b}$ is not in the column space of 
$\vm{A}$. On the other hand,  noting that the column space of $\vm{A}$ and 
$\vm{A}\vm{A}^\top$ are identical, a solution via the duality approach (with quadratic $H$) exists 
if and only if 
the linear system in~\eref{eq:Ax=b} has at least one solution.}

\subsection{System of quadratic equations}
Consider the following system of quadratic equations:
\begin{subequations}\label{eq:quadratic}
\begin{align}
x^2 + y^2 &= 3, \qquad x^2 - y^2 = 1, \label{eq:quadratic_a} \\
\intertext{with exact solution sets:}
(x,y) &= (\pm\sqrt{2},\pm 1) \label{eq:quadratic_b}.
\end{align}
\end{subequations}
Since the solution to~\eref{eq:quadratic_a} is nonunique, we 
choose $H(x,y;\bar{x},\bar{y},\beta) = 
\revised{\beta} \,
\bigl[ (x-\bar{x})^2 + (y-\bar{y})^2 \bigr] $
as the convex potential, where $(\bar{x},\bar{y})$ is 
referred to as a base state and $\beta \in \Re$ is a constant.
\revised{The choice of the base state permits to judiciously target a particular solution of the primal
problem, if the primal problem has multiple solutions as is the case for the system~\eqref{eq:quadratic} (see also~\cite{Singh:2024:HCN}).}
For an example comprising a system of a quadratic and a linear 
algebraic equation, see~\cite[App.~B]{Kouskiya:2024:IBD}. 
The Lagrangian saddle-point problem 
\revised{associated with~\eqref{eq:quadratic} is}:
\begin{equation}\label{eq:L_quadratic}
\max_{\vm{\lambda}} \min_{x,y}
\left[ {L}(\vx,\vm{\lambda}) :=  H(x,y;\bar{x},\bar{y},\beta) 
+ \lambda_1 (3 - x^2 - y^2) + \lambda_2 (1 - x^2 + y^2) \right] .
\end{equation}
On setting $\frac{\partial {L}}{\partial \vx} = \vm{0}$, we obtain
the DtP map: 
\begin{equation}\label{eq:dtp_quadratic}
x := x_H = \frac{\beta \bar{x}}{\beta - \lambda_1 - \lambda_2}, \quad
y := y_H = \frac{\beta \bar{y}}{\beta + \lambda_1 - \lambda_2} .
\end{equation}
\revised{The dual functional is convex. The selection of $\beta$ provides a weighting so that in the 
dual space ($\lambda_1$ and $\lambda_2$), the singularity in the DtP map occurs far away from the origin. This permits an optimization algorithm to converge with a zero initial guess for the dual variables.}
On substituting $x_H$ and $y_H$ in~\eref{eq:L_quadratic}, we obtain
the dual function:
\begin{equation}
S(\vm{\lambda}) = L \bigl( x_H(\vm{\lambda}),y_H(\vm{\lambda}),\vm{\lambda} \bigr).
\end{equation}

For $\beta = 10$, and the base state $(\bar{x},\bar{y}) = (1,1)$, the
plot of $S(\vm{\lambda})$ is presented in~\fref{fig:g_quadratic}. We observe
that $g$ is maximized at $\lambda_1 = \lambda_2 \approx 1.467$, and hence
the primal solution from~\eref{eq:dtp_quadratic} is:
$x = x_H = \sqrt{2}$ and $y = y_H = 1$, which is in agreement with 
one of the exact solutions given in~\eref{eq:quadratic_b}.
\begin{figure}[!hbt]
\centering
\includegraphics[width=0.55\textwidth]{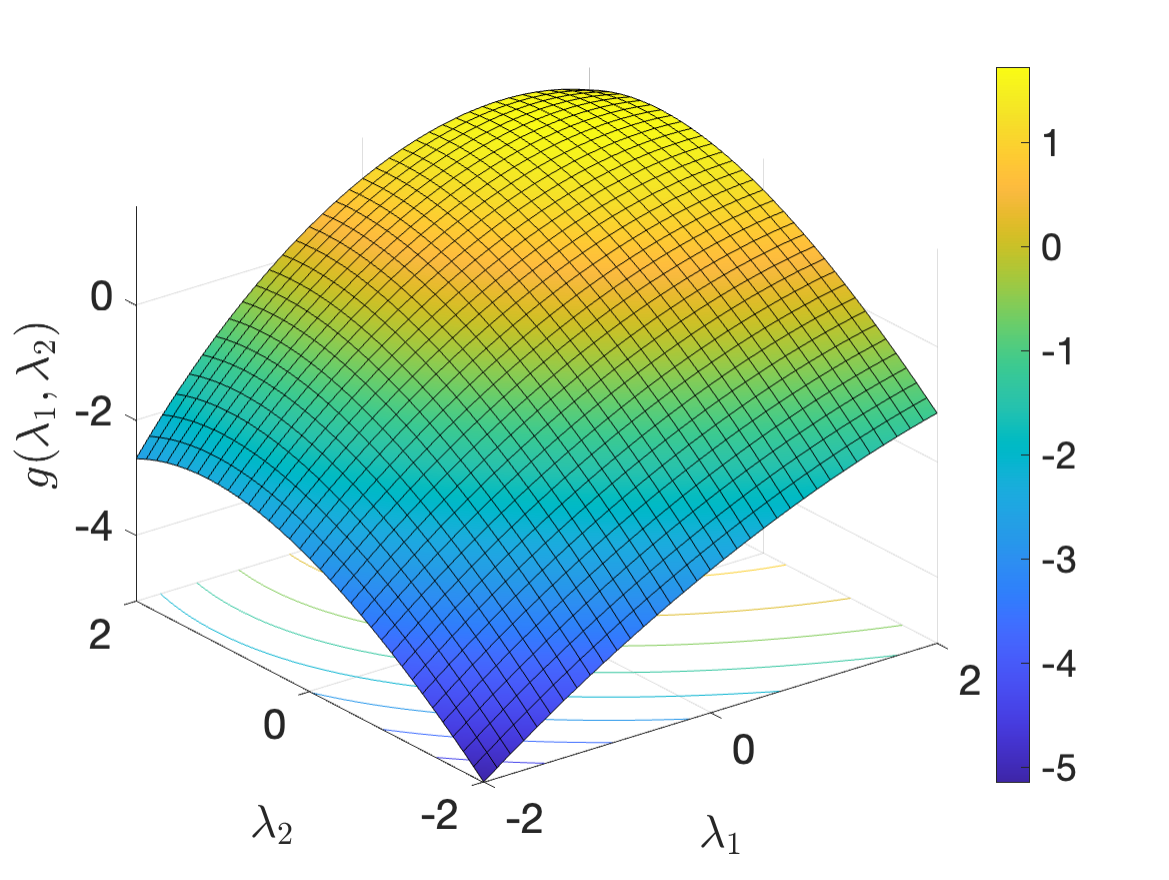}
\caption{Plot of the dual function $S(\vm{\lambda})$ 
         for the system of two quadratic equations.}\label{fig:g_quadratic}
\end{figure}

\subsection{Entropy shape functions for a polygon}
Consider a convex polygon $P \subset \Re^2$ with nodal (vertex) coordinates
$\{\vx_i\}_{i=1}^n$, where $\vx_i \equiv (x_i,y_i)$.
Let $\vm{\phi}(\vx) : P \to \Re_+^n$ be the $n$ nonnegative 
generalized barycentric coordinates (shape functions)
for a polygon~\cite{Floater:2015:GBC}.  For each fixed $\vx \in P$, these 
nonnegative shape functions ($\phi_i(\vx) \geq 0$)
satisfy the constant and linear reproducing conditions, namely
\begin{subequations}\label{eq:reproducing_conditions}
\begin{align}
\sum_{i=1}^n \phi_i(\vx) = 1, \label{eq:reproducing_conditions_a} \\
\sum_{i=1}^n \phi_i(\vx) \vx_i = \vx. \label{eq:reproducing_conditions_b} 
\end{align}
\end{subequations}
If $n > 3$ (polygon with more than three vertices), the solution 
to~\eref{eq:reproducing_conditions} is nonunique. A feasible solution using 
the convex potential $H$ as the Shannon
entropy~\cite{Shannon:1948:MTC,Jaynes:1957:ITA} is 
proposed in~\cite{Sukumar:2004:COP}. The primal optimization problem is 
posed as:
\begin{equation}\label{eq:H_entropy}
\max_{ { \vm{\phi} \in \Re_+^n }} 
  \left[ H(\vm{\phi})  := - \sum_{i=1}^n 
                           \phi_i  \ln \phi_i \right],
\end{equation}
subject to the linear constraints in~\eref{eq:reproducing_conditions}.

If $\lambda_0$ and $\vm{\lambda} \in \Re^2$ are the Lagrange multipliers 
associated with the constraints in~\eref{eq:reproducing_conditions_a} 
and~\eref{eq:reproducing_conditions_b}, respectively, then the
Lagrangian saddle-point problem is:
\begin{equation}\label{eq:L_maxent}
\min_{\lambda_0,\vm{\lambda}} \max_{\vm{\phi} \in \Re_+^n} 
\left[ {L}(\vm{\phi}, \lambda_0, \vm{\lambda}) = 
- \sum_{i=1}^n \phi_i \ln \phi_i +
\lambda_0 \Bigg( 1 - \sum_{i=1}^n \phi_i  \Bigg) + \vm{\lambda} \cdot \sum_{i=1}^n \phi_i (\vx - \vx_i) \right].
\end{equation}
On using the stationary conditions, we obtain
\begin{equation}\label{eq:dL_maxent}
\dfrac{\partial {L}}{\partial \vm{\phi}} = \vm{0}
\implies -\ln \phi_i - \ln Z - \vm{\lambda} \cdot (\vx_i - \vx) = 0 \quad 
(i = 1,2,\dots n),
\end{equation}
where $\lambda_0 = \ln Z - 1$ ($Z$ is the partition function).
On simplifying~\eref{eq:dL_maxent} and using~\eref{eq:reproducing_conditions_a}
leads us to the DtP map:
\begin{equation}\label{eq:dtp_maxent}
\phi_i(\vm{\lambda}) := \phi_i^H(\vm{\lambda}) =
 \frac{Z_i(\vm{\lambda})}{Z(\vm{\lambda})}, \quad
Z(\vm{\lambda}) = \sum_{j=1}^n Z_j(\vm{\lambda}), \quad
Z_i(\vm{\lambda}) = \exp \bigl( - \vm{\lambda} \cdot (\vx_i - \vx) \bigr) .
\end{equation}
On substituting $\phi_i^H(\vm{\lambda})$ from~\eref{eq:dtp_maxent}
in~\eref{eq:L_maxent} yields the dual function:
\begin{equation}
S(\vm{\lambda}) = \max_{ \vm{\phi} \in \Re_+^n } L (\vm{\phi},\lambda_0,\vm{\lambda})
                = {L}(\vm{\phi}^H(\vm{\lambda}),\lambda_0,\vm{\lambda})
                = \ln Z(\vm{\lambda}) ,
\end{equation}
and the dual variational principle can now be stated as:
\begin{equation}
\vm{\lambda}^\ast = \argmin \limits_{\vm{\lambda} \in \Re^2} \ S(\vm{\lambda}),
\end{equation}
which is an unconstrained convex optimization problem. Note that the dual
problem has two unknowns, whereas the primal problem has
$n \ge 3$ unknowns. Finally, we point the reader to a related primal-dual
connection in meshfree methods. 
Moving least squares basis functions are computed via
a dual formulation (minimizing a quadratic 
form~\cite{Belytschko:1996:MMO}), whose
primal formulation is the minimization of a quadratic form 
subject to the linear constraints 
in~\eref{eq:reproducing_conditions}~\cite{Black:1999:MAF}.

\section{Dual formulation for differential equations}\label{sec:dual}
We present the derivation of the dual 
functional and the dual variational principle for an
ODE and then for PDEs. First, as an introduction to the
duality method, we select an initial-value
problem. Then we consider the transient convection-diffusion equation 
in one spatial dimension, so that as special cases we 
can obtain the dual functional for the Laplace equation,
one-dimensional steady-state convection-diffusion equation,
transient heat equation, and the one-dimensional transport equation.

\subsection{Initial-value problem}\label{subsec:IVP}
We begin by illuminating 
the scheme in the simplest possible setting, which also explicitly demonstrates how an initial-value problem (IVP) 
can be robustly solved 
using a boundary-value problem in time. 
While simple, we are not aware of any existing variational principle for this problem. 
For $a,  u_0 \in \mathbb{R}$ 
\revised{and $T > 0$ denoting the terminal
time}, consider the following IVP:
\begin{subequations} \label{eq:ODE}
  \begin{align}
    \dot{u}(t) &= a u(t),  \ \  t \in (0,T), \label{eq:ODE_a} \\
    \intertext{with the initial condition}
     u(0) &= u_0. \label{eq:ODE_b}
    \end{align}
\end{subequations}
If $a < 0$ in~\eqref{eq:ODE_a}, the problem is dissipative.

The goal is to design a variational principle whose Euler--Lagrange equation is the ODE in an appropriate sense. 
To do so, we treat the 
ODE in~\eref{eq:ODE_a}
as a constraint and optimize the 
arbitrarily chosen auxiliary convex potential (density) function,
$H(u) = \frac{1}{2} u^2$,
by the method of Lagrange multipliers
with an important distinction -- instead of looking for a solution in the larger $(u,\lambda)$ space, we
treat the optimality condition with respect to the
primal field $u$ as generating a functional relationship for the primal field $u$ in terms of the dual
field $\lambda$ -- and in doing so, we 
solve the problem in the smaller space of just the dual variables (this idea is general, and applies to the considerations of Sec.~\ref{sec:formalism}). A crucial ingredient that enables this approach
is the free choice of the auxiliary potential. 
We begin by 
creating a functional in the fields $(u, \lambda)$ from \eqref{eq:ODE} with the dual field $\lambda$ treated like a 
test function and
the auxiliary function $H$ appended to it:
\begin{equation}\label{eq:L_IVP}
\begin{aligned}
   & \int_0^T [ H(u) + \lambda(\dot{u} - au) ] \, dt 
    = \int_0^T \left[ \dfrac{u^2}{2} + \lambda 
        ( \dot{u} - a  u ) \right] dt \\
    &= u(T) \lambda(T)  -  u(0) \lambda(0) + 
        \int_0^T \left[ \dfrac{u^2}{2}  
       -  u ( \dot{\lambda}   + a \lambda ) \right] dt \\
     &=  u(T) \lambda(T) - u_0 \lambda(0)
       + \int_0^T \left[ \dfrac{u^2}{2}  
       -  u ( \dot{\lambda}   + a \lambda ) \right] dt ,  
  \end{aligned}
\end{equation}
where integration by parts has been used and the 
initial condition has been incorporated to arrive at the
last equality. We now omit the term $u(T) 
\lambda(T)$ -- this 
choice is motivated by the sole design
consideration that the Euler--Lagrange equation of
the variational principle that we propose should 
recover the primal problem $\bigl($see discussion around \eqref{eq:uh_ivp}$\bigr)$
-- and refer to the functional that results as $L$: 
\begin{equation}\label{eq:Lhat_IVP}
L [u,\lambda] = \int_0^T {\cal L} (u,
{\cal D}) \, dt
- u_0 \lambda(0) , \quad
{\cal L}(u,{\cal D}) =  
\dfrac{u^2}{2} - u( \dot{\lambda} +
         a \lambda ) , \quad {\cal D} = (\lambda,\dot{\lambda}),
\end{equation}
where ${\cal L}$ is the Lagrangian density for the primal-dual problem.  

Next we require that 
\begin{equation}\label{eq:ODE_DtP}
    \dfrac{{\partial \cal L}}{\partial u}(u,{\cal D}) = 0
    \implies 
   u - \dot{\lambda} - a \lambda = 0 \mbox{ for all } {\cal D},
\end{equation}
which generates the dual-to-primal (DtP) mapping
${\cal D} \mapsto u_H({\cal D})$:
\begin{equation}\label{eq:uH}
u := u_H({\cal D}) = \dot{\lambda} + a \lambda .
\end{equation}
This solvability requirement serves as the main design specification for the choice of the auxiliary potential.
Substituting $u_H$ from~\eref{eq:uH} in~\eref{eq:Lhat_IVP}
generates the required dual functional (`action'):
\begin{equation}\label{eq:S_IVP}
S[\lambda] := L[u_H({\cal D}), \lambda] = 
- \dfrac{1}{2} \int_0^T \left[ u_H({\cal D})\right]^2 \, dt 
- u_0 \lambda(0) =  \int_0^T {\cal L}(u_H({\cal D}), {\cal D}) \, dt - u_0 \lambda(0),
\end{equation}
with the dual Lagrangian density given simply by ${\cal L}(u_H({\cal D}), {\cal D})$. To obtain the strong form posed in~\eref{eq:ODE}, we 
begin by setting $\delta S[\lambda; \delta \lambda] = 0$ for all 
variations, and use integration
by parts and the DtP map given in~\eref{eq:uH}
along with the condition that $\delta \lambda$ vanishes 
at the terminal time (defining the set of admissible variations) to obtain
\begin{equation}\label{eq:uh_ivp}
\dot{u}_H (t)= a u_H (t), \quad  u_H(0) = u_0
\end{equation}
(where for ease of notation we use
$u_H(t) := u_H\bigl({\cal D}(t)\bigr)$.
This last requirement $\bigl( \delta \lambda (T) = 0\bigr)$ indicates that one needs to impose a Dirichlet boundary condition on $\lambda(T)$ (not
necessarily vanishing), which was 
the sole motivation to
omit the term $u(T)\lambda(T)$ in~\eref{eq:L_IVP}. The key enabling feature of the above consistency of the dual functional with the primal problem is the condition $\frac{\partial {\cal L}}{\partial u} = 0$:
\begin{equation}\label{eq:ode_consist}
\begin{aligned}
    \delta S[\lambda; \delta \lambda] & = \int_0^T \frac{\partial {\cal L}}{\partial {\cal D}}(u_H, {\cal D}) \, \delta {\cal D} \, dt - u_0 \delta \lambda(0)\\
     & = - \int_0^T \left( u_H \delta \dot{\lambda} + 
      a u_H \delta \lambda  \right) \, dt - u_0 \delta \lambda(0) 
      = \int_0^T \left( \dot{u}_H - au_H \right) \delta \lambda \, dt + \big(u_H (0) - u_0 \big) \delta \lambda(0),\\
      \textrm{so that } \ \   0 & = \delta S[\lambda; \delta \lambda] \quad  \forall \delta \lambda, \ \delta \lambda(T) = 0  \ \ \Longleftrightarrow \eqref{eq:uh_ivp}, 
\end{aligned}
\end{equation}
on using the DtP mapping condition~\eqref{eq:ODE_DtP} and the affine dependence of ${\cal L}$ on ${\cal D}$ along with the 
Dirichlet boundary condition on $\lambda(T)$, 
and this persists regardless of the nonlinearities of ${\cal L}(u, {\cal D})$ in $u$. Substituting \eqref{eq:ODE_DtP} in \eqref{eq:ode_consist} along with using the arbitrariness of $\delta \lambda$ constrained by $\delta \lambda (T) = 0$  makes contact with the demonstration in~\ref{appendix:B},
where the weak form (the 
third expression in the chain of equalities 
in~\eqref{eq:ode_consist} set to zero)
is worked out explicitly in terms of dual variables.

On substituting the DtP map from~\eref{eq:uH} 
in~\eref{eq:uh_ivp}, the dual
Euler--Lagrange (second-order) boundary-value problem is given by:
\begin{subequations}\label{eq:dual_ODE}
    \allowdisplaybreaks
    \begin{align}
        \ddot{\lambda} - a^2 \lambda &= 0, \\
        \dot{\lambda}(0) + a \lambda(0) & = u_0,  \label{eq:dual_ODE_b} \\
        \lambda(T) &=  \lambda_T \quad \textrm{(arbitrarily chosen)}. \label{eq:final_bc}
    \end{align}
\end{subequations}
This is an \emph{elliptic} boundary-value-problem in  $(0,T)$ for the 
function $\lambda$. The reason why the condition at the
terminal time $t = T$ in~\eqref{eq:final_bc} can be imposed, even though the primal problem \eqref{eq:ODE} is an IVP that does not admit a `final-time' boundary condition, is that
\begin{equation}
u_H(T) = \dot{\lambda}(T) + a \lambda(T),
\end{equation}
and even with $\lambda(T)$ specified, 
$\dot{\lambda}(T)$ is free to adjust to recover the unique 
value of $u_H(T)$. Here,
$\dot{\lambda}(T) = u_0 e^{aT} - a \lambda_T$ 
has to be satisfied, which depends on both $u_0$ and $\lambda_T$.

Clearly, the solution $t \mapsto u_H(t) = u_H({\cal D}(t))$ cannot depend on 
$\lambda_T$ by uniqueness of the primal solution, even though
${\cal D}(t) = \bigl( \lambda(t), \dot{\lambda}(t) \bigr)$ depends on 
$\lambda_T$ in a non-trivial manner $\bigl($see~\eqref{eq:ode_1} and~\eqref{eq:ode_2} below$\bigr)$,
and we verify it in this specific case since it is perhaps
not obvious from visual inspection. 
On writing the dual solution as
\begin{equation}\label{eq:ode_1}
\lambda(t) = c_1 e^{|a|t} + c_2 e^{-|a|t}
\end{equation}
for constants $c_1, c_2$, the boundary conditions in~\eref{eq:dual_ODE} yield
\begin{equation}\label{eq:ode_2}
    \begin{bmatrix}
        (|a| + a)   &  -(|a| - a) \\ 
        e^{|a|T}  & e^{-|a|T}
    \end{bmatrix}
    \begin{Bmatrix}
        c_1 \\ 
        c_2
    \end{Bmatrix} = 
    \begin{Bmatrix}
        u_0 \\
        \lambda_T
    \end{Bmatrix},
\end{equation}
and on using~\eref{eq:uH}, $u_H(t)$ can be expressed as  
\begin{subequations}
\begin{align}
    u_H(t) &= 
    \begin{Bmatrix}
        (|a|+a) e^{|a|t}
        &
        -(|a| - a) e^{-|a|t}
    \end{Bmatrix}
    \raisebox{-0.94em}{$
    \begin{Bmatrix}
        c_1 \\
        c_2
    \end{Bmatrix}
    $} 
    \notag \\
    & = 
    \frac{1} {\Bigl( (|a|+a)e^{-|a|T} + (|a| - a) e^{|a|T}\Bigr)} 
    \begin{Bmatrix}
        (|a|+a)e^{|a|t} &
        - (|a| - a)e^{-|a|t}
    \end{Bmatrix}
    \raisebox{-0.94em}{$
    \begin{Bmatrix}
        u_0 e^{-|a|T} + \lambda_T (|a| - a) \\
        - u_0 e^{|a|T} + \lambda_T (|a| + a)
    \end{Bmatrix} 
    $} . \notag
\end{align}
\end{subequations}
The term involving $\lambda_T$ in $u_H(t)$ evaluates
to 
\begin{equation*}
\lambda_T (|a|^2 - a^2) e^{|a|t}   -  \lambda_T (|a|^2 - a^2) e^{-|a|t} = 0  \quad (!), 
\end{equation*}
and it can be verified that the rest of the terms evaluate to
\begin{equation*}
u_H(t) = u_0e^{at}.
\end{equation*}

In closing this section, we mention that the dual variational scheme applies seamlessly to nonlinear systems of ODEs, as demonstrated in \cite{Kouskiya:2024:HCH} for Euler's system for the motion of a rigid body about a fixed point with and without damping.

\subsection{Transient convection-diffusion equation}\label{subsec:transient_cd}
Consider the strong form of the
transient convection-diffusion model problem:
\begin{subequations}\label{eq:transient_cd}
\begin{align}
\kappa \dfrac{\partial^2 u}{\partial x^2}
 - \alpha  \dfrac{\partial u}{\partial x} &= \dfrac{\partial u}{\partial t}
 \ \ \textrm{in } \Omega = \Omega_0 \times \Omega_t = (0,1) \times (0,1), 
    \label{eq:transient_cd_a} \\
u(0,t) &= \bar{u}_1, \quad u(1,t) = \bar{u}_2, \label{eq:transient_cd_b} \\
u(x,0) &= u_0(x), \label{eq:transient_cd_c}
\end{align}
\end{subequations}
where $\kappa \geq 0$ is the diffusion coefficient, $\alpha$ is the convection
coefficient, 
$\bar{u}_1$ and $\bar{u}_2$ are boundary data (constants) and $u_0(x)$
is the prescribed initial condition.  In primal form, the convection-diffusion equation does not possess a variational structure, and therefore weak formulations have to be constructed using the strong form.  In the 
strongly convective regime, standard Galerkin methods produce spurious
oscillations. To remedy this deficiency, upwinding and/or stabilized methods
such as SUPG~\cite{Brooks:1982:SUP} and GLS~\cite{Hughes:1989:FEF} have been proposed, but they involve a user-specified tuning parameter and do not possess
a direct correspondence with a variational principle. 
Our goal is to
devise a variational scheme that is based on dual fields, so that standard 
Galerkin discretization can be directly used (without upwinding or stabilization, nor concerns about the inf-sup condition) to compute accurate solutions. \revised{Once the dual fields are solved, the
DtP map is used to compute the primal solutions.
The dual approach stands to provide significant advantages over existing methods to numerically solve problems that do not possess an exact, primal variational structure.}

On introducing the field $q = \partial u / \partial x$, we 
recast~\eref{eq:transient_cd_a} as a system of two first-order PDEs:
\begin{subequations}\label{eq:transient_cd_uq}
\begin{align}
\kappa \dfrac{\partial q}{\partial x}
- \alpha q &= \dfrac{\partial u}{\partial t}, \\
q &= \dfrac{\partial u}{\partial x}.
\end{align}
\end{subequations}
We refer to $(u,q)$ as the primal fields. Let {$\lambda$} and {$\mu$}
be the dual (Lagrange multiplier) fields that correspond to the two 
constraints in~\eref{eq:transient_cd_uq}. Now, we choose the convex 
potential (density) function $H(u,q)$ as:
\begin{equation}\label{eq:H_cd}
H(u,q) = \frac{1}{2}(u^2 + q^2) ,
\end{equation}
and consider the functional 
\begin{subequations}\label{eq:saddle_cd}
\begin{align}
   \widehat{L}[u,q,\lambda,\mu] = \int_0^1 \int_0^1 H(u,q) \, dx \, dt
   \ \ &+ \int_0^1 \int_0^1 
     \left[ \dfrac{\partial u}{\partial t} 
            - \kappa \dfrac{\partial q}{\partial x} 
            + \alpha q 
         \right] \lambda \, dx \, dt 
   - \int_0^1 \int_0^1 
     \left[ \dfrac{\partial u}{\partial x} - q \right] \mu \, dx \, dt .
     \label{eq:saddle_cd_b}
\end{align}
\end{subequations}
On using the divergence theorem on terms that involve the partial
derivatives of $u$ and $q$ in $\widehat{L}(u,q,\lambda,\mu)$, we obtain
\begin{align}\label{eq:Lag_cd_1}
\begin{split}
\widehat{L}[u,q,\lambda,\mu ] = \ \ & 
     \int_0^1 \int_0^1 \dfrac{u^2 + q^2}{2} dx \, dt  
     + \int_0^1 \int_0^1 
     \left[  \kappa q \dfrac{\partial \lambda}{\partial x} 
             + \alpha q \lambda - u \dfrac{\partial \lambda}{\partial t} 
             - q \mu - u \dfrac{\partial \mu}{\partial x} \right] dx \, dt
             \\
     &+ \int_0^1 \left[ u \mu - \kappa q \lambda \right]_{x=0}^{x=1}
     dt
     + \int_0^1 \left[ u \lambda \right]_{t=0}^{t=1} \, dx  .
\end{split}
\end{align}
There are six boundary terms in~\eref{eq:Lag_cd_1}: initial and
Dirichlet boundary conditions on $u$ are present in
three of them (these appear as natural boundary conditions in 
$\widehat{L}$),
and the remaining three terms involve $\lambda$.
On following the treatment of the boundary terms for the
IVP in~\sref{subsec:IVP}, we can
use any Dirichlet
data (zero or nonzero) on $\lambda$ for these boundary terms,
and the corresponding terms are dropped from $\widehat{L}$. 
We refer to
the functional that results as ${L}$. 
Now, on incorporating the boundary conditions 
from~\eref{eq:transient_cd_b}
and the initial condition from~\eref{eq:transient_cd_c}, 
we obtain 
\begin{align}\label{eq:Lag_cd_2}
\begin{split}
{L}[u,q,\lambda,\mu ] = \ \ & \int_0^1 \int_0^1 \dfrac{u^2 + q^2}{2}
                              dx \, dt 
     + \int_0^1 \int_0^1 
       \left[ - u \left( \dfrac{\partial \lambda}{\partial t} 
                         + \dfrac{\partial \mu}{\partial x} \right) 
              - q \left( - \kappa \dfrac{\partial \lambda}{\partial x} 
                         - \alpha \lambda
                         + \mu \right) \right] dx \, dt
             \\
     &+ \int_0^1 \left[ \bar{u}_2 \mu(1,t) - \bar{u}_1 \mu(0,t) \right] dt
      - \int_0^1  u_0(x) \lambda(x,0) \, dx ,
\end{split}
\end{align}
which remains valid for nonhomogeneous Dirichlet data for 
$\lambda$ (see the numerical example
in~\sref{subsec:laplace_solution}).
The functional ${L}$ in~\eqref{eq:Lag_cd_2} is an 
example of a 
pre-dual functional as introduced in prior works~\cite{Acharya:2022:VPN,Acharya:2023:DVP,Acharya:2024:VDS,Acharya:2024:APD,Acharya:2024:HCC, Singh:2024:HCN,Kouskiya:2024:HCH,Kouskiya:2024:IBD}.
We now pose the primal-dual/mixed variational 
problem as:
\begin{equation}\label{eq:Lhat}
\sup_{\lambda, \mu} \inf_{u, q} \ {L} [ u,q,\lambda,\mu ].
\end{equation}
Define
\begin{equation}\label{eq:FH_cd}
{\cal L}(u,q,{\cal D}) := \dfrac{u^2 + q^2}{2} 
                        - u \left( \dfrac{\partial \lambda}{\partial t} 
                                   + \dfrac{\partial \mu}{\partial x} \right) 
                        - q \left( - \kappa \dfrac{\partial \lambda}{\partial x} 
                                   - \alpha \lambda + \mu \right),
\quad {\cal D} = \left(\lambda,\frac{\partial \lambda}
{\partial x},\dot{\lambda}, \mu,
\frac{\partial \mu}{\partial x}\right),
\end{equation}
which is the space-time integrand that appears in~\eref{eq:Lag_cd_2}, and corresponds to the \emph{Lagrangian density} 
of variational field theories and particle mechanics.
On setting
\begin{subequations}\label{eq:dtp_cd}
\begin{align}
\dfrac{\partial {\cal L}}{\partial u} = 0, \quad  &
\dfrac{\partial {\cal L}}{\partial q} = 0, \label{eq:dtp_cd_a}
\intertext{we obtain the DtP map:}
u := u_H({\cal D}) = \dfrac{\partial \lambda}{\partial t} 
                        + \dfrac{\partial \mu}{\partial x}, \quad &
q := q_H({\cal D}) = \mu - \alpha \lambda 
     - \kappa \dfrac{\partial \lambda}{\partial x} .
\label{eq:dtp_cd_b}
\end{align}
\end{subequations}
Using~\eref{eq:dtp_cd_b} in~\eqref{eq:Lhat},
we have (in this quadratic setting)
\begin{equation}\label{eq:g_cd}
\inf_{u,q} {L}[u,q,\lambda,\mu ]  = {L} \bigl[ u_H({\cal D} ),
q_H({\cal D})\bigr]
 =: S [ \lambda,\mu ],
\end{equation}
where $S[ \lambda,\mu ] $ is the dual functional 
({\em action integral}), and now we 
state the dual variational principle:
\begin{subequations}\label{eq:dual_cd}
\begin{align}
\sup_{\lambda \in {\sf S}_\lambda, \ \mu \in {\sf S}_\mu} & S [ \lambda,\mu ] , \\
S [ \lambda,\mu ]  = \ \ - \dfrac{1}{2} \int_0^1 \int_0^1 
                  \left[ \bigl( u_H({\cal D}) \bigr)^2 + 
               \bigl( q_H({\cal D}) \bigr)^2 \right] dx\,dt   
      &+ \int_0^1 \left[ \bar{u}_2 \mu(1,t) - \bar{u}_1 \mu(0,t) \right] dt
      - \int_0^1  u_0(x) \lambda(x,0) \, dx , \\
  {\sf S}_\lambda = \{\lambda : \lambda \in H^1(\Omega), \ \lambda(0,t) = 
                       \lambda(1,t) &= \lambda(x,1) = 0 \}, \quad
  {\sf S}_\mu = \{\mu : \mu \in H^1(\Omega) \}  ,\label{eq:dual_cd_c}
\end{align}
\end{subequations}
where $H^1(\Omega)$ is the Sobolev space that contains
functions in $\Omega$ with square-integrable derivatives 
up to order $1$,
and $u_H({\cal D})$ and $q_H({\cal D})$ are given 
in~\eref{eq:dtp_cd_b}. In~\eref{eq:dual_cd_c},
we have chosen homogeneous Dirichlet boundary conditions
on $\lambda$:
$\lambda(0,t) = \lambda(1,t) = \lambda(x,1) = 0$.
To show that~\eref{eq:dual_cd} implies the
Euler--Lagrange equations of the primal problem
 and the imposed initial/boundary conditions (strong form)
given in~\eref{eq:transient_cd}, one proceeds along familiar
lines~\cite{Acharya:2024:HCC}: start with~\eref{eq:dual_cd} and
set $\delta S[ \lambda,\mu;\delta \lambda ] = 0$
and $\delta S[ \lambda,\mu;\delta \mu ] = 0$, apply the divergence theorem, and use~\eref{eq:dtp_cd_b} and the fundamental 
lemma of calculus of variations to arrive at the 
desired result.

The direct way to see the above consistency check is to realize that the dual functional, $S[\lambda, \mu]$, up to boundary terms, is the space-time integral of the Lagrangian ${\cal L}(u_H({\cal D}), q_H({\cal D}))$, 
and to 
note that $\frac{\partial {\cal L}}{\partial u} = 0$ and that ${\cal L}$ is necessarily affine in ${\cal D}$ with the coefficient of ${\cal D}$ formed from the primal strong form through integration by parts. These steps 
continue to hold
for nonlinear primal PDEs, and form the core idea of the duality-based method presented in this article.

\subsection{Laplace equation}\label{subsec:laplace}
Consider the Laplace equation in $\Omega = (0,1)$
with Dirichlet boundary conditions $u(0) = \bar{u}_1$ and
$u(1) = \bar{u}_2$. We choose the primal fields $u$ and $q$ that satisfy
\begin{equation}\label{eq:mixed_laplace}
q^\prime = 0, \qquad q = u^\prime,
\end{equation}
where $(\cdot)^\prime := d(\cdot)/dx$. On 
using~\eref{eq:dtp_cd_b} and~\eref{eq:dual_cd}, setting $\alpha = 0$ 
and $\kappa=1$, 
and dropping the time dependence, the DtP map and the
dual functional for this Laplace boundary-value problem are:
\begin{subequations}\label{eq:dual_laplace}
\begin{align}
u &:= u_H({\cal D}) = \mu^\prime, 
\quad q := q_H({\cal D}) = \mu -  \lambda^\prime,
\quad {\cal D} = (\lambda,\lambda^\prime,\mu,
\mu^\prime),
\label{eq:dual_laplace_a}\\
S [ \lambda,\mu ] &= -\dfrac{1}{2} \int_0^1 
                  \left[ \bigl( u_H({\cal D}) \bigr)^2 + 
                         \bigl( q_H({\cal D}) \bigr)^2 \right] dx
                   +  \bar{u}_2 \mu(1) - \bar{u}_1 \mu(0) , 
                   \ \ \lambda \in {\sf S}_\lambda, \ 
                   \mu \in {\sf S}_\mu, \label{eq:dual_laplace_b} \\ 
  {\sf S}_\lambda &= \{\lambda : \lambda \in H^1(0,1), \ \lambda(0) =
                       \lambda(1) = 0 \}, \quad
  {\sf S}_\mu = \{\mu : \mu \in H^1(0,1)\} . \label{eq:dual_laplace_c}
\end{align}
\end{subequations}

Let $\bar{u}_1 = 0$ and $\bar{u}_2 = 1$ be the Dirichlet data 
so that $\mu^\prime(0) = u(0) = 0$ and $\mu^\prime(1) = u(1) = 1$.
On substituting~\eref{eq:dual_laplace_a} in~\eref{eq:mixed_laplace} and
noting that the Dirichlet boundary data for $\lambda$ can be arbitrary, we find
that the dual fields satisfy
\begin{subequations}\label{eq:DUAL_laplace}
\begin{align}
\mu^\prime - \lambda^{\prime \prime} &= 0 \ \ \textrm{and} \ \ 
\mu - \lambda^\prime = \mu^{\prime \prime} \implies
\lambda^{\prime \prime \prime \prime} = 0 \ \ \textrm{and} \ \ 
\mu^{\prime \prime \prime} = 0,  \label{eq:DUAL_laplace_a} \\
\mu^\prime(0) &= 0, \ \ \mu^\prime(1) = 1, \ \ 
\lambda(0) = \lambda_0, \ \ \lambda(1) = \lambda_1, \label{eq:DUAL_laplace_b}
\end{align}
\end{subequations}
where $\lambda_0$ and $\lambda_1$ are arbitrary. Equation~\eref{eq:DUAL_laplace_a}
provides the Euler--Lagrange equations for the dual fields.
First, let us select $\lambda_0 = \lambda_1 = 0$, the choice we made earlier 
to define ${\sf S}_\lambda$ in~\eqref{eq:dual_laplace_c}.  Now, 
on using~\eqref{eq:DUAL_laplace_b}, the exact solution
for $\lambda$ (cubic function)
and $\mu$ (quadratic function) can be written as:
$\mu(x) = a_0 + x^2/2$ and 
$\lambda(x) = x(1-x)(b_0 + b_1 x)$, where $a_0$, $b_0$ and $b_1$ are constants.
Using the DtP map in~\eref{eq:dual_laplace_a},
the primal fields are:
\begin{equation*}
u(x) = x, \quad 
q(x) = (a_0 - b_0) + (2 b_0 - 2b_1 ) x + \left( \frac{1}{2} + 3b_1 \right)x^2 ,
\end{equation*}
and to recover the exact solution $q(x) = 1$ we must have
$b_1 = - 1/6$, $b_0 = b_1 = -1/6$ and $a_0 = 1 + b_0 = 5/6$.  Adding
any fixed cubic function to $\lambda(x)$ will not affect the recovery
of the exact solution, which is consistent with the fact that the boundary data
for $\lambda$ can be arbitrary (need not be homogeneous) in 
the set ${\sf S}_\lambda$ given in~\eref{eq:dual_laplace_c}.

\subsection{One-dimensional, steady state, convection-diffusion equation}\label{subsec:cd}
On choosing $\kappa = 1$ and dropping the time dependence 
in~\eref{eq:transient_cd}, we obtain the steady-state one-dimensional 
convection-diffusion problem: $u^{\prime \prime} - \alpha u^\prime = 0$, with 
boundary conditions $u(0) = \bar{u}_1$, $u(1) = \bar{u_2}$. On 
using~\eref{eq:dtp_cd_b} and~\eref{eq:dual_cd}, 
the DtP map and the dual functional are: 
\begin{subequations}\label{eq:dual_steadystate_cd}
\begin{align}
u &:= u_H({\cal D}) = \mu^\prime, \quad
q := q_H({\cal D}) = \mu - \alpha \lambda 
                        - \lambda^\prime, \quad 
{\cal D} = (\lambda,\lambda^\prime, \mu,
\mu^\prime),
\label{eq:dual_steadystate_cd_a}
\\
S [ \lambda,\mu ] &= -\dfrac{1}{2} \int_0^1 
                  \left[ \bigl( u_H({\cal D}) \bigr)^2 + 
                         \bigl( q_H({\cal D}) \bigr)^2 \right] dx
                   +  \bar{u}_2 \mu(1) - \bar{u}_1 \mu(0) , 
                   \ \ \lambda \in {\sf S}_\lambda, \ 
                   \mu \in {\sf S}_\mu, \label{eq:dual_steadystate_cd_b} \\ 
  {\sf S}_\lambda &= \{\lambda : \lambda \in H^1(0,1), \ \lambda(0) =
                       \lambda(1) = 0 \}, \quad
  {\sf S}_\mu = \{\mu : \mu \in H^1(0,1)\} . \label{eq:dual_steadystate_cd_c}
\end{align}
\end{subequations}

\subsection{Transient heat equation}\label{subsec:transient_heat}
Consider the initial-boundary value problem (IBVP)
 of heat conduction, where we choose
both Dirichlet and Neumann boundary conditions. Setting $\alpha = 0$
in~\eref{eq:transient_cd} and changing the second Dirichlet boundary
condition to a zero flux boundary condition, the IVBP reads:
\begin{subequations}\label{eq:transient_heat}
\begin{align}
\kappa \dfrac{\partial^2 u}{\partial x^2} &= \dfrac{\partial u}{\partial t}
 \ \ \textrm{in } \Omega = \Omega_0 \times \Omega_t = (0,1) \times (0,1), 
    \label{eq:transient_heat_a} \\
u(0,t) &= \bar{u}_1, \quad \kappa \dfrac{\partial u}{\partial x}(1,t) = 0, 
\label{eq:transient_heat_b} \\
u(x,0) &= u_0(x), \label{eq:transient_heat_c}
\end{align}
\end{subequations}
where $\kappa \geq 0$ is the thermal conductivity coefficient.

The dual formulation for this problem is derived in~\cite{Kouskiya:2024:HCH}.
We proceed by following the steps carried out to derive the dual functional 
in~\sref{subsec:transient_cd}. 
Referring to~\eref{eq:Lag_cd_1}, we now note that the boundary
term $\kappa q(1,t) \lambda(1,t)$ vanishes since
$\kappa q(1,t)$ is zero due to the Neumann boundary condition. Hence,
$\lambda(1,t)$ is unconstrained. In addition, since the boundary
term $u(1,t) \mu(1,t)$ is free, we choose to impose the
Dirichlet boundary condition
$\mu(1,t) = 0$. With these modifications in place, 
we use~\eref{eq:dtp_cd_b} and~\eref{eq:dual_cd} to obtain the
DtP map and the dual functional for the transient 
heat conduction problem:
\begin{subequations}\label{eq:dual_heat}
\begin{align}
u &:= u_H({\cal D}) = \dfrac{\partial \lambda}{\partial t} 
                        + \dfrac{\partial \mu}{\partial x}, \quad 
q := q_H({\cal D}) = \mu - \kappa \dfrac{\partial \lambda}{\partial x}  , \quad
{\cal D} = \left( \lambda, \frac{\partial \lambda}
{\partial x}, \dot{\lambda}, \mu,
\frac{\partial u}{\partial x} \right), 
\label{eq:dual_heat_a} \\
S [ \lambda,\mu ]  &= - \dfrac{1}{2} \int_0^1 \int_0^1 
                  \left[ \bigl( u_H({\cal D}) \bigr)^2 + 
                         \bigl( q_H({\cal D}) \bigr)^2 \right] dx\,dt
      - \int_0^1 \bar{u}_1 \mu(0,t) \, dt
      - \int_0^1  u_0(x) \lambda(x,0) \, dx , \label{eq:dual_heat_b} \\
  {\sf S}_\lambda &= \{\lambda : \lambda \in H^1(\Omega), \ \lambda(0,t) = 
                       \lambda(x,1) = 0 \}, \quad
  {\sf S}_\mu = \{\mu : \mu \in H^1(\Omega), \ \mu(1,t) = 0 \} .
  \label{eq:dual_heat_c}
\end{align}
\end{subequations}

\section{Variational formulation and discrete equations}\label{sec:varform_cd}
The numerical formulation and implementation of the dual variational form
for the transient convection-diffusion problem (see~\sref{subsec:transient_cd})
is presented.  Let $\vm{D} := (\lambda,\mu)$ denote the dual fields.
On setting the first variation of the dual functional
in~\eref{eq:dual_cd} to zero, we obtain:
\begin{subequations}\label{eq:dL=0_first}
\begin{align}
\delta S [ \lambda,\mu; \delta \lambda ] &=
                   - \int_0^1 \int_0^1 
                  \left[ u_H \delta u_H(\lambda,\mu;\delta \lambda) +
                         q_H \delta q_H(\lambda,\mu;\delta \lambda)
                  \right] dx \, dt
                  - \int_0^1  u_0(x) \delta \lambda(x,0) \, dx = 0
                   , \\
\delta S [ \lambda,\mu; \delta \mu ] &=
                   - \! \int_0^1 \int_0^1 
                  \left[ u_H \delta u_H(\lambda,\mu;\delta \mu) +
                         q_H \delta q_H(\lambda,\mu;\delta \mu)
                  \right] dx \, dt
                  + \int_0^1 \left[ \bar{u}_2 \delta \mu(1,t) 
                  - \bar{u}_1 \delta \mu(0,t) \right] dt = 0,
\end{align}
\end{subequations}
and on taking the variation of the DtP map in~\eref{eq:dtp_cd_b}
and substituting it in~\eref{eq:dL=0_first}, we obtain the variational equations
\begin{subequations}\label{eq:dL=0_second}
\begin{align}
                  - \int_0^1 \int_0^1 
                  \left\{ \
                  \left[
                   \dfrac{\partial \lambda}{\partial t} 
                    + \dfrac{\partial \mu}{\partial x} \right]
                  \dfrac{\partial (\delta \lambda)}{\partial t} 
                  + 
                  \left[ \mu -  \alpha \lambda - 
                  \kappa \dfrac{\partial \lambda}{\partial x} 
                         \right]
                          \left[ 
                          - \alpha \delta\lambda 
                          - \kappa \dfrac{\partial (\delta \lambda)}
                           {\partial x} 
                         \right] \
                  \right\} dx\,dt
                  - \int_0^1  u_0(x) \delta \lambda(x,0) \, dx &= 0
                  , \\
                  - \int_0^1 \int_0^1 
                  \left\{ \
                  \left[
                   \dfrac{\partial \lambda}{\partial t} 
                    + \dfrac{\partial \mu}{\partial x} \right]
                  \dfrac{\partial (\delta \mu)}{\partial x} 
                  + 
                  \left[ \mu - \alpha \lambda
                  - \kappa \dfrac{\partial \lambda}{\partial x} 
                         \right] \delta \mu
                  \right\} dx\,dt
                  + \int_0^1 \left[ \bar{u}_2 \delta \mu(1,t) 
                                    - \bar{u}_1 \delta \mu(0,t) \right] dt
                  &= 0,
\end{align}
\end{subequations}
which after rearranging leads to the statement of the dual variational form. 
Find $\vm{D} \in {\sf S}_\lambda \times {\sf S}_\mu$, 
such that
\begin{subequations}\label{eq:variational_form_cd}
\begin{align}
\begin{split}
& a_{11}(\lambda,\delta \lambda) + a_{12}(\mu,\delta \lambda) = \ell_1(\delta \lambda) \ \ \forall \delta \lambda \in {\sf S}_\lambda , \\
& a_{21}(\lambda, \delta \mu) + a_{22}(\mu,\delta \mu) = \ell_2(\delta \mu)
  \ \ \forall \delta \mu \in {\sf S}_\mu ,
\end{split}\\
\intertext{where the bilinear forms $a_{ij}(\cdot,\cdot)$, and the 
linear forms $\ell_1(\cdot$ and $\ell_2(\cdot)$ are given by}
\begin{split}
& a_{11}(\lambda,\delta \lambda) =  \int_0^1 \int_0^1 
                  \left[
                   \dfrac{\partial \lambda}{\partial t} 
                   \dfrac{\partial (\delta \lambda)}{\partial t} 
                   + \left( \alpha \lambda + \kappa 
                     \dfrac{\partial \lambda}{\partial x} \right) 
                   \left( \alpha \delta \lambda + \kappa 
                     \dfrac{\partial (\delta \lambda)}{\partial x} \right) 
                   \right] dx \, dt, \\
& a_{12}(\mu,\delta \lambda) =  \int_0^1 \int_0^1
                   \left[
                   \dfrac{\partial \mu}{\partial x} 
                   \dfrac{\partial (\delta \lambda)}{\partial t} 
                   - \mu 
                   \left( \alpha \delta \lambda + \kappa 
                     \dfrac{\partial (\delta \lambda)}{\partial x} \right) 
                   \right] dx \, dt, \\
& a_{21}(\lambda,\delta \mu) =  \int_0^1 \int_0^1
                  \left[
                   \dfrac{\partial \lambda}{\partial t} 
                   \dfrac{\partial (\delta \mu)}{\partial x} 
                   - \left(
                   \alpha \lambda + \kappa \dfrac{\partial \lambda}{\partial x} 
                   \right) \delta \mu \right] dx \, dt,\\
& a_{22}(\mu,\delta \mu) =  \int_0^1 \int_0^1
                   \left[
                   \dfrac{\partial \mu}{\partial x} 
                   \dfrac{\partial (\delta \mu)}{\partial x} 
                   + \mu \delta \mu \right]  dx \, dt ,
\end{split} \\
& \ell_1(\delta \lambda) = - \int_0^1  u_0(x) \delta \lambda(x,0)\, dx, \quad
  \ell_2(\delta \mu) = \int_0^1 \left[ \bar{u}_2 \delta \mu(1,t) 
                       - \bar{u}_1 \delta \mu(0,t) \right] dt, \quad
\intertext{and}
& {\sf S}_\lambda = \{\lambda : \lambda \in H^1(\Omega), \ \lambda(0,t) 
                      = \lambda(1,t) = \lambda(x,1) = 0 \}, \quad
  {\sf S}_\mu = \{\mu : \mu \in H^1(\Omega) \}.
  \label{eq:variational_form_cd_d}
\end{align} 
\end{subequations}
Homogeneous Dirichlet boundary conditions are chosen for
$\lambda$ in~\eref{eq:variational_form_cd_d}. However, we
remind the reader that 
$\lambda$ admits 
nonzero Dirichlet boundary conditions on the left, right and top
boundaries.

\subsection{Uniqueness of solutions for the dual variational system}\label{subsec:uniqueness}

To establish uniqueness of solutions to the dual variational equations, we
begin by considering two solutions $(\lambda_1, \mu_1)$ and $(\lambda_2, \mu_2)$ 
of~\eqref{eq:variational_form_cd}, and denote their difference as
\begin{equation*}
\lambda_1 - \lambda_2 =: \lambda_d, \qquad \qquad 
\mu_1 - \mu_2 =: \mu_d .
\end{equation*}
Clearly, $(\lambda_d, \mu_d)$ is a test function belonging to ${\sf S}_\lambda \times {\sf S}_\mu$, which also 
satisfies~\eqref{eq:variational_form_cd}:
\begin{equation}\label{eq:quad_form}
    \begin{aligned}
        a_{11} (\lambda_d,\lambda_d) + a_{12}(\mu_d,\lambda_d) + a_{21}(\lambda_d,\mu_d) + a_{22}(\mu_d,\mu_d) & = 0 \\
        \implies \int^1_0 \int_0^1
        \left[ \left( \frac{\partial \lambda_d}{\partial t} + \frac{\partial \mu_d}{\partial x} \right)^2 + \left( \mu_d - \kappa \frac{\partial \lambda_d}{\partial x} - \alpha \lambda_d\right)^2 \right]  dx \, dt & = 0.
    \end{aligned}
\end{equation}
Thus,
\begin{equation}\label{eq:uniq_1}
    \frac{\partial \lambda_d}{\partial t} + \frac{\partial \mu_d}{\partial x}  = 0 \ \ \textrm{a.e.\ in} \ \Omega 
\end{equation}
and
\begin{equation}\label{eq:uniq_2}
     \mu_d - \kappa \frac{\partial \lambda_d}{\partial x} - \alpha \lambda_d = 0 \ \ \textrm{a.e.\ in} \ \Omega.
\end{equation}
Noting the boundary conditions $\lambda_d(0,t) 
= \lambda_d(1,t) = 0$,
on multiplying~\eqref{eq:uniq_1} by $\lambda_d$ and integrating over $\Omega$ yields
\[
\int_0^1 \int_t^1 
\left[ \lambda_d \frac{\partial \lambda_d}{\partial t} - \mu_d\frac{\partial \lambda_d}{\partial x} 
\right] dx \, dt = 0, \qquad 0 \leq t \leq 1,
\]
and substituting for $\mu_d$ from \eqref{eq:uniq_2} one obtains
\begin{equation*}
    \begin{aligned}
    \int_0^1 \int_t^1 \left[ \frac{1}{2} \frac{\partial}{\partial t} \lambda_d^2 - \frac{\partial \lambda_d}{\partial x} \left( \kappa \frac{\partial \lambda_d} {\partial x} + \alpha  \lambda_d \right) \right] dx \, dt &= 0\\
       \implies
       \int_0^1 \int_t^1 \left[ \frac{1}{2} \frac{\partial}{\partial t} \lambda_d^2 - \kappa \left( \frac{\partial \lambda_d} {\partial x} \right)^2 - \alpha \frac{1}{2}\frac{\partial} {\partial x} \lambda_d^2 \right]
        dx \, dt & = 0\\
       \implies \int_t^1 \frac{\partial}{\partial t} \int_0^1 \frac{1}{2} \lambda_d^2 \, dx \, dt & = \int_0^1 \int_t^1 \kappa \left( \frac{\partial \lambda_d}{\partial x} \right)^2  dx \, dt \ge 0
    \end{aligned}
\end{equation*}
due to the boundary conditions
satisfied by $\lambda_d$ and $\kappa \geq 0$. We thus have, with the definition 
\[
A(t) := \int_0^1 \frac{1}{2} \lambda_d^2 (x,t) \, dx \geq 0,
\]
that
\[
A(1) - A(t)  \geq 0 \implies 0 = A(1) \geq A(t)  \geq 0 
\]
on using the fact that $\lambda_d(x,1) = 0$, so that $A(t) = 0$,
which implies that
\[
\lambda_d = 0 \ \ \textrm{a.e.\ in} \ \Omega, \\
\]
and from \eqref{eq:uniq_2} that
\[
\mu_d = 0 \ \ \textrm{a.e.\ in} \ \Omega,
\]
due to the arbitrariness of $t \in [0,1]$, and uniqueness follows.

\noindent \underline{Remark}: We note that for $\kappa = 0$, the primal system reduces to the linear transport equation, which is known to admit unique solutions that transport discontinuities in $u$ if present in initial or boundary data. It is interesting that the DtP mapping involves gradients of the dual fields to define the primal solution and therefore the dual second-order boundary-value problem in space-time cannot be elliptic, in order to represent such discontinuous primal solutions. Indeed, this is the case, as shown 
in~\cite[Secs.~2,\,3]{Kouskiya:2024:HCH}. Our uniqueness proof above is also a generalization of similar proofs, done separately for the heat and linear transport equations in~\cite{Kouskiya:2024:HCH},
where the degenerate ellipticity of the corresponding dual boundary-value problems is 
also shown. Since the highest order derivatives of the primal operator for the convection-diffusion equation considered here
correspond to the heat ($\kappa > 0$) and linear transport ($\kappa = 0$) 
problems, the dual problem herein
may be considered to be degenerate elliptic as well (and certainly for the special cases when either $\kappa$ or $\alpha$ vanish).
It is known that degenerate elliptic problems can have unique solutions~\cite{punzo2009uniqueness}.

\subsection{Discrete equations}\label{subsec:Kd=f}
We discretize the variational form using a standard Galerkin method.
On borrowing notation from neural networks (subscript $\theta$ points to the unknown weights and biases in a network), the numerical approximation for the dual fields are
written as:
\begin{equation}\label{eq:ansatz_cd}
\lambda_\theta(x,t) = \sum_{j=1}^{N_\lambda} \phi_j^\lambda (x,t) \, d_j
:= \vm{N}_\lambda(x,t) \, \vm{d}_\lambda \in {\cal S}_\lambda, \quad
\mu_\theta(x,t) = \sum_{j = 1}^{N_\mu} \phi_j^\mu (x,t) \, d_j^\mu
:= \vm{N}_\mu (x,t) \, \vm{d}_\mu \ \in {\cal S}_\mu ,
\end{equation}
where the sets $\{\phi_j^\lambda \}_{j=1}^{N_\lambda}$
and $\{\phi_j^\mu\}_{j=1}^{N_\mu}$ are approximating
functions for $\lambda$ and $\mu$, respectively. 
For convenience, we use $\lambda_\theta(x,t)$ and 
$\mu_\theta(x,t)$ to denote the trial functions in a 
neural network- or B-spline-based Galerkin method.
We ensure that by construction,
$\lambda_\theta(x,t)$ 
a priori satisfies the homogeneous Dirichlet boundary conditions 
in~\eref{eq:variational_form_cd_d}.  Let
$\vm{d} := \{ \vm{d}_\lambda \ \vm{d}_\mu \}^\top$ 
denote the vector that contains the unknown
coefficients. On selecting
trial functions of the form~\eref{eq:ansatz_cd} and
choosing admissible variations $\delta \lambda = \phi_i^\lambda$ and 
$\delta \mu = \phi_i^\mu$, and substituting them 
in~\eref{eq:variational_form_cd}, we obtain the following
system of linear equations:
\begin{subequations}\label{eq:Kd=f}
\begin{align}
\vm{K} \vm{d} &= \vm{f}, 
\quad \vm{K} = 
\begin{bmatrix} \vm{K}_{\lambda \lambda} & \vm{K}_{\lambda \mu} \\ 
                \vm{K}_{\mu \lambda} & \vm{K}_{\mu \mu} 
\end{bmatrix},
\quad 
\vm{f} = 
\begin{Bmatrix} \vm{f}_\lambda \\ \vm{f}_\mu \end{Bmatrix}, \label{eq:Kd=f_a}\\ 
\begin{split}\label{eq:Kd=f_b}
K_{\lambda \lambda}^{ij} &= \int_0^1 \int_0^1 
                  \left[
                   \dfrac{\partial \phi_j^\lambda}{\partial t} 
                   \dfrac{\partial \phi_i^\lambda}{\partial t} 
                   + \left( \alpha \phi_j^\lambda + \kappa 
                     \dfrac{\partial \phi_j^\lambda}{\partial x} \right) 
                   \left( \alpha \phi_i^\lambda + \kappa 
                     \dfrac{\partial \phi_i^\lambda}{\partial x} \right) 
                   \right] dx \, dt \Rightarrow
                   \vm{K}_{\lambda \lambda} = \int_0^1 \int_0^1 
                   \left[ 
                   \dfrac{\partial \vm{N}_\lambda^\top}{\partial t} 
                   \dfrac{\partial \vm{N}_\lambda}{\partial t} 
                   + \vm{C}_\lambda^\top \vm{C}_\lambda \right] 
                   dx \, dt , \\
K_{\lambda \mu}^{ij} &= \int_0^1 \int_0^1
                   \left[
                   \dfrac{\partial \phi_j^\mu}{\partial x} 
                   \dfrac{\partial \phi_i^\lambda}{\partial t} 
                   - \phi_j^\mu 
                   \left( \alpha \phi_i^\lambda + \kappa 
                     \dfrac{\partial \phi_i^\lambda}{\partial x} \right) 
                   \right] dx \, dt \Rightarrow
                   \vm{K}_{\lambda \mu} = \int_0^1 \int_0^1 
                   \left[
                   \dfrac{\partial \vm{N}_\lambda^\top}{\partial t} 
                   \dfrac{\partial \vm{N}_\mu}{\partial x} 
                   - \vm{C}_\lambda^\top \vm{N}_\mu \right] dx \, dt, \\
K_{\mu \lambda}^{ij} &= \int_0^1 \int_0^1
                   \left[
                   \dfrac{\partial \phi_j^\lambda}{\partial t} 
                   \dfrac{\partial \phi_i^\mu}{\partial x} 
                   - \left(
                   \alpha \phi_j^\lambda + \kappa 
                   \dfrac{\partial \phi_j^\lambda}{\partial x} 
                   \right) \phi_i^\mu \right] dx \, dt \Rightarrow
                   \vm{K}_{\mu \lambda} = \int_0^1 \int_0^1 
                   \left[
                   \dfrac{\partial \vm{N}_\mu^\top}{\partial x} 
                   \dfrac{\partial \vm{N}_\lambda}{\partial t} 
                   - \vm{N}_\mu^\top \vm{C}_\lambda 
                   \right] dx \, dt, \\
K_{\mu \mu}^{ij} &= \int_0^1 \int_0^1
                   \left[
                   \dfrac{\partial \phi_j^\mu}{\partial x} 
                   \dfrac{\partial \phi_i^\mu}{\partial x} 
                   + \phi_j^\mu \phi_i^\mu \right]  dx \, dt \Rightarrow
                   \vm{K}_{\mu \mu} = \int_0^1 \int_0^1 
                   \left[
                   \dfrac{\partial \vm{N}_\mu^\top}{\partial x} 
                   \dfrac{\partial \vm{N}_\mu}{\partial x} 
                   + \vm{N}_\mu^\top \vm{N}_\mu \right]  dx \, dt ,
\end{split} \\
\begin{split}\label{eq:Kd=f_c}
f_\lambda^i &= - \int_0^1  u_0(x) \phi_i^\lambda (x,0)\, dx \Rightarrow
\vm{f}_\lambda = - \int_0^1  u_0(x) \vm{N}_\lambda^\top (x,0)\, dx , \\
f_\mu^i &= \int_0^1 \left[ \bar{u}_2 \phi_i^\mu (1,t) 
                    - \bar{u}_1 \phi_i^\mu (0,t) \right] dt \Rightarrow
\vm{f}_\mu = \int_0^1 \left[ \bar{u}_2 \vm{N}_\mu^\top (1,t) 
                    - \bar{u}_1 \vm{N}_\mu^\top (0,t) \right] dt,
\end{split}
\intertext{where}
\vm{C}_\lambda &= \alpha \vm{N}_\lambda + 
                   \kappa \dfrac{\partial \vm{N}_\lambda}{\partial x}. \label{eq:Kd=f_d}
\end{align}
\end{subequations}
Observe that $\vm{K}$ is symmetric since 
$\vm{K}_{\mu \lambda} = \vm{K}_{\lambda \mu}^\top$. 
\revised{When}
the admissible trial functions for the dual fields in~\eref{eq:ansatz_cd} consist of linearly independent basis
functions,  it can be seen that the uniqueness proof in~\sref{subsec:uniqueness}
translates to uniqueness of solutions for the discrete
problem. Hence, $\vm{K}$ is invertible and 
the system of linear equations, $\vm{K}\vm{d} = \vm{f}$,
has a unique solution. It is interesting to observe that because of this invertibility of $\vm{K}$, the quadratic form, $\vm{d}^\top \vm{K} \vm{d}$, of the discrete problem obtained by the substitution $(\lambda_d, \mu_d) \to (\lambda_\theta, \mu_\theta)$ on the left-hand side
of~\eqref{eq:quad_form} is positive definite, even though the corresponding PDE is not elliptic, in general.

\section{Numerical examples}\label{sec:numerical_examples}
We apply the dual variational principle to solve the Laplace equation,
steady-state one-dimensional convection-diffusion equations, transient 
convection-diffusion
equation, and the transient heat equation. The transient problems are solved
as a space-time Galerkin method with a terminal boundary condition.
In previous applications of the dual formulation~\cite{Kouskiya:2024:HCH,Kouskiya:2024:IBD,Singh:2024:HCN},
the dual variables have been approximated using 
linear $C^0$ finite elements to solve linear and nonlinear
PDEs.  From~\eref{eq:dtp_cd_b}, since the primal
field $u$ depends on the space-time derivatives of $(\lambda, \mu)$, an $L^2$ projection
of the primal fields $u_H$ onto linear finite elements was
carried out to obtain a $C^0$-continuous solution. Herein,
we adopt smooth approximations for the dual fields so that they
directly provide at least $C^0$ primal fields. In one dimension, 
RePU activation function on a shallow network and univariate
B-splines are adopted as approximating functions,
whereas tensor-product B-splines are used for two-dimensional
(space-time) problems. \revised{The 
numerical implementation is carried out using symbolic computations 
in \texttt{Matlab}\texttrademark\ (R2024a Release).}

\subsection{Laplace equation}\label{subsec:laplace_solution}
Consider the Laplace equation, $u^{\prime \prime} = 0$ in 
$\Omega = (0,1)$, with boundary conditions $u(0) = 0, \ (1) = 1$. 
The exact solution is: $u(x) = x$. From~\eref{eq:DUAL_laplace_a}, 
we deduce that a quadratic approximation for $\mu(x)$ and a cubic
approximation for $\lambda(x)$ should recover the exact solution.  As a
verification test, we choose the following ansatz:
\begin{equation*}
\mu_\theta(x) = a_0 + a_1x + a_2x^2, \quad
\lambda_\theta(x) = x(1-x)(b_0 + b_1x), 
\end{equation*}
which are kinematically admissible since $\lambda(x)$ vanishes at $x = 0$ 
and $x = 1$.

On setting $\kappa = 1$, $\alpha = 0$, $\bar{u}_1 = 0$, $\bar{u}_2 = 1$, and
dropping the time dependence in~\eref{eq:Kd=f}, we obtain the
 expressions for the stiffness matrix and force vector.  The system of linear equations and its solution are:
\begin{align*}
\begin{bmatrix}
1 & \frac{1}{2} & \frac{1}{3}  & 0 & 0 \\
\frac{1}{2} & \frac{4}{3} & \frac{5}{4} & \frac{1}{6} & \frac{1}{12} \\
\frac{1}{3} & \frac{5}{4} & \frac{23}{15} & \frac{1}{6} & \frac{1}{10} \\
0 & \frac{1}{6} & \frac{1}{6} & \frac{1}{3} & \frac{1}{6} \\ 
0 & \frac{1}{12} & \frac{1}{10} & \frac{1}{6} & \frac{2}{15}
\end{bmatrix}
\begin{Bmatrix} a_0 \\ a_1 \\ a_2 \\ b_0 \\ b_1 \end{Bmatrix}
&= \begin{Bmatrix} 1 \\ 1 \\ 1 \\ 0 \\ 0 \end{Bmatrix}
\quad 
\Rightarrow
\quad
\begin{Bmatrix} a_0 \\ a_1 \\  a_2 \\ b_0 \\ b_1 \end{Bmatrix}
= 
\begin{Bmatrix}
\frac{5}{6} \\ 0 \\ \frac{1}{2} \\ -\frac{1}{6} \\ -\frac{1}{6}
\end{Bmatrix}.
\end{align*}
The solution for the dual fields are:
\begin{equation*}
\mu_\theta(x) = \frac{5}{6} + \frac{1}{2}x^2, \quad
\lambda_\theta(x) = -\frac{1}{6}x(1 - x^2) ,
\end{equation*}
and the solution for the primal fields are:
\begin{equation*}
u_\theta(x) = \mu_\theta^\prime(x) = x, \quad
q_\theta(x) = \mu_\theta (x) - \lambda_\theta^\prime(x) = 1,
\end{equation*}
and therefore the exact solution for $u$ and $q = u'$ are recovered.

To demonstrate that changing the boundary conditions on $\lambda(x)$ does not
affect the numerical solution, we add an arbitrary cubic function to
the earlier choice for $\lambda(x)$. Here, the ansatz for the dual fields are:
\begin{equation*}
\mu_\theta(x) = a_0 + a_1x + a_2x^2. \quad
\lambda_\theta(x) = x(1-x)(b_0 + b_1x) + 1 - 3x + x^2 + 3x^3.
\end{equation*}
The dual field $\lambda_\theta(x)$ satisfies the Dirichlet boundary conditions
$\lambda_\theta(0) = 1$ and $\lambda_\theta(1) = 2$. Now,
the system of linear equations and its solution are:
\begin{align*}
\begin{bmatrix}
1 & \frac{1}{2} & \frac{1}{3}  & 0 & 0 \\
\frac{1}{2} & \frac{4}{3} & \frac{5}{4} & \frac{1}{6} & \frac{1}{12} \\ 
\frac{1}{3} & \frac{5}{4} & \frac{23}{15} & \frac{1}{6} & \frac{1}{10} \\
0 & \frac{1}{6} & \frac{1}{6} & \frac{1}{3} & \frac{1}{6} \\ 
0 & \frac{1}{12} & \frac{1}{10} & \frac{1}{6} & \frac{2}{15}
\end{bmatrix}
\begin{Bmatrix} a_0 \\  a_1 \\  a_2  \\ b_0 \\ b_1 \end{Bmatrix}
&= \begin{Bmatrix} 2 \\ \frac{29}{12} \\ \frac{23}{10} \\ \frac{11}{6} 
                     \\ \frac{16}{15} \end{Bmatrix}
\quad \Rightarrow \quad
\begin{Bmatrix} a_0 \\ a_1 \\  a_2 \\ b_0 \\ b_1 \end{Bmatrix}
= 
\begin{Bmatrix}
\frac{11}{6} \\ 0 \\ \frac{1}{2} \\ \frac{23}{6} \\ \frac{17}{6}
\end{Bmatrix} .
\end{align*}
The solution for the dual fields are:
\begin{equation*}
\mu_\theta (x) = \frac{11}{6} + \frac{1}{2}x^2, \quad
\lambda_\theta (x) = \frac{x(1-x)(23 + 17x)}{6} + 1 - 3x + x^2 + 3x^3,
\end{equation*}
and the solution for the primal fields are:
\begin{equation*}
u_\theta(x) = \mu_\theta^\prime(x) = x, \quad
q_\theta(x) = \mu_\theta (x) - \lambda_\theta^\prime(x) = 1,
\end{equation*}
and hence once again the exact solutions for $u$ and $u^\prime$ are recovered.

\subsection{Steady-state one-dimensional convection-diffusion problem}
\label{subsec:cd_solution}
Consider the following one-dimensional convection-diffusion model problem~\cite{Sukumar:2022:EIB}:
\begin{subequations}
\begin{align}
u^{\prime \prime} - \alpha u^\prime &= 0 \ \ \textrm{in } \Omega = (0,1) \\
u(0) &= 0, \ \ u(1) = 1,
\intertext{with exact solution}
u(x) &= \frac{e^{\alpha x} - 1} {e^\alpha - 1},
\end{align}
\end{subequations}
where $\alpha$ is the Peclet number (ratio of convective rate to diffusive rate).  As the Peclet number increases, a boundary layer develops near the right end $x = 1$. 

\noindent \underline{Remark}: As noted earlier, in the strongly convective regime, standard Galerkin methods produce spurious oscillations~\cite{Hughes:1987:FEM}. On using 
the dual variational form, we can use a standard Galerkin method to compute
accurate solutions.

The dual functional for this problem is given
in~\eref{eq:dual_steadystate_cd}. On setting $\kappa = 0$
and dropping the time dependence in~\eref{eq:Kd=f}, we can write
the stiffness matrix (symmetric)
and the force vector as:
\begin{align}\label{eq:Kd=f_steadystate_cd}
\begin{split}
\vm{K} &= \begin{bmatrix}
    \vm{K}_{\lambda \lambda} & \vm{K}_{\lambda \mu} \\
    \vm{K}_{\mu \lambda } & \vm{K}_{\mu \mu}
    \end{bmatrix}, \quad
    \vm{f} = \begin{Bmatrix}
    \vm{f}_{\lambda} \\ \vm{f}_{\mu}
    \end{Bmatrix}, \\
\vm{K}_{\lambda \lambda} &= \int_0^1 \alpha^2
                   \vm{N}_\lambda^\top \vm{N}_\lambda \, dx, \quad
                   \vm{K}_{\lambda \mu} = - \int_0^1 
                    \alpha \vm{N}_\lambda^\top \vm{N}_\mu \,dx , \\
                   \vm{K}_{\mu \lambda} 
                   &= - \int_0^1 
                    \alpha \vm{N}_\mu^\top \vm{N}_\lambda \, dx ,
                    \quad
                   \vm{K}_{\mu \mu} = \int_0^1 
                   \left[
                   \dfrac{\partial \vm{N}_\mu^\top}{\partial x} 
                   \dfrac{\partial \vm{N}_\mu}{\partial x} 
                   + \vm{N}_\mu^\top \vm{N}_\mu \right]  dx,  \\
\vm{f}_\lambda &= \vm{0}, \quad \vm{f}_\mu = \vm{N}_\mu^\top (1) ,
\end{split}
\end{align}
where $\vm{N}_\lambda$ and $\vm{N}_\mu$ define row vectors that
contain the approximating functions that are 
used to form $\lambda_\theta(x)$
and $\mu_\theta(x)$, respectively. 
We present numerical solutions using a shallow neural network with RePU activation function and univariate B-splines.

\subsubsection{Neural network solutions}
The RePU ($\textrm{ReLU}^k$) activation function~\cite{He:2024:EAP} in
a neural network is given by
\begin{equation}\label{eq:repu}
\sigma(x;p) = \texttt{RePU}(x;p) = \texttt{ReLU}^p(x) = 
[\max(0,x)]^p.
\end{equation}
It is known
that deep neural networks with ReLU and $\textrm{ReLU}^2$ activation
functions can represent Lagrange finite elements of any order 
in arbitrary dimensions~\cite{He:2024:DNN}. 
Here we consider a shallow neural network (single hidden layer) with $2n$
neurons. The interval $\Omega := [0,1]$ is discretized with the knot sequence:
\begin{equation}\label{eq:knot}
\Xi := [x_0, \,  x_1, \dots, \, x_n],
\end{equation}
where $x_0 = 0$ and $x_n$ = 1.  The bias and weight are fixed in
the hidden layer.  A linear output is used to define the two dual fields $\lambda_\theta(x)$ and $\mu_\theta(x)$. 

Smooth approximations are chosen for the 
dual fields $\mu_\theta(x)$ (degree $p$) and
$\lambda_\theta(x)$ (degree $q$), which are written in the form:
\begin{subequations}\label{eq:repu_trial_cd}
\begin{align}
\mu_\theta(x) &= \sum_{i=0}^{n-1} \sigma(x-x_i;p) a_i^{+}
                 + \sum_{i=1}^{n} \sigma(x_i-x;p) a_i^{-}
              := \sum_{k=1}^{2n} \phi_k^\mu(x) \, \vm{a}_k
                 \ \in C^{p-1}(\Omega),\label{eq:repu_trial_cd_a} \\
\begin{split} \label{eq:repu_trial_cd_b}
\lambda_\theta(x) &= (1-x)\lambda_\theta^+(x) + x \lambda_\theta^-(x) \\
                  &=
                  (1-x)
                  \left[ \sum_{i=0}^{n-1} \sigma(x-x_i;q) \,
                  b_i^{+} \right] + x
                  \left[ \sum_{i=1}^{n} \sigma(x_i-x;q) \, b_i^{-} \right]
                  := \sum_{k=1}^{2n} \phi_k^\lambda(x) \, \vm{b}_k
                  \ \in C^{q-1}(\Omega) . % \ \ (q \geq 2) .
\end{split}
\end{align}
\end{subequations}
The shallow neural network approximations that appear in~\eref{eq:repu_trial_cd} are atypical in the neural network literature. Instead of nonlinear approximations with both
bias and weights that are unknowns, we fix the
bias and weight ($z = wx + b$, $w = \pm 1$, $b = \mp x_i$)
in the hidden layer so that the approximations are linear in the 
unknowns (weights) from the output layer.
\revised{The two RePU functions with arguments $x - x_i$ 
and $x_i - x$ are
associated with all interior knot locations
$\{x_i\}_{i=1}^{n-1}$.  In
$\Omega = [0,1]$, the RePU functions associated with
$i = 0$ and $i = n$ are $x^p$ and
$(1 - x)^p$, respectively.
Consider the approximation $\mu_\theta(x)$ 
(with $p = 1$) given in~\eqref{eq:repu_trial_cd_a}. For $n = 2$, the approximating functions in $\Omega$ are:
$x$, $1-x$, 
$\textrm{ReLU}(x-\frac{1}{2})$ and 
$\textrm{ReLU}(\frac{1}{2}-x)$. Equivalently, these
functions can be considered to be elements of the set 
${\cal V}_\theta := \left\{ x, \, 1-x, \, x - \frac{1}{2}, \, |x-\frac{1}{2}| \, \right\}$.
If ${\cal V} = \{1, \, x\}$ is a basis for the
vector space of polynomials of degree 1, then 
the elements of ${\cal V}_\theta$ span 
${\cal V}$ but constitute a frame (and not a basis),
since ${\cal V}_\theta$ contains three affine functions.  
Similar steps and arguments
for $p = 2$ (with $n = 3$) lead us to the result that the
$\textrm{RePU}$ functions in~\eqref{eq:repu_trial_cd_a} are linearly dependent. By extension, we deduce that in general the trial functions in~\eqref{eq:repu_trial_cd} consist of linearly dependent functions.   However, contrary to this being viewed as a deficiency, it is in fact a 
strength of the dual approach that even linearly 
dependent functions can be chosen in the ansatz for the dual fields.  Based on 
the discussion surrounding~\eqref{eq:lambdadual_Ax} and~\eqref{eq:xdual_Ax}, any feasible solution for the dual linear system (using any generalized inverse) when used in the DtP map yields a numerically computed primal solution. 
This relies on the fact that the linear system is consistent: right-hand side vector $\vm{f}$ must belong
to the column space of the stiffness matrix.  As noted earlier, ${\cal V}_\theta$ contains linearly dependent functions that span ${\cal V}$. If we remove these additional functions in ${\cal V}_\theta$, then 
the resulting set will be a basis that spans degree 1 polynomials. Hence,
$\vm{K}$ will be invertible and the linear system will have a unique solution. So we can deduce that the linear system is consistent.
This is
confirmed by the numerical solutions that are presented 
in~\fref{fig:alpha1_p2q3_cd} 
to~\fref{fig:alpha50_p2q3_cd}.}
Observe that $\lambda_\theta(x)$ in~\eref{eq:repu_trial_cd_b} satisfies the 
Dirichlet boundary conditions since $\lambda_\theta(0) = \lambda_\theta(1) = 0$.
Plots of RePU functions that are used to construct 
$\mu_\theta(x)$ ($p = 2$) and $\lambda_\theta(x)$ ($q = 3$) 
are presented in~\fref{fig:repu_mulambda}.
These functions are polynomials on either side of the knot location $x_i$. The functions
$\textrm{ReLU}^2(x - x_i)$ and $\textrm{ReLU}^2(x_i - x)$ in Figs.~\ref{fig:repu_mulambda_a} and~\ref{fig:repu_mulambda_b}, respectively, contribute to $\mu_\theta(x)$. A  barycentric convex
combination of the contributions due to
$\textrm{ReLU}^3(x - x_i)$ and $\textrm{ReLU}^3(x_i - x)$ in Figs.~\ref{fig:repu_mulambda_c} and~\ref{fig:repu_mulambda_d}, respectively, are used to form $\lambda_\theta(x)$. 
\begin{figure}[!tbh]
\centering
\begin{subfigure}{0.48\textwidth}
\includegraphics[width=\textwidth]{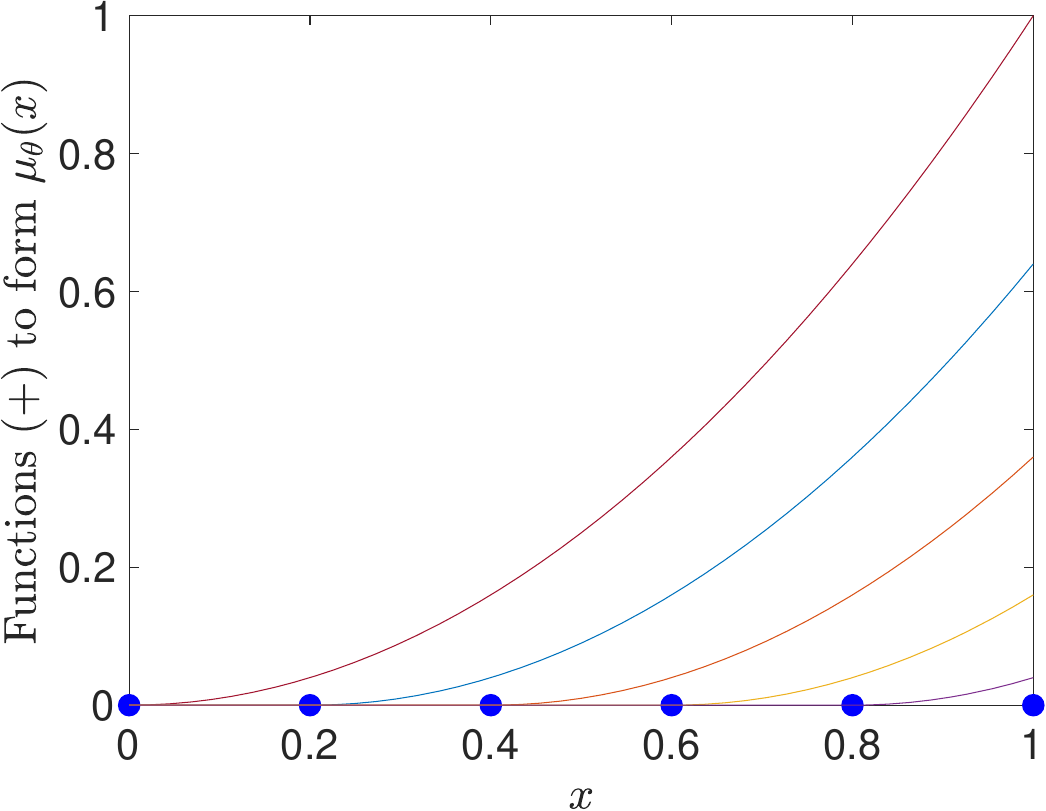}
\subcaption{}\label{fig:repu_mulambda_a}
\end{subfigure} \hfill
\begin{subfigure}{0.48\textwidth}
\includegraphics[width=\textwidth]{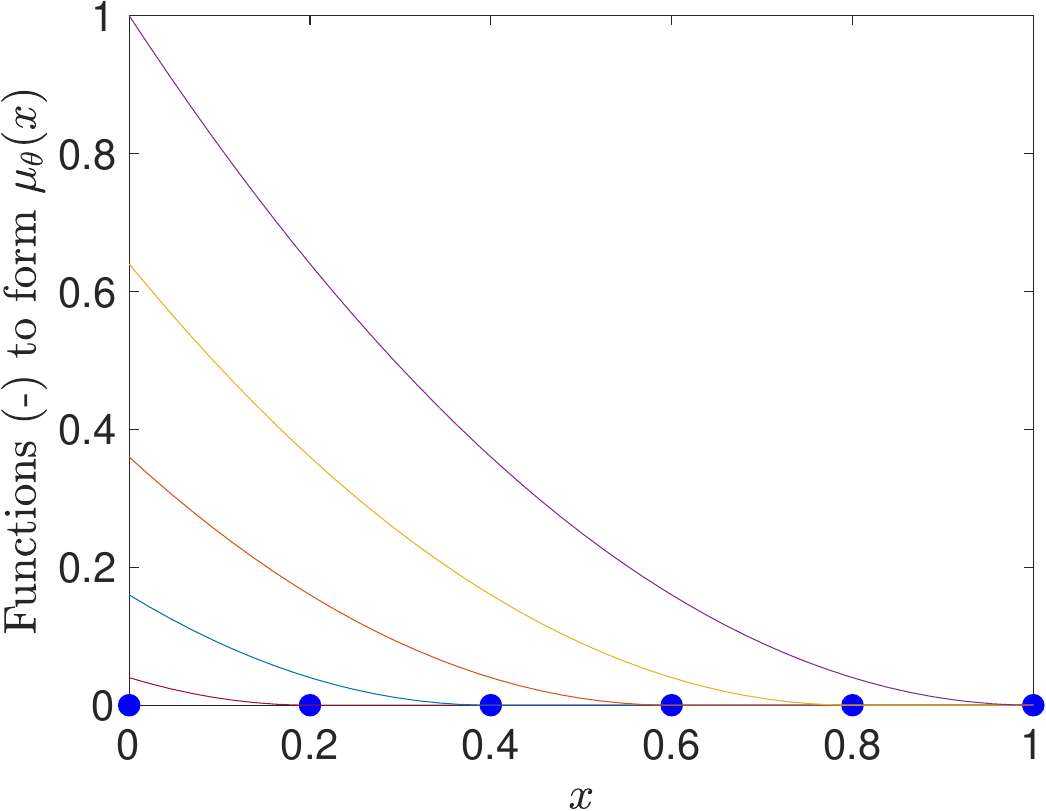}
\subcaption{}\label{fig:repu_mulambda_b}
\end{subfigure} \hfill
\begin{subfigure}{0.48\textwidth}
\includegraphics[width=\textwidth]{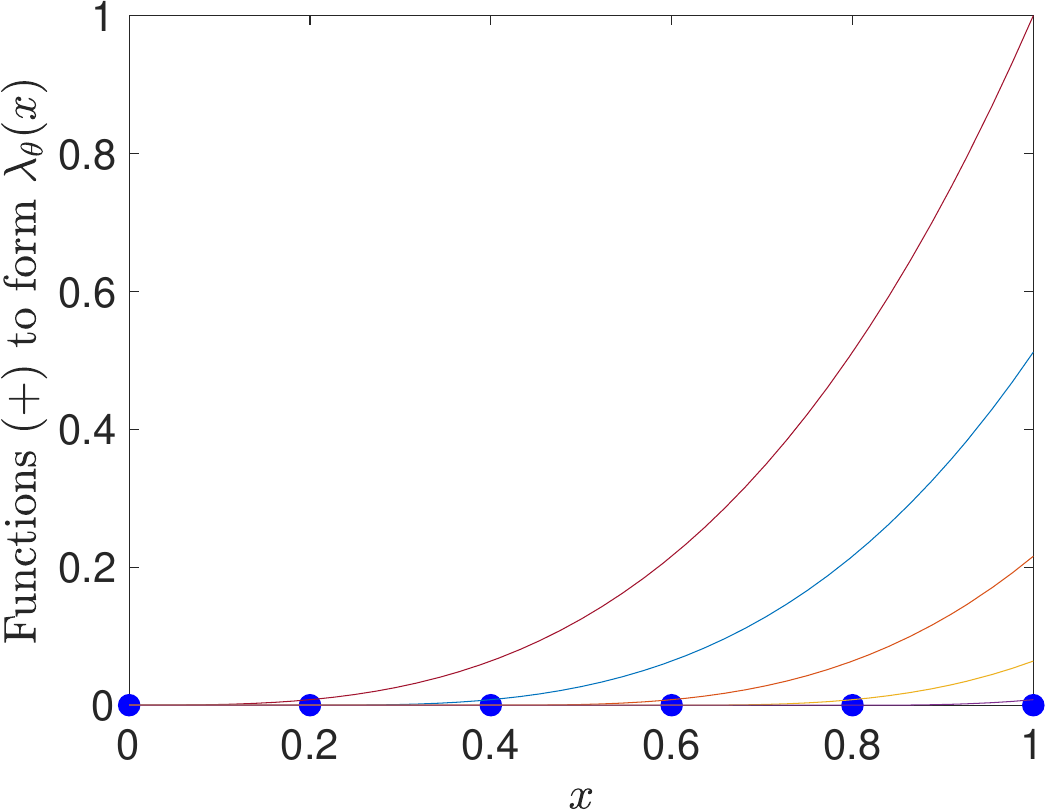}
\subcaption{}\label{fig:repu_mulambda_c}
\end{subfigure} \hfill
\begin{subfigure}{0.48\textwidth}
\includegraphics[width=\textwidth]{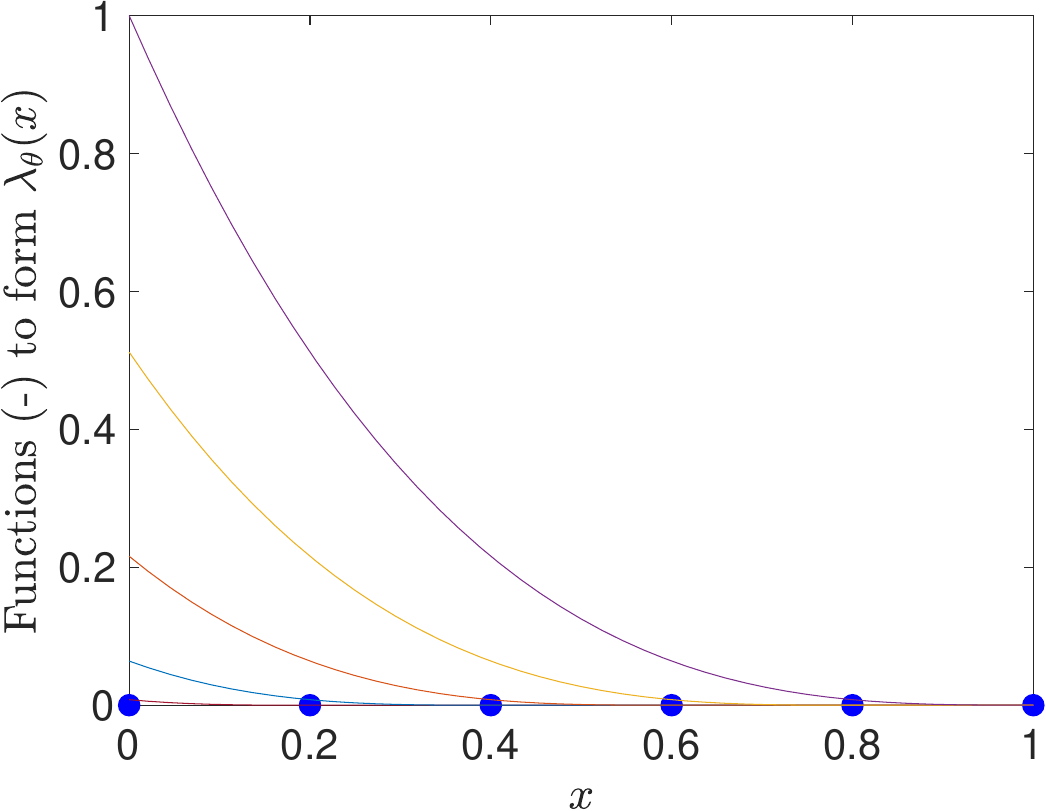}
\subcaption{}\label{fig:repu_mulambda_d}
\end{subfigure}
\caption{Plots of RePU functions that are used to form the dual fields. 
         Filled circles in blue denote the location of the knots,
         $x_i$ $(i = 0,1,\dots,5)$. Knots are uniformly spaced.
         Linear combination of the functions
         (a) $\sigma(x-x_i;p=2)$ and (b) $\sigma(x_i-x;p=2)$ are used
         to form $\mu_\theta(x)$.
         Linear combination of the functions (c) $\sigma(x-x_i;q=3)$ 
         are used to form $\lambda_\theta^+(x)$ and linear
         combination of the functions (d) $\sigma(x_i-x;q=3)$
         are used to form
         $\lambda_\theta^-(x)$. The dual field $\lambda_\theta(x) =
         (1-x) \lambda_\theta^+(x) + x \lambda_\theta^-(x)$.}
         \label{fig:repu_mulambda}
\end{figure}

% p = 2, q = 3 (alpha = 1)

For $\alpha = 1$, 
the neural network solution is compared to the exact
solution in~\fref{fig:alpha1_p2q3_cd}.  
In the neural network computations,
eleven distinct
knot locations are used ($n = 10$). Degrees $p = 2$ are $q = 3$ are
chosen to form $\mu_\theta(x)$ and $\lambda_\theta(x)$, 
respectively. The dual fields are presented
in~\fref{fig:alpha1_p2q3_cd_b} and the primal fields,
which are computed from the dual fields using the DtP
map in~\eqref{eq:dual_steadystate_cd_a}, are
shown in~\fref{fig:alpha1_p2q3_cd_c}.
We observe
that the maximum errors in Figs.~\ref{fig:alpha1_p2q3_cd_d} 
and~\ref{fig:alpha1_p2q3_cd_e} for $u$ and $u^\prime$
are about $10^{-3}$ and $10^{-4}$, respectively. 
\begin{figure}
\centering
\begin{subfigure}{0.33\textwidth}
\includegraphics[width=\textwidth]{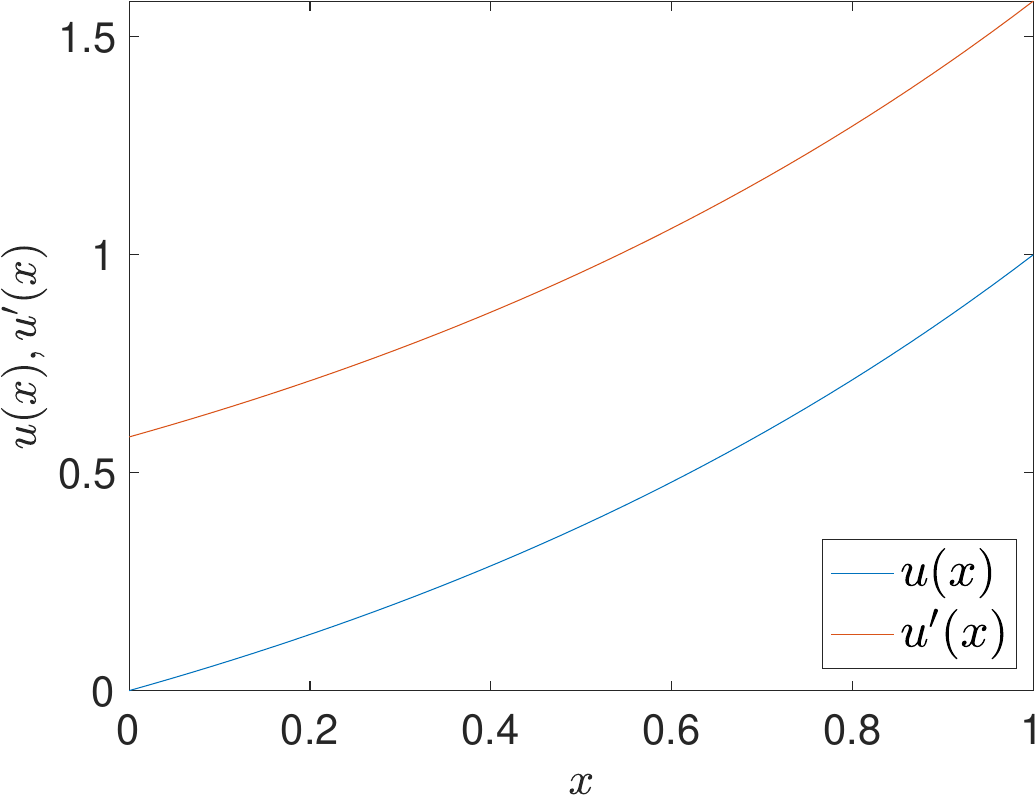}
\subcaption{}\label{fig:alpha1_p2q3_cd_a}
\end{subfigure}
\begin{subfigure}{0.33\textwidth}
\includegraphics[width=\textwidth]
{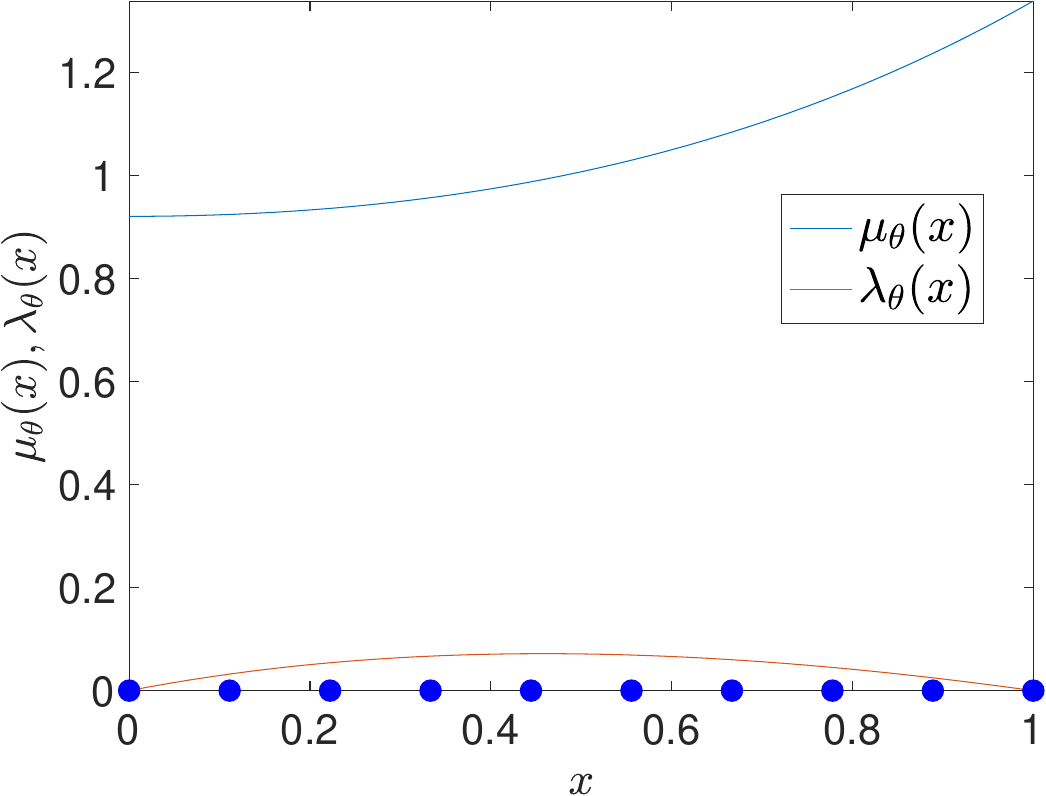}
\subcaption{}\label{fig:alpha1_p2q3_cd_b}
\end{subfigure} \hfill
\begin{subfigure}{0.33\textwidth}
\includegraphics[width=\textwidth]{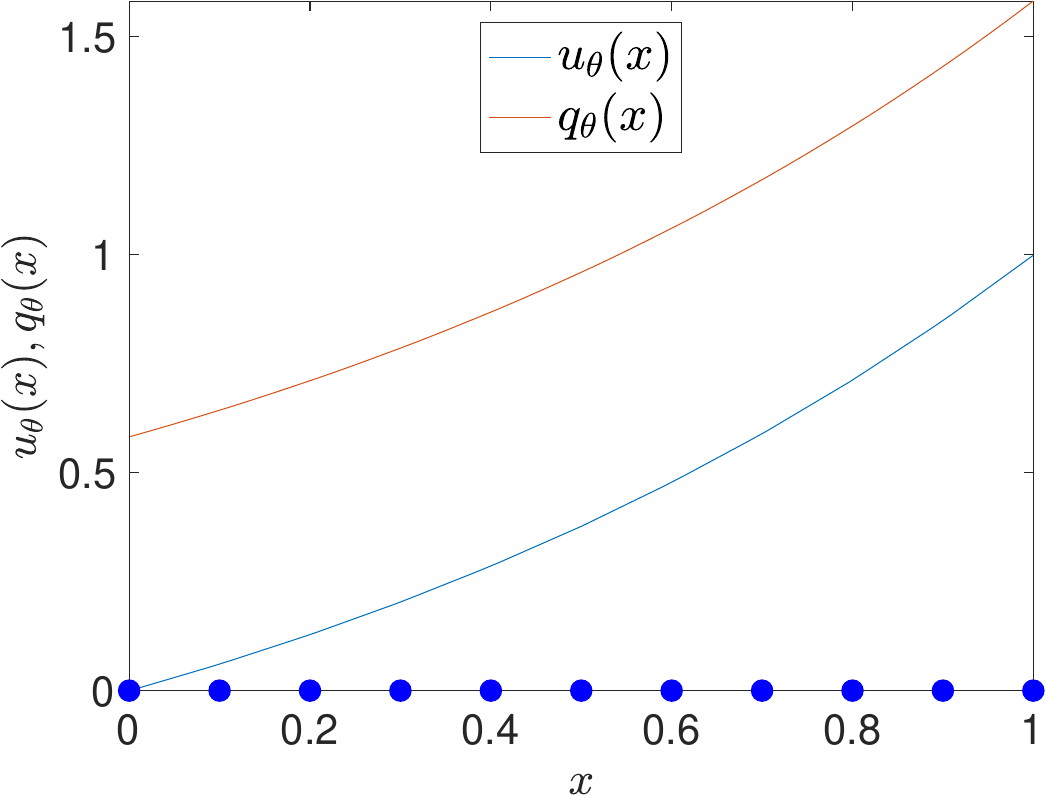}
\subcaption{}\label{fig:alpha1_p2q3_cd_c}
\end{subfigure}
\begin{subfigure}{0.48\textwidth}
\includegraphics[width=\textwidth]{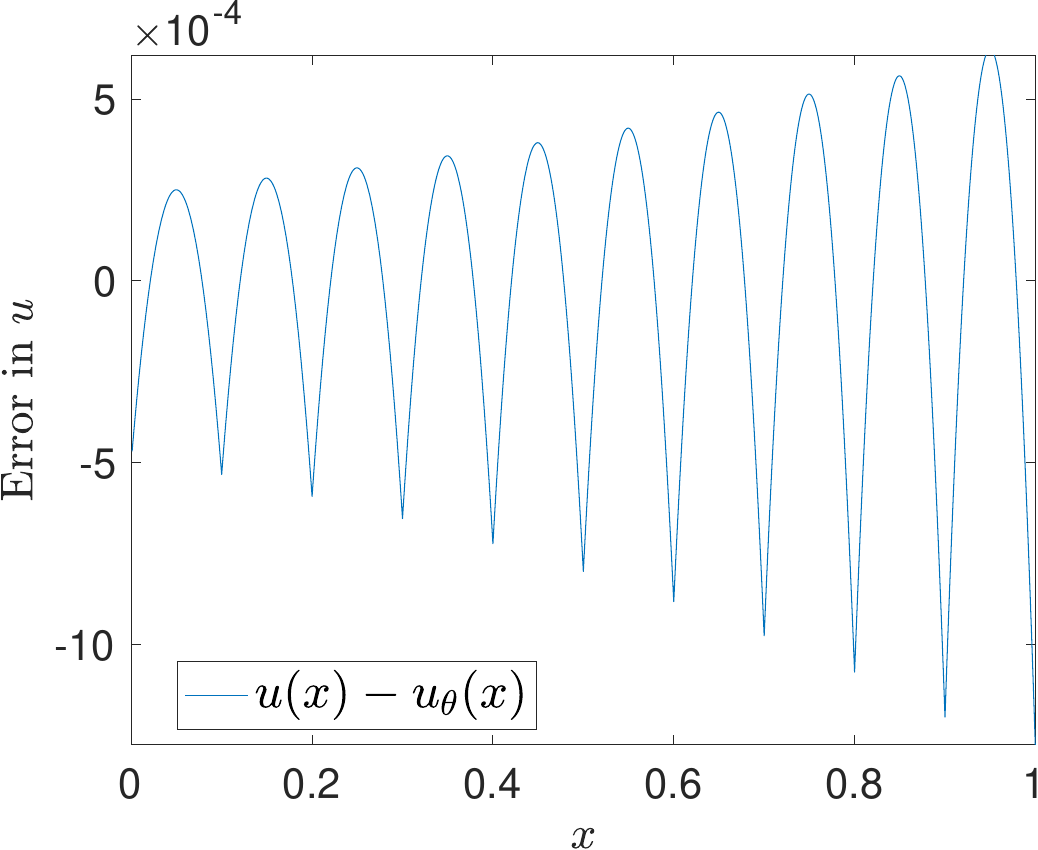}
\subcaption{}\label{fig:alpha1_p2q3_cd_d}
\end{subfigure}
\begin{subfigure}{0.48\textwidth}
\includegraphics[width=\textwidth]{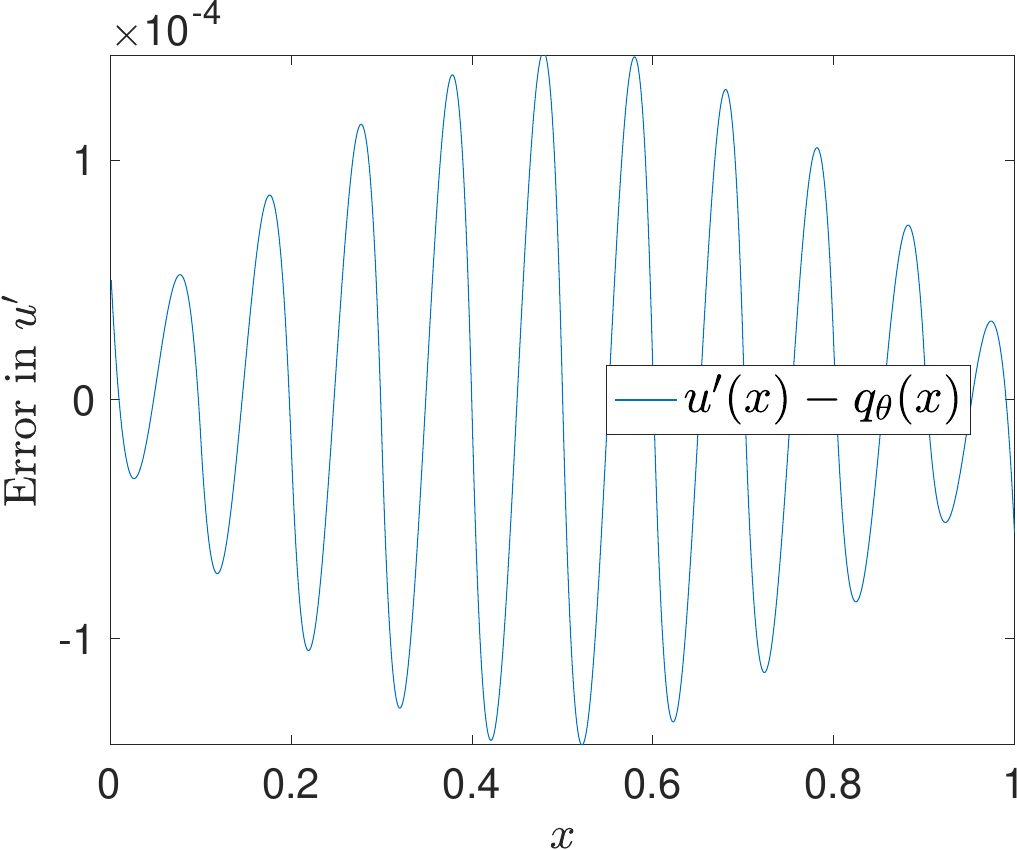}
\subcaption{}\label{fig:alpha1_p2q3_cd_e}
\end{subfigure}
\caption{Neural network solution for the steady-state 
         convection-diffusion problem ($\alpha = 1$) using
         $n = 10$, $p = 2$ and $q = 3$. 
         (a) Exact solutions for $u$ and $u^\prime$;
         (b) Dual fields, $\mu_\theta(x)$
             and $\lambda_\theta(x)$;
         (c) Primal fields, 
             $u_\theta(x)$ and $q_\theta (x)$;
         (d) Error in $u$; and (e) Error in $u^\prime$.}
         \label{fig:alpha1_p2q3_cd}
\end{figure}
%
% p = 3, q = 4 (alpha = 1)
Similar trends are observed in the results presented 
in~\fref{fig:alpha1_p3q4_cd} for
$\mu_\theta(x)$ ($p = 3$) and
$\lambda_\theta(x)$ ($q = 4$).  The maximum errors in 
Figs.~\ref{fig:alpha1_p3q4_cd_b} and~\ref{fig:alpha1_p3q4_cd_c}
for $u$ and $u^\prime$ are $10^{-5}$
and $4 \times 10^{-5}$, respectively.  In the results shown in 
both Figs.~\ref{fig:alpha1_p2q3_cd} and~\ref{fig:alpha1_p3q4_cd}, the accuracy
of $q_\theta = u_\theta^\prime(x)$ is seen to be comparable to that of $u_\theta(x)$.  This is so since the degree of the piecewise 
polynomial $\lambda_\theta(x)$ (contributes to 
$q_\theta = u_\theta^\prime$, which is the constraint equation
that is weakly enforced)
is one greater than that of $\mu_\theta(x)$ 
($u_\theta = \mu_\theta^\prime$).
\begin{figure}[!tbh]
\centering
\begin{subfigure}{0.33\textwidth}
\includegraphics[width=\textwidth]{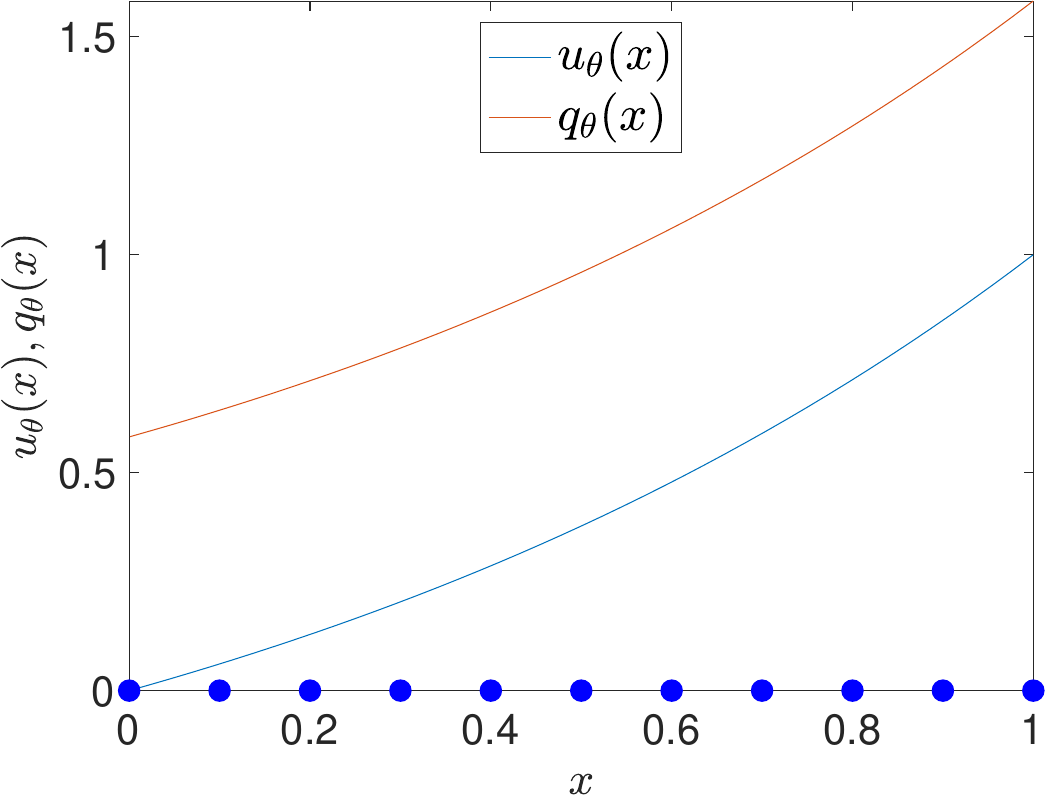}
\subcaption{}\label{fig:alpha1_p3q4_cd_a}
\end{subfigure} \hfill
\begin{subfigure}{0.33\textwidth}
\includegraphics[width=\textwidth]{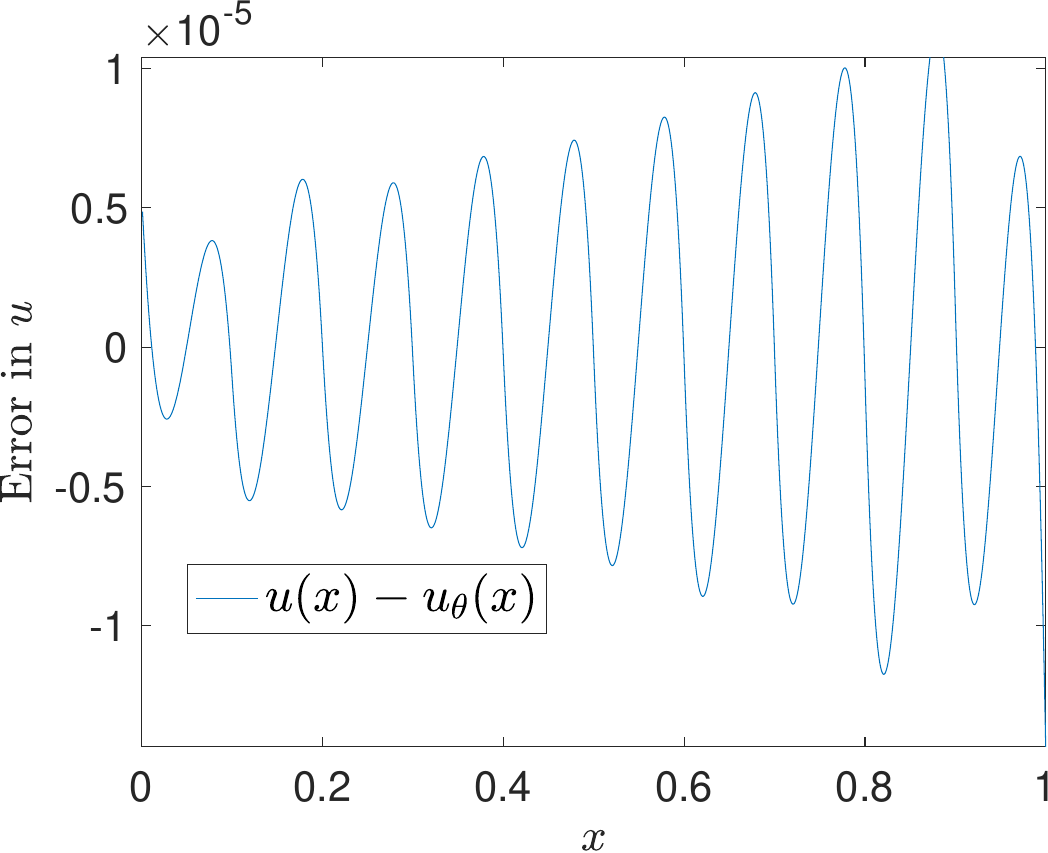}
\subcaption{}\label{fig:alpha1_p3q4_cd_b}
\end{subfigure} \hfill
\begin{subfigure}{0.33\textwidth}
\includegraphics[width=\textwidth]{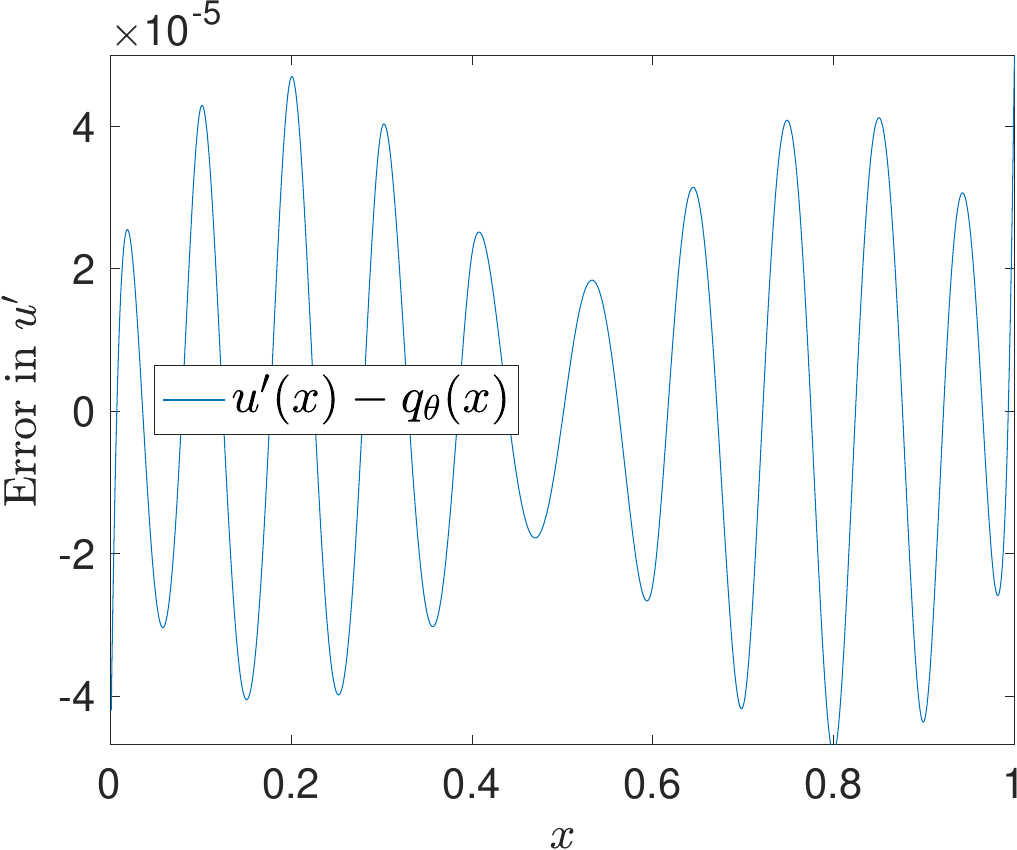}
\subcaption{}\label{fig:alpha1_p3q4_cd_c}
\end{subfigure}
\caption{Neural network solution for the steady-state convection-diffusion 
         problem ($\alpha = 1$) using $n = 10$, $p = 3$, and $q = 4$. 
         (a) $u_\theta(x), \ q_\theta(x)$;
         (b) $u - u_\theta$; and
         (c) $u^\prime - q_\theta$.}
         \label{fig:alpha1_p3q4_cd}
\end{figure}
% p = 2, q = 3 (alpha = 10)
Next we set $\alpha = 10$ and choose $n = 30$ in the numerical discretization. The neural network solution is compared to the exact
solution in~\fref{fig:alpha10_p2q3_cd} for $\mu_\theta(x)$ ($p = 2$)
and $\lambda_\theta(x)$ ($q = 3$). The exact solution
displays a sharp ascent close to $x = 1$ in~\fref{fig:alpha10_p2q3_cd_a}. The maximum errors in 
Figs.~\ref{fig:alpha10_p2q3_cd_c} and~\ref{fig:alpha10_p2q3_cd_d} for
$u$ and $u^\prime$ are about $6 \times 10^{-3}$
and $7 \times 10^{-5}$, respectively. 
A convergence study is conducted using
$n = [ 2,\,4,\,8,\,16,\,32,\,64]$ and 
three different choices for $p$ and $q$: 
$p = 2$, $q = 3$; 
$p = 2$, $q = 4$; and $p = 3$, $q = 4$.
The relative $L^2$ norm and relative $H^1$ seminorm of 
the error 
$u - u_\theta$ 
are defined as:
\begin{equation*}
E_u = \sqrt{ \frac{\int_0^1 (u - u_\theta)^2 \, dx} {\int_0^1 u^2 \, dx} }, \qquad
E_q = \sqrt{ \frac{\int_0^1 (u^\prime - q_\theta)^2 \, dx} {\int_0^1 (u^\prime)^2 \,dx} }.
\end{equation*}
Plots of the relative errors versus the number 
of degrees of freedom are presented
in~\fref{fig:repu_cd_alpha10_convergence}.  
The rates of convergence (slopes) in $(u,u^\prime)$
are $(2,3)$, $(2,3.2)$ and $(2.9,4)$.
\revised{Note that errors of ${\cal O}(10^{-5})$ are observed 
in the convergence plots shown in~\fref{fig:repu_cd_alpha10_convergence}. 
Since the Galerkin method is convergent, it is systematically improvable and these errors would decrease with refinement. Since we use a shallow neural network with a single hidden layer that has a pre-determined (fixed) bias and weight, the unknown coefficients only arise from the output layer, which are obtained via the solution of a singular but consistent
linear system of equations. This involves neither training nor tuning of network parameters as is needed in PINNs. 
}
\begin{figure}[!tbh]
\centering
\begin{subfigure}{0.48\textwidth}
\includegraphics[width=\textwidth]{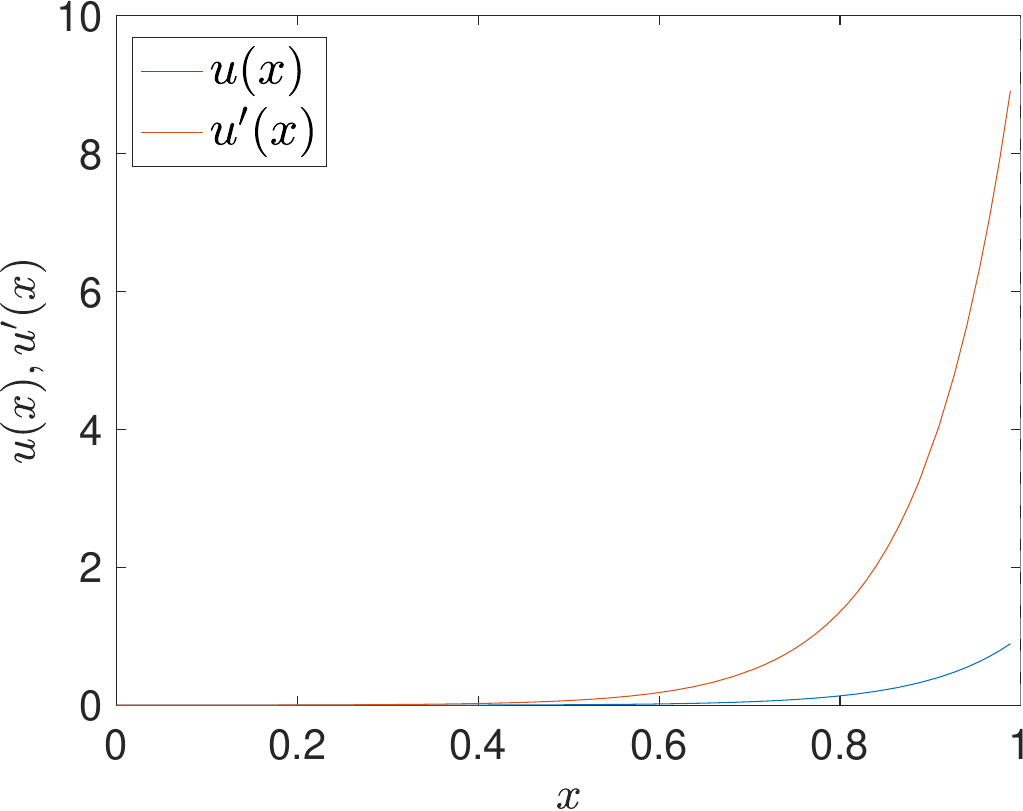}
\subcaption{}\label{fig:alpha10_p2q3_cd_a}
\end{subfigure} \hfill
\begin{subfigure}{0.48\textwidth}
\includegraphics[width=\textwidth]{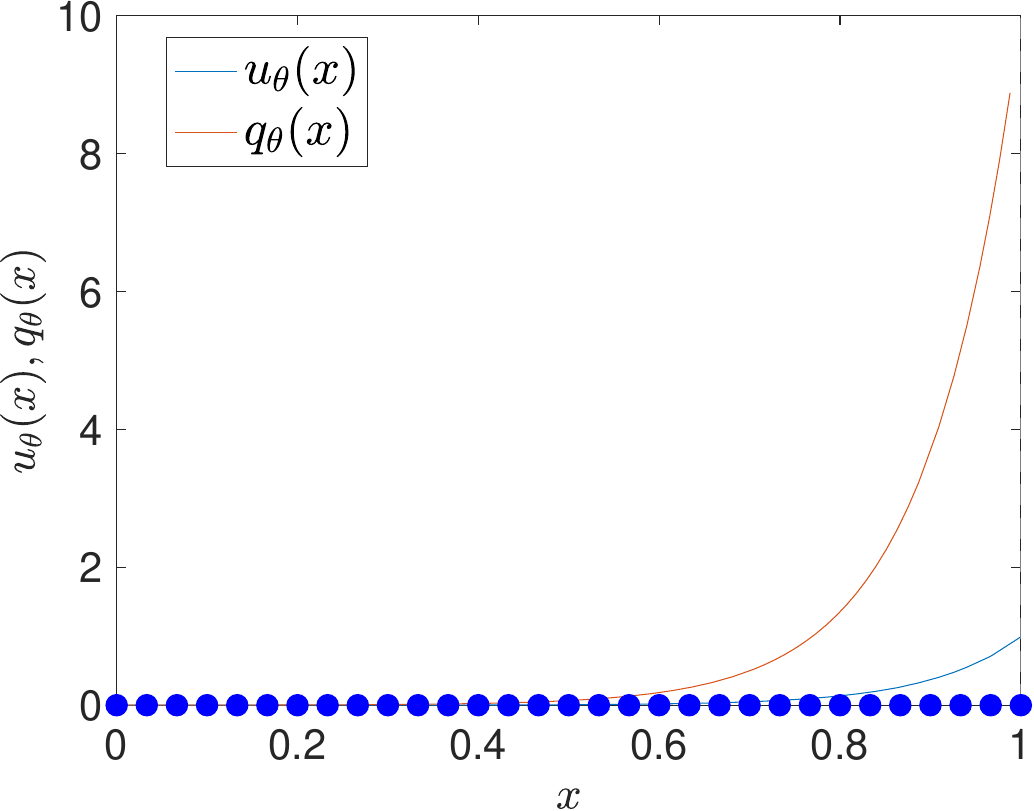}
\subcaption{}\label{fig:alpha10_p2q3_cd_b}
\end{subfigure} \hfill
\begin{subfigure}{0.48\textwidth}
\includegraphics[width=\textwidth]{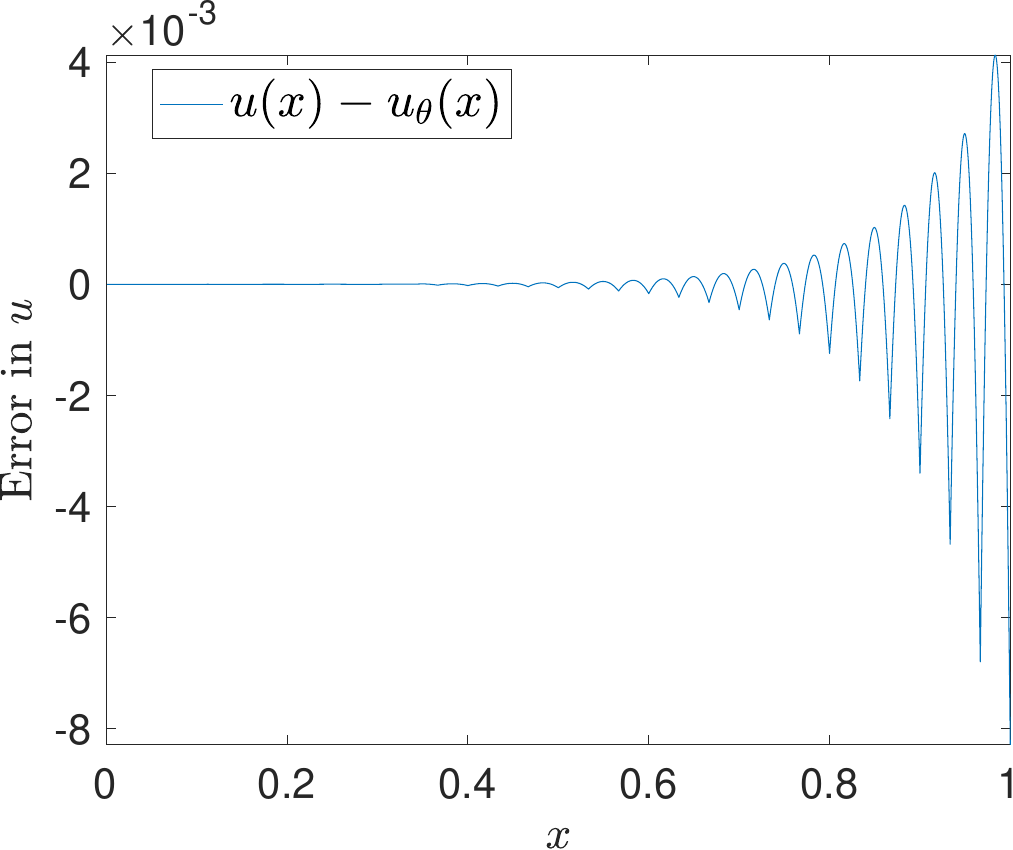}
\subcaption{}\label{fig:alpha10_p2q3_cd_c}
\end{subfigure} \hfill
\begin{subfigure}{0.48\textwidth}
\includegraphics[width=\textwidth]{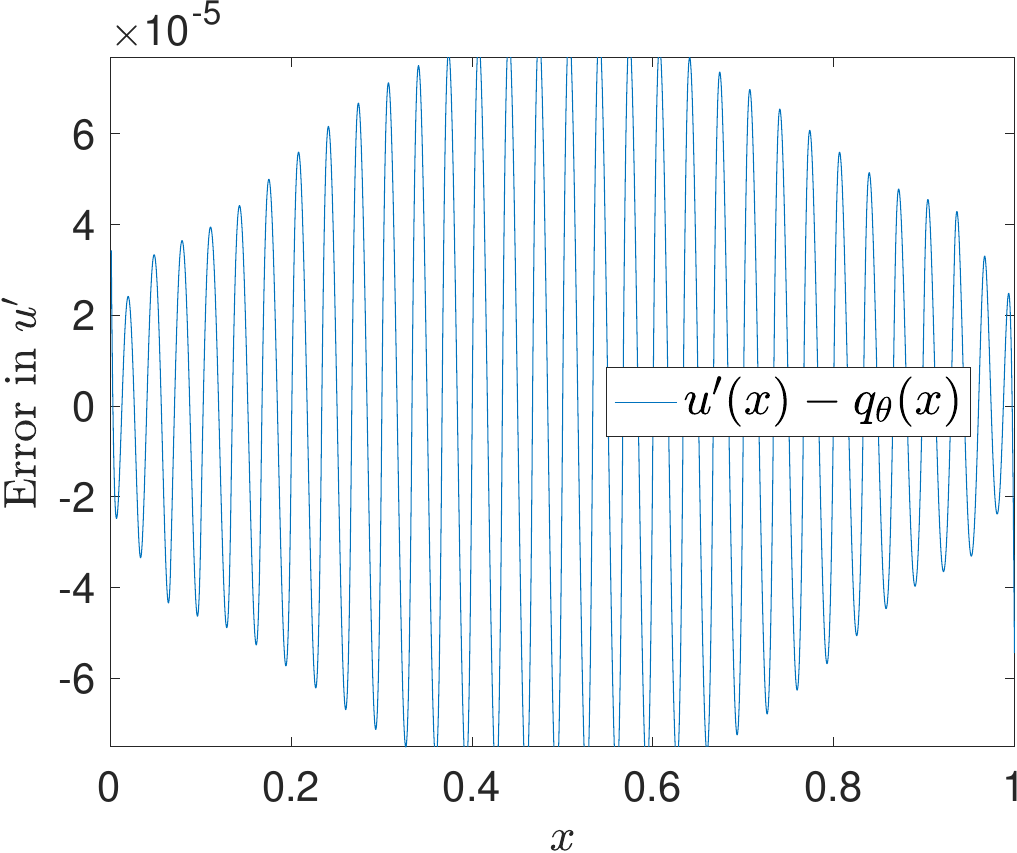}
\subcaption{}\label{fig:alpha10_p2q3_cd_d}
\end{subfigure}
\caption{Neural network solution for the steady-state 
         convection-diffusion problem ($\alpha = 10$) using
         $n = 30$, $p = 2$ and $q = 3$. (a) Exact solutions;
         (b) Neural network solutions;
         (c) $u - u_\theta$; and (d) $u^\prime - q_\theta$.}
         \label{fig:alpha10_p2q3_cd}
\end{figure}
% Convergence study for alpha = 10
\begin{figure}[!tbh]
\centering
\begin{subfigure}{0.33\textwidth}
\includegraphics[width=\textwidth]{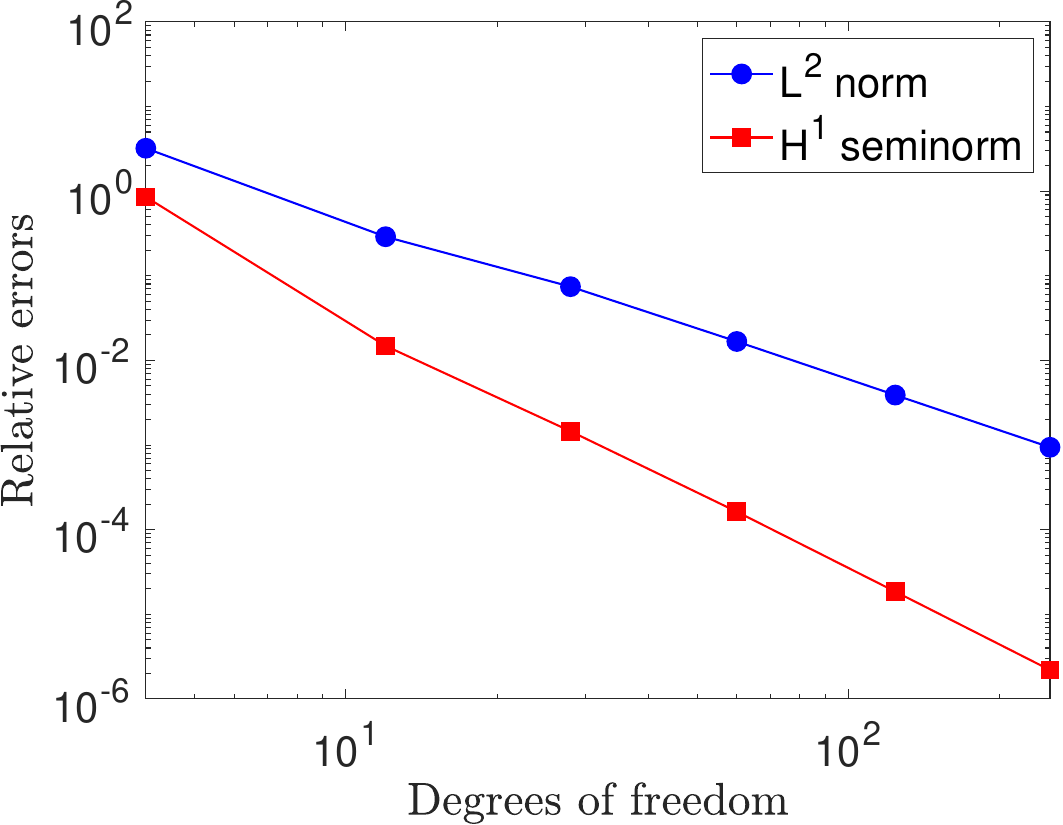}
\subcaption{}\label{fig:repu_cd_alpha10_convergence_a}
\end{subfigure} \hfill
\begin{subfigure}{0.33\textwidth}
\includegraphics[width=\textwidth]{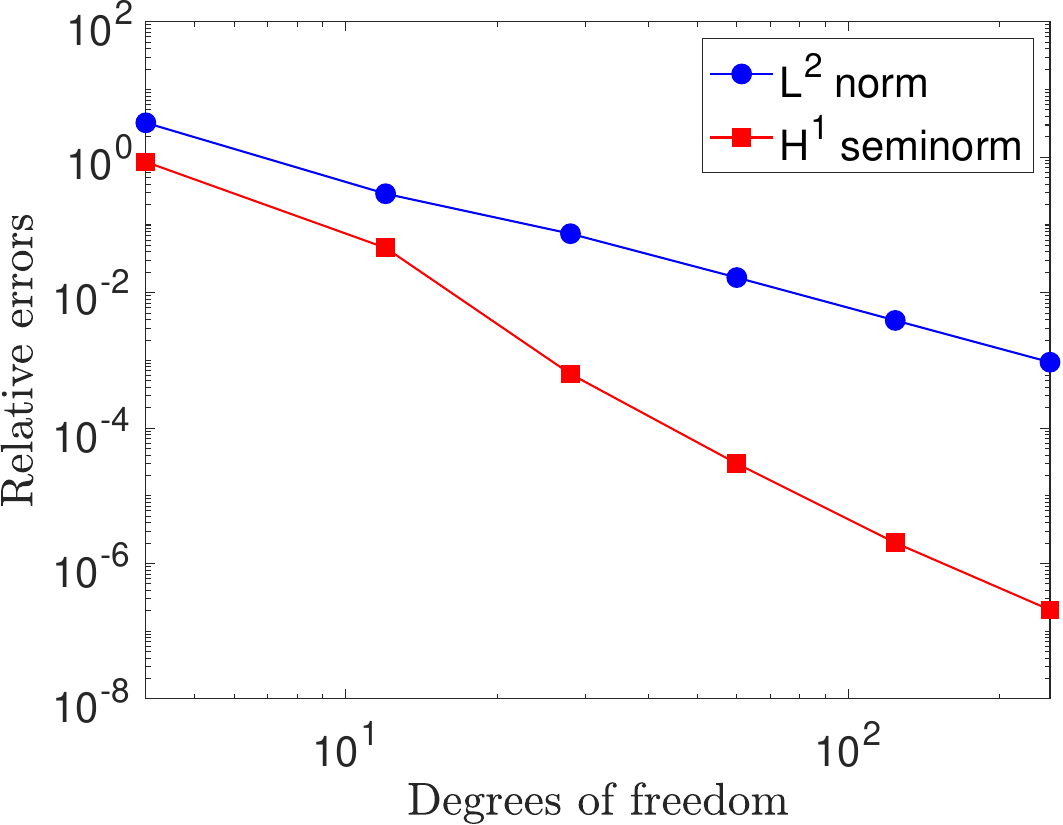}
\subcaption{}\label{fig:repu_cd_alpha10_convergence_b}
\end{subfigure} \hfill
\begin{subfigure}{0.33\textwidth}
\includegraphics[width=\textwidth]{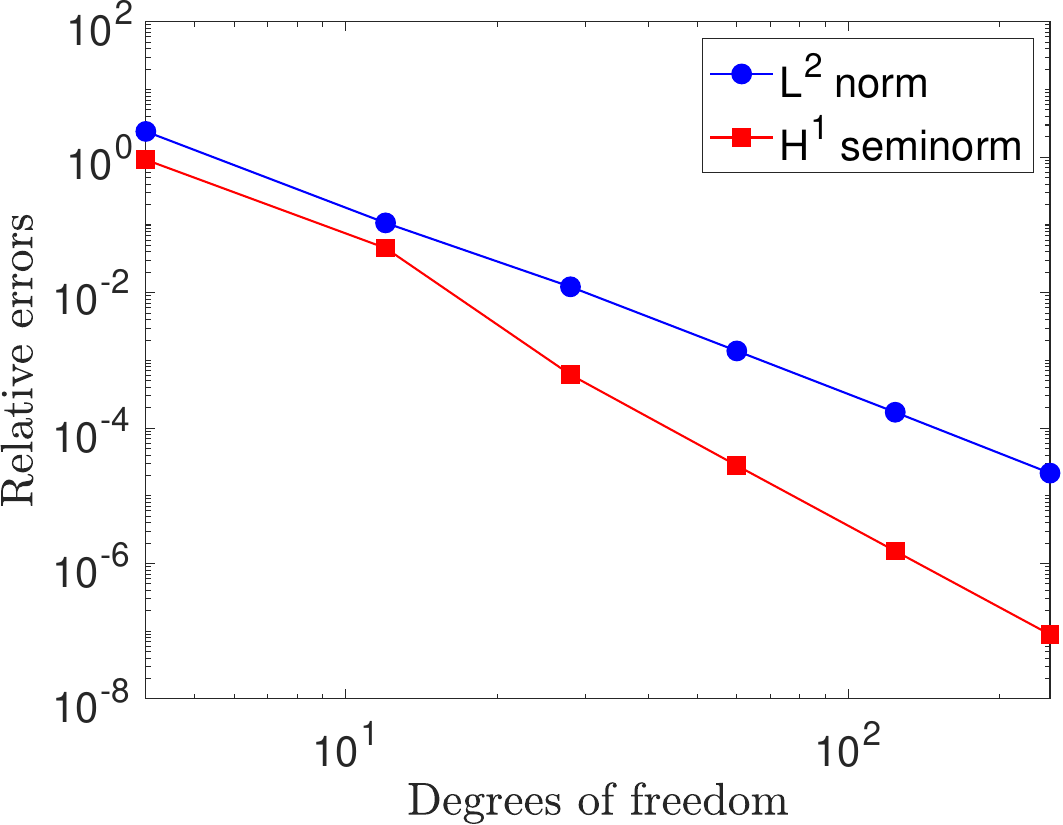}
\subcaption{}\label{fig:repu_cd_alpha10_convergence_c}
\end{subfigure} 
\caption{Convergence study
         with neural network approximants for the steady-state
         convection-diffusion problem ($\alpha = 10$). 
         The dual field $\mu_\theta(x)$ is approximated
         using RePU functions of degree $p$ and the dual field
         $\lambda_\theta(x)$ is formed by barycentric convex
         combination of RePU functions of degree $q$. 
         (a) $p = 2$ and $q = 3$; 
         (b) $p = 2$ and $q = 4$; and
         (c) $p = 3$ and $q = 4$.
         The rates of convergence in $u$ and $u^\prime$
         are of order $p$ and $p+1$, respectively.
         }\label{fig:repu_cd_alpha10_convergence}
\end{figure}

Finally, we set $\alpha = 50$, and the solution
$u$ develops a sharp boundary layer in the vicinity of $x = 1$ . The neural network solution is compared to the exact
solution in~\fref{fig:alpha50_p2q3_cd}. We 
choose $n = 50$, and $\mu_\theta(x)$ ($p = 2$ and
$\lambda_\theta(x)$ ($q = 3$) are selected as piecewise quadratic and
piecewise cubic polynomials, respectively. The maximum errors in 
Figs.~\ref{fig:alpha50_p2q3_cd_c} and~\ref{fig:alpha50_p2q3_cd_d} for
$u$ and $u^\prime$ are about $6 \times 10^{-2}$
and $1.5 \times 10^{-2}$, respectively. The larger errors in $u$
and $u^\prime$ are due to the steep gradient of $u$ in 
the vicinity of $x = 1$, which can be captured by either increasing $n$
($h$-refinement) or using a high-order approximation ($p$-refinement).
Observe that
$u \in (0,1]$ and $u^\prime \in (0,50]$ 
for $x \in [0.8,1]$. We revisit this problem using B-splines in the 
next section.
\begin{figure}[!tbh]
\centering
\begin{subfigure}{0.45\textwidth}
\includegraphics[width=\textwidth]{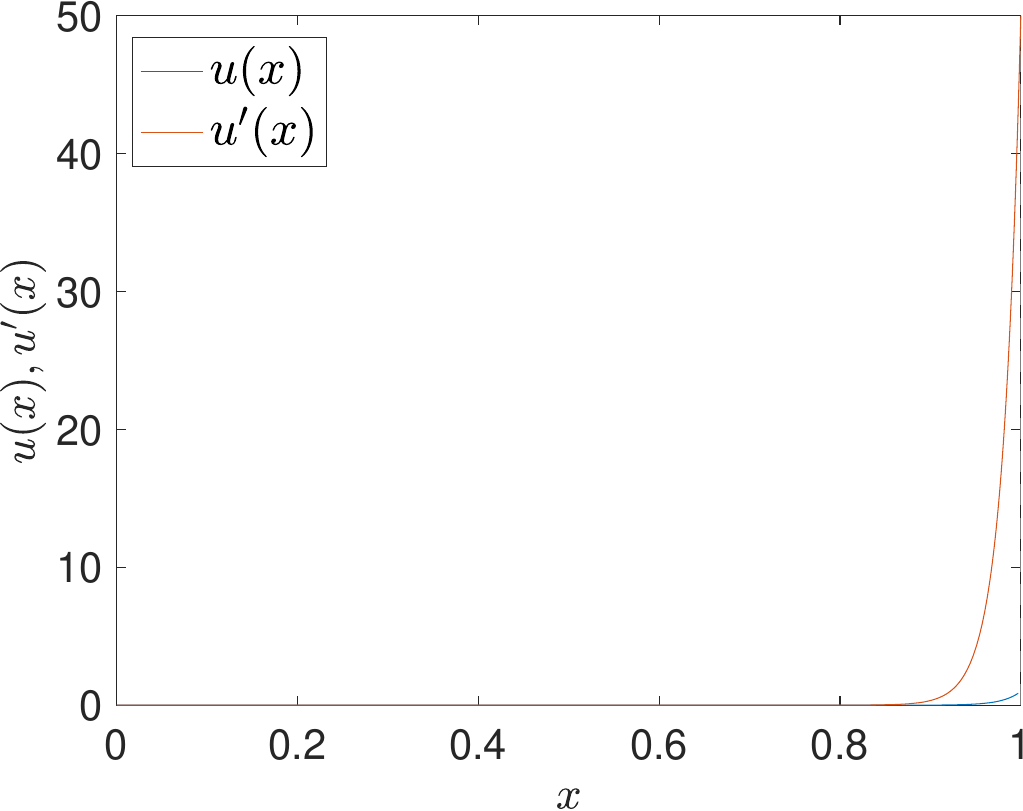}
\subcaption{}\label{fig:alpha50_p2q3_cd_a}
\end{subfigure} \hfill
\begin{subfigure}{0.45\textwidth}
\includegraphics[width=\textwidth]{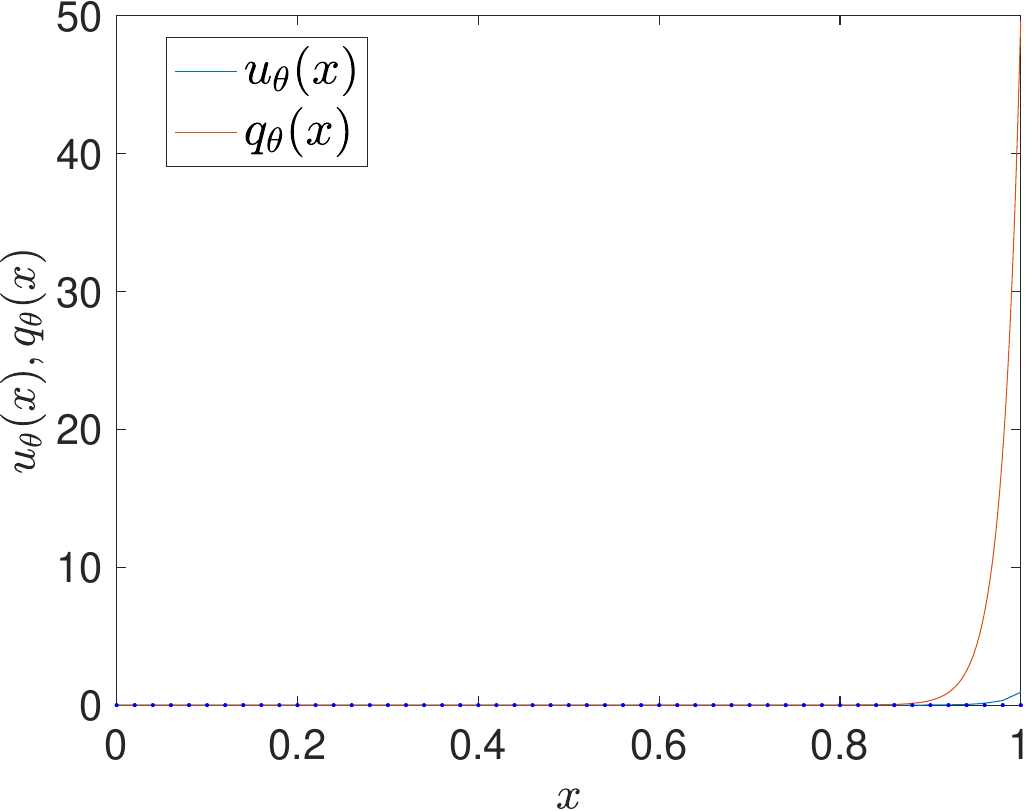}
\subcaption{}\label{fig:alpha50_p2q3_cd_b}
\end{subfigure} \hfill
\begin{subfigure}{0.45\textwidth}
\includegraphics[width=\textwidth]{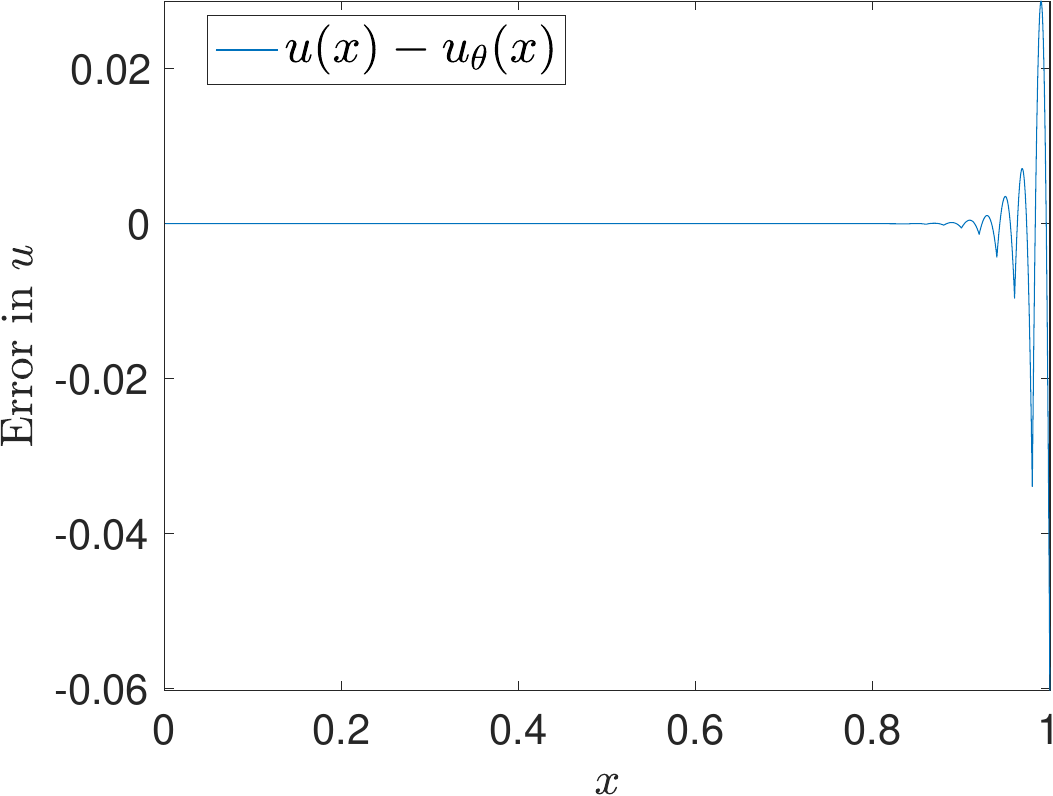}
\subcaption{}\label{fig:alpha50_p2q3_cd_c}
\end{subfigure} \hfill
\begin{subfigure}{0.45\textwidth}
\includegraphics[width=\textwidth]{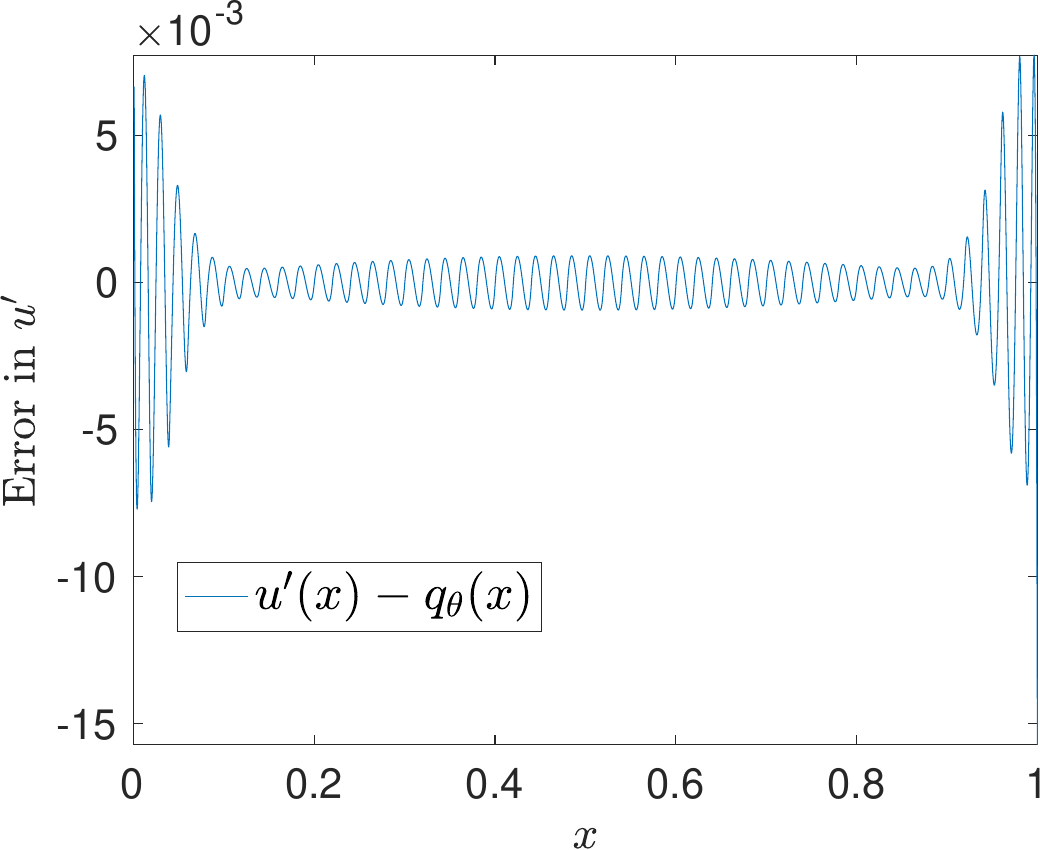}
\subcaption{}\label{fig:alpha50_p2q3_cd_d}
\end{subfigure}
\caption{Neural network solution for the steady-state 
         convection-diffusion problem ($\alpha = 50$) using
         $n = 50$, $p = 2$ and $q = 3$. 
         (a) Exact solutions;
         (b) Neural network solutions;
         (c) $u - u_\theta$; and (d) $u^\prime - q_\theta$.}
         \label{fig:alpha50_p2q3_cd}
\end{figure}

\subsubsection{Univariate B-spline solutions}\label{subsubsec:bsp}
We adopt a univariate B-spline approximation~\cite{deBoor:2001:PGS}
to represent the dual fields $\mu(x)$ and $\lambda(x)$. 
\revised{Unlike the RePU functions used in the previous section,
B-splines are compactly-supported and are linearly
independent.}
An open uniform knot vector is used:
\begin{equation}\label{eq:bsp_knot}
  \Xi := [\underset{p+1}{\underbrace{0,0,\dots,0}}, x_1, \dots, x_{n-1},
             \underset{p+1}{\underbrace{1,1,\dots,1}}],
\end{equation}
where the end knot points are repeated $p+1$ times so that $C^0$ continuity and interpolation is ensured at $x = 0$ and $x = 1$. B-splines are used to construct smooth approximations for the dual fields
$\mu(x)$ and $\lambda(x)$:
\begin{subequations}\label{eq:bsp_mulambda}
\begin{align}
  \mu_\theta(x) &= \sum_{i=1}^{n+p} B_i^p(x) a_i
                              \in C^{p-1}(\Omega), 
                                 \label{eq:bsp_mulambada_a} \\
  \lambda_\theta(x) &= \sum_{i=1}^{n+q} B_i^q(x) b_i
                                   \in C^{q-1}(\Omega), 
                                   \label{eq:bsp_mulambda_b}
\end{align}
\end{subequations}
where $B_i^p(x)$ is the $p$-th degree B-spline basis function
with control points ($\vm{a}$ and $\vm{b}$) as the unknowns. 
B-spline computations are performed using the Cox--de Boor algorithm. Note that on setting $b_1 = b_{n+q} = 0$ 
in~\eref{eq:bsp_mulambda_b}, we have
$\lambda_\theta(0) = \lambda_\theta(1) = 0$, which
ensures that the trial functions are admissible.
In one dimension, the approximations 
in~\eref{eq:bsp_mulambda} can be viewed as a shallow network within 
the KAN architecture~\cite{Liu:2024:KAN}. 

For $n = 20$, the nonzero 
B-spline basis functions that are used
to form $\mu_\theta(x)$ ($p = 2$) and $\lambda_\theta(x)$ ($q = 3$)
are presented in~\fref{fig:bsp_p2q3_mulambda}. 
We set the Peclet number $\alpha$ to $50$. For the numerical 
computations, we choose $n = 20$, and vary 
$p$ (degree of $\mu_\theta$) 
and $q$ (degree of $\lambda_\theta$).
The B-spline basis functions to form
$\mu_\theta(x)$ ($p = 2$) and $\lambda_\theta(x)$ ($q = 3$) are 
shown in Figs.~\ref{fig:bsp_p2q3_mulambda_a} 
and~\ref{fig:bsp_p2q3_mulambda_b}, respectively.
\begin{figure}[!thb]
\centering
\begin{subfigure}{0.48\textwidth}
\includegraphics[width=\textwidth]{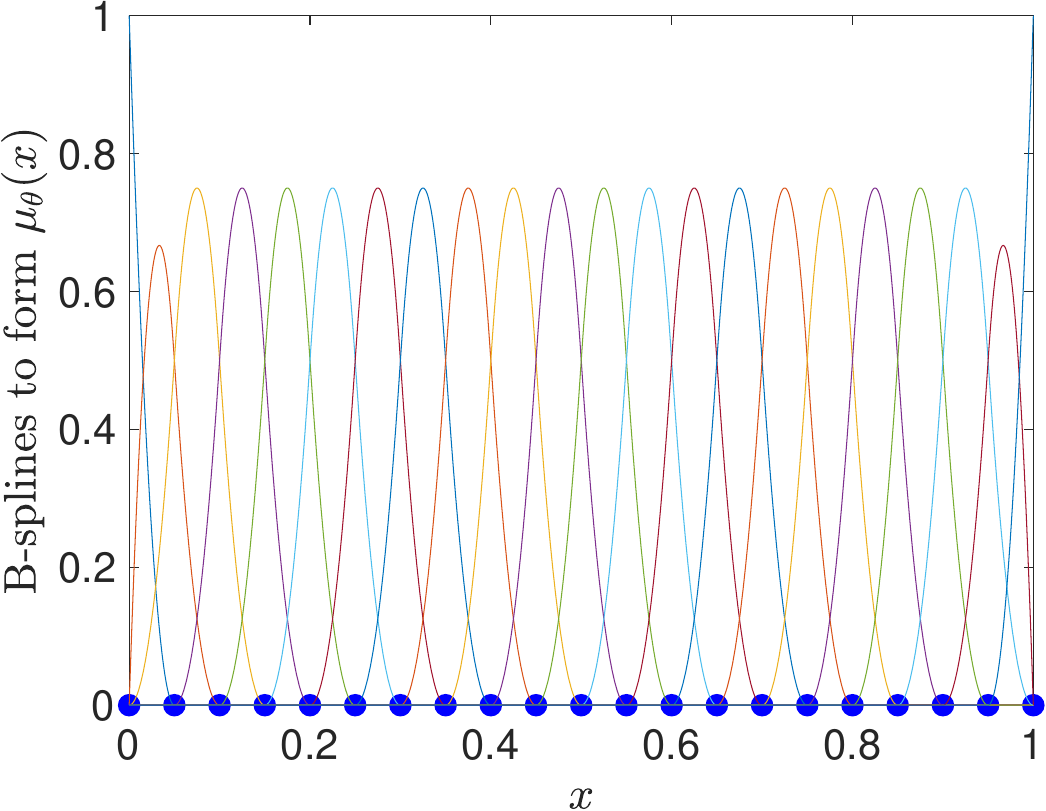}
\subcaption{}\label{fig:bsp_p2q3_mulambda_a}
\end{subfigure} \hfill
\begin{subfigure}{0.45\textwidth}
\includegraphics[width=\textwidth]{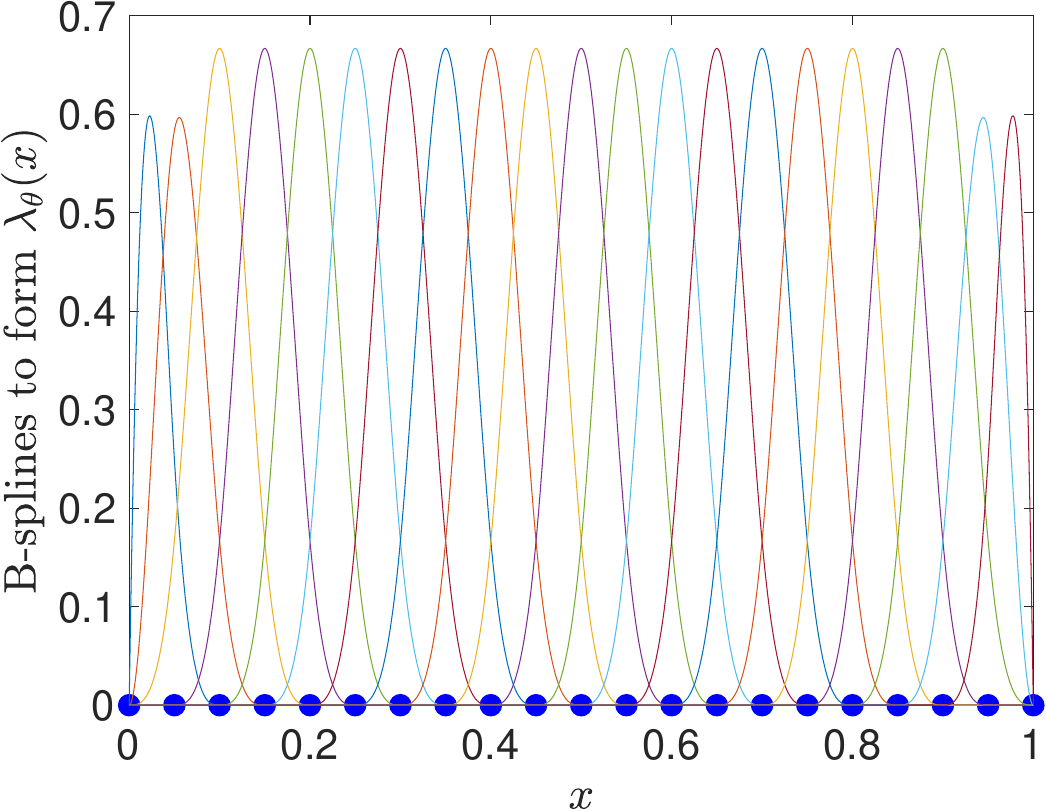}
\subcaption{}\label{fig:bsp_p2q3_mulambda_b}
\end{subfigure}
\caption{Plots of B-spline basis functions that are used to form the 
         dual fields for $n = 20$. 
         (a) $\mu_\theta(x)$ ($p = 2$) and 
         (b) $\lambda_\theta(x)$ ($q = 3$).}
         \label{fig:bsp_p2q3_mulambda}
\end{figure}
% p = 2, q = 3; q = 5, q = 6 (alpha = 50)
For $p = 2$ and $q = 3$, the dual fields are presented
in~\fref{fig:alpha50_p2q3_bsp_cd_dualfields}. 
The primal fields are computed from the dual fields using the DtP
map in~\eqref{eq:dual_steadystate_cd_a}.
The exact solution for $u$ and $u^\prime$ is presented
in~\fref{fig:alpha50_p2q3_cd_a}. The 
errors in the B-spline solution for two choices of 
$p$ and $q$
are presented in~\fref{fig:alpha50_p2q3p5q6_bsp_cd}. 
For $p = 2$ and $q = 3$, we see from
Figs.~\ref{fig:alpha50_p2q3p5q6_bsp_cd_a} and~\ref{fig:alpha50_p2q3p5q6_bsp_cd_b} that
the maximum error in $u$ and $u^\prime$ are 
about $0.2$ and $2$, respectively. For $p = 5$ and $q = 6$,
the maximum error in $u$ and $u^\prime$ from
Figs.~\ref{fig:alpha50_p2q3p5q6_bsp_cd_c} and~\ref{fig:alpha50_p2q3p5q6_bsp_cd_d} are
$4 \times 10^{-3}$ and $8 \times 10^{-3}$, respectively.
\begin{figure}[!tbh]
\centering
\begin{subfigure}{0.48\textwidth}
\includegraphics[width=\textwidth]{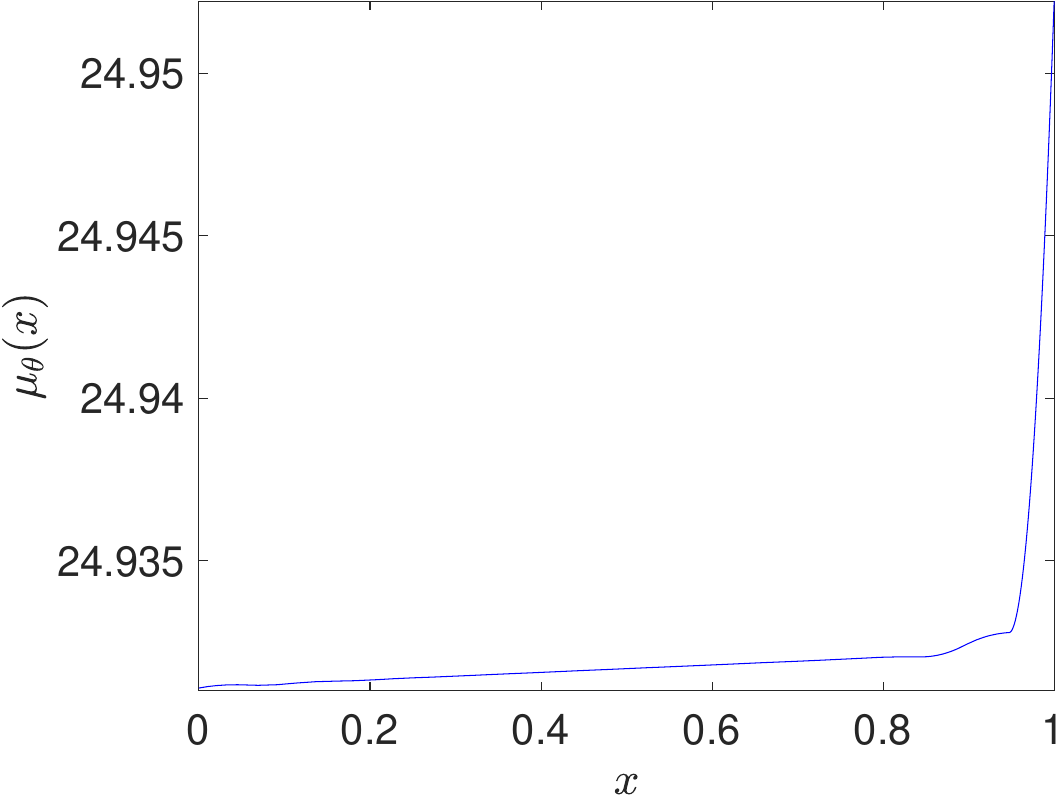}
\subcaption{}\label{fig:alpha50_p2q3_bsp_cd_dualfields_a}
\end{subfigure} \hfill
\begin{subfigure}{0.45\textwidth}
\includegraphics[width=\textwidth]{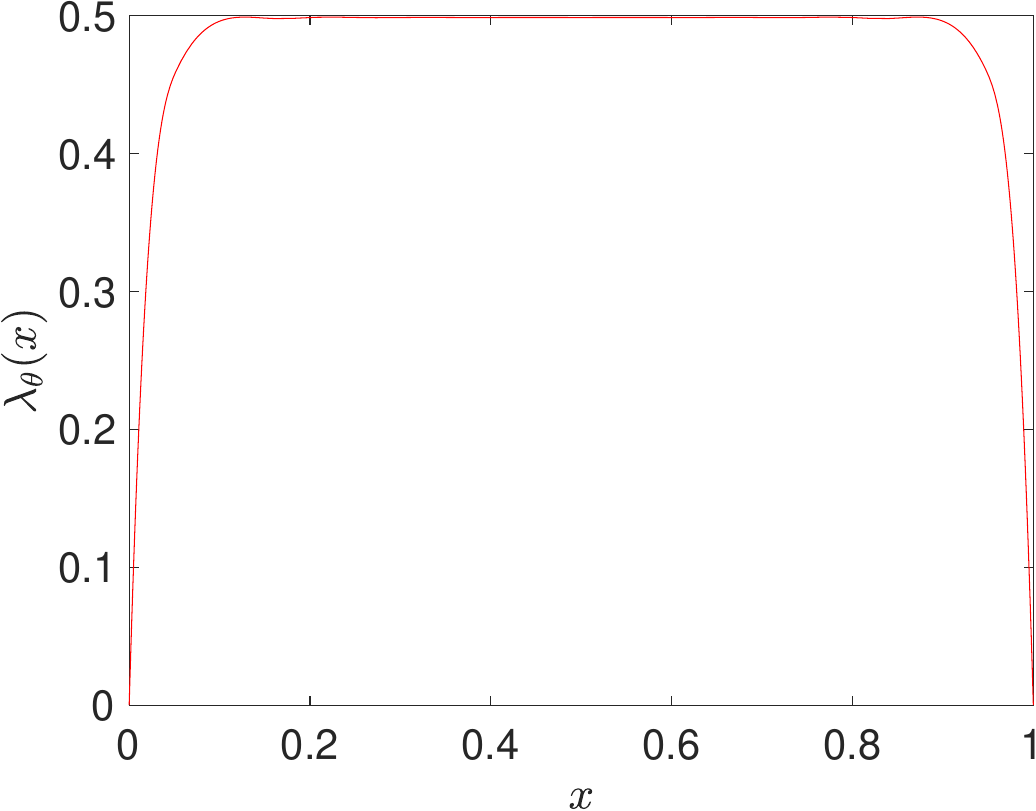}
\subcaption{}\label{fig:alpha50_p2q3_bsp_cd_dualfields_b}
\end{subfigure}
\caption{B-spline solution (dual fields)
         of the steady-state convection-diffusion
         equation ($\alpha = 50$) for $n = 20$.
         (a) $\mu_\theta(x)$ ($p = 2$) and 
         (b) $\lambda_\theta(x)$ ($q = 3$).
         }
         \label{fig:alpha50_p2q3_bsp_cd_dualfields}
\end{figure}
\begin{figure}[!tbh]
\centering
\begin{subfigure}{0.48\textwidth}
\includegraphics[width=\textwidth]{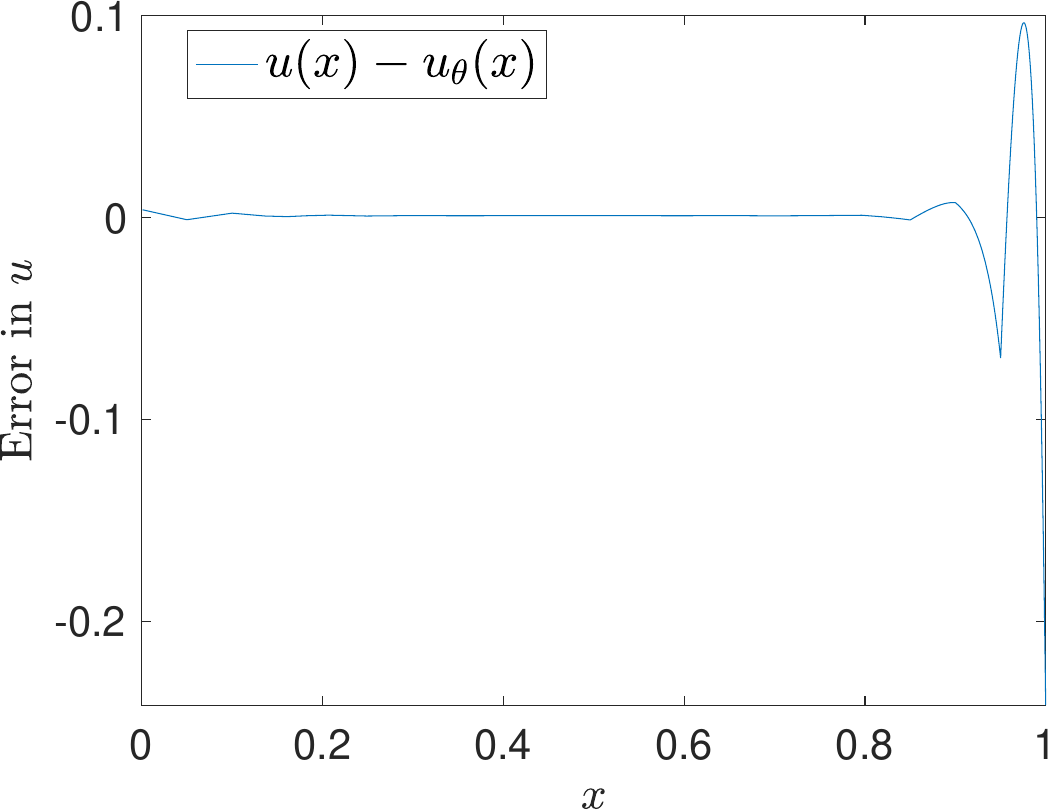}
\subcaption{}\label{fig:alpha50_p2q3p5q6_bsp_cd_a}
\end{subfigure} \hfill
\begin{subfigure}{0.48\textwidth}
\includegraphics[width=\textwidth]{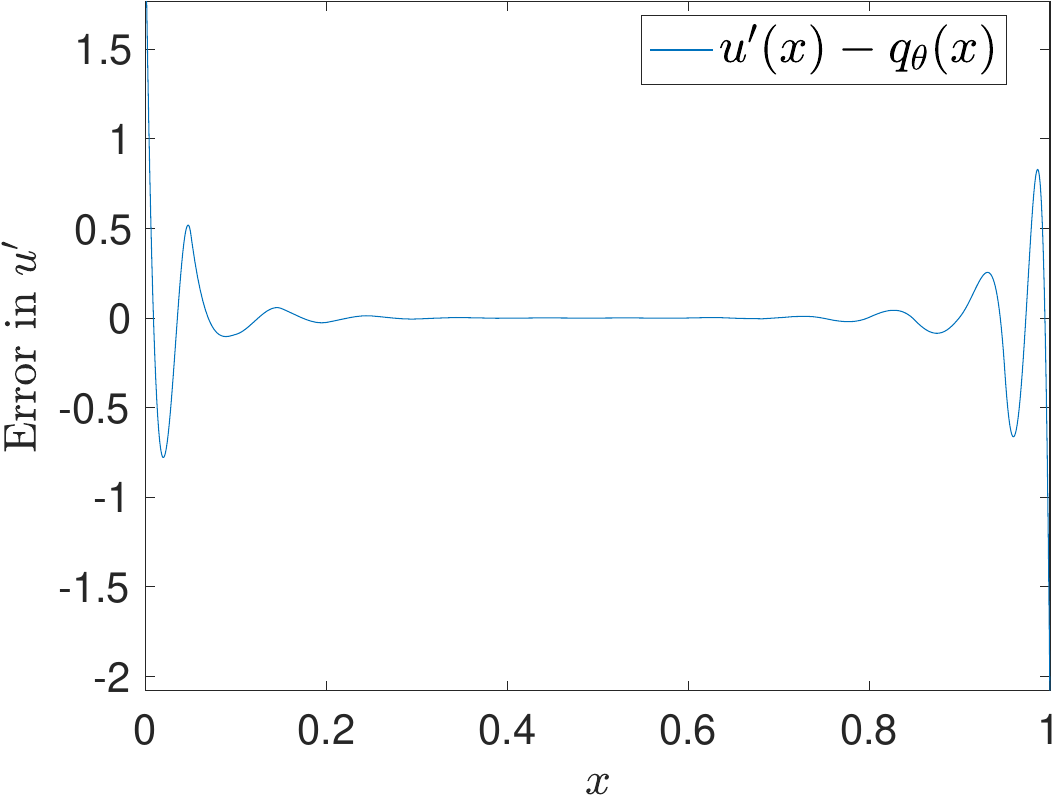}
\subcaption{}\label{fig:alpha50_p2q3p5q6_bsp_cd_b}
\end{subfigure} \hfill
\begin{subfigure}{0.48\textwidth}
\includegraphics[width=\textwidth]{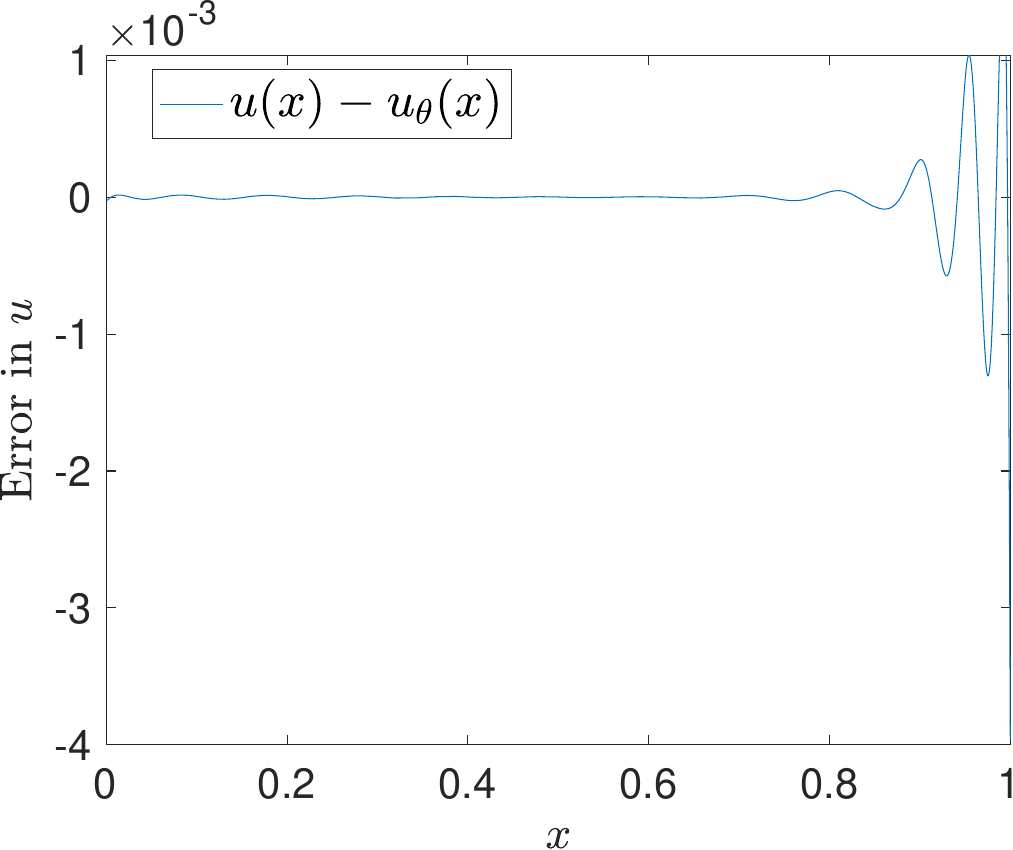}
\subcaption{}\label{fig:alpha50_p2q3p5q6_bsp_cd_c}
\end{subfigure} \hfill
\begin{subfigure}{0.48\textwidth}
\includegraphics[width=\textwidth]{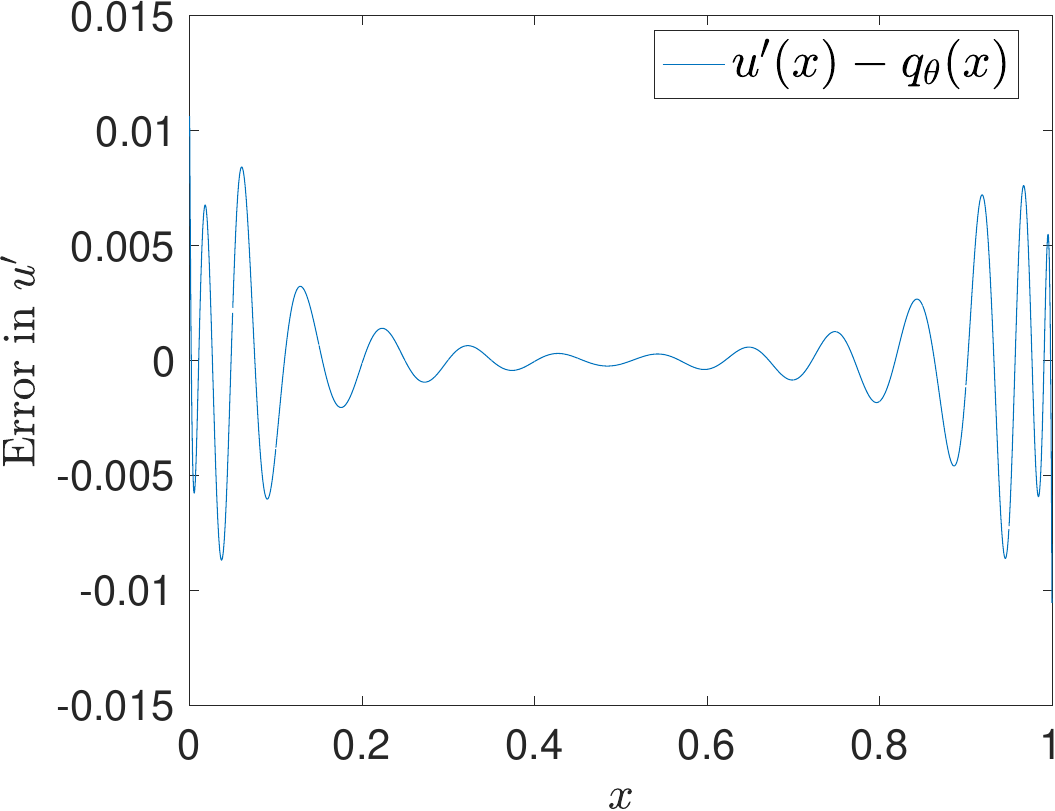}
\subcaption{}\label{fig:alpha50_p2q3p5q6_bsp_cd_d}
\end{subfigure}
\caption{B-spline solution for the steady-state 
         convection-diffusion problem ($\alpha = 50$).
         For $n = 20, \ p = 2, \ q = 3$:
         (a) $u - u_\theta$ and (d) $u^\prime - q_\theta$.
         For $n = 20, \ p = 5, \ q = 6$:
         (c) $u - u_\theta$ and (d) $u^\prime - q_\theta$.}
         \label{fig:alpha50_p2q3p5q6_bsp_cd}
\end{figure}
% p = 7, q = 8 (alpha = 50)
For $n = 20$, B-spline basis functions to form
$\mu_\theta(x)$ ($p = 7$) and $\lambda_\theta(x)$ ($q = 8$) are shown 
in~\fref{fig:alpha50_p7q8_bsp_cd_a} 
and~\fref{fig:alpha50_p7q8_bsp_cd_b}, respectively. The 
errors in $u$ and
$u^\prime$ are presented in Figs.~\ref{fig:alpha50_p7q8_bsp_cd_c} 
and~\ref{fig:alpha50_p7q8_bsp_cd_d}, respectively.  
The maximum errors in $u$ and $u^\prime$ are  
$1.25 \times 10^{-4}$ and $1.25 \times 10^{-3}$, respectively.
Since $||u||_\infty = 1$ and $||u^\prime||_\infty = 50$,
the relative
errors for $u$ and $u^\prime$ that B-splines deliver are proximal.
\begin{figure}[!tbh]
\centering
\begin{subfigure}{0.48\textwidth}
\includegraphics[width=\textwidth]{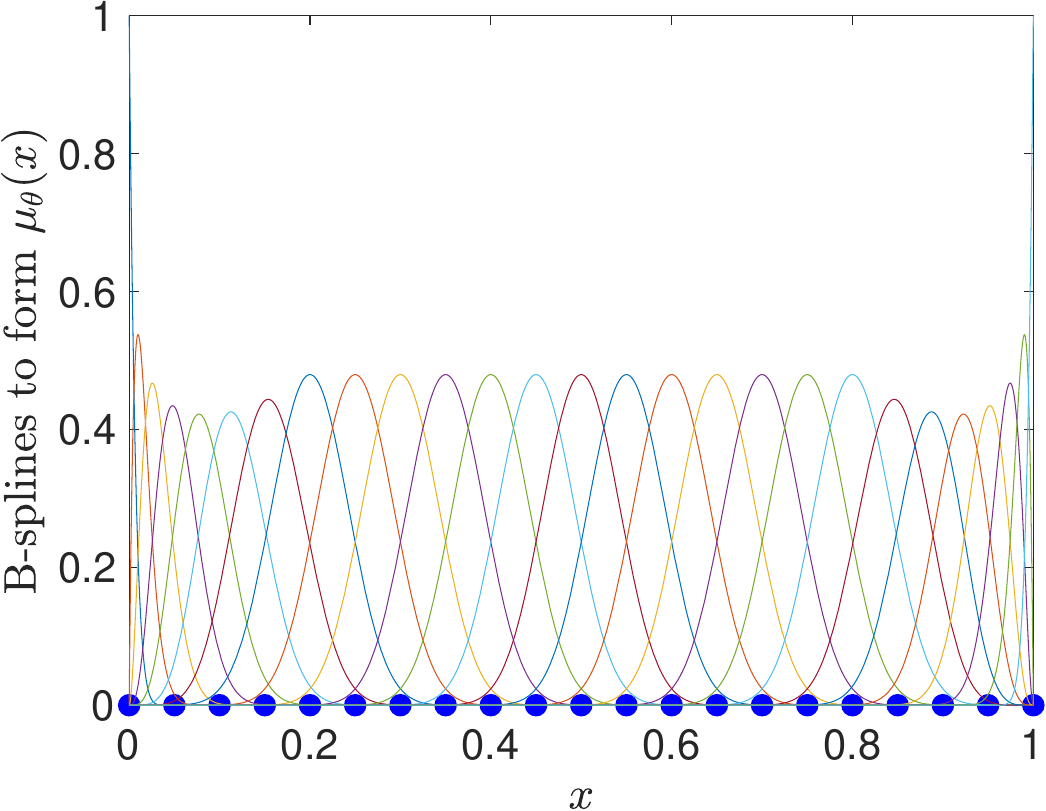}
\subcaption{}\label{fig:alpha50_p7q8_bsp_cd_a}
\end{subfigure} \hfill
\begin{subfigure}{0.48\textwidth}
\includegraphics[width=\textwidth]{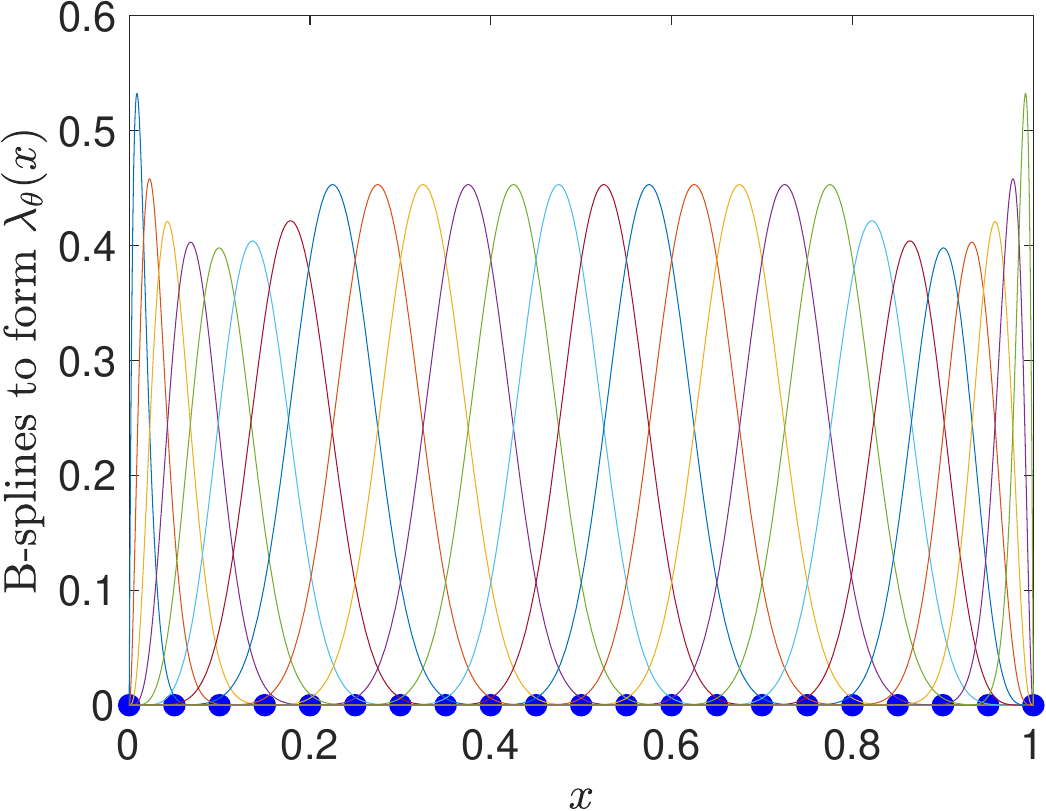}
\subcaption{}\label{fig:alpha50_p7q8_bsp_cd_b}
\end{subfigure} \hfill
\begin{subfigure}{0.48\textwidth}
\includegraphics[width=\textwidth]{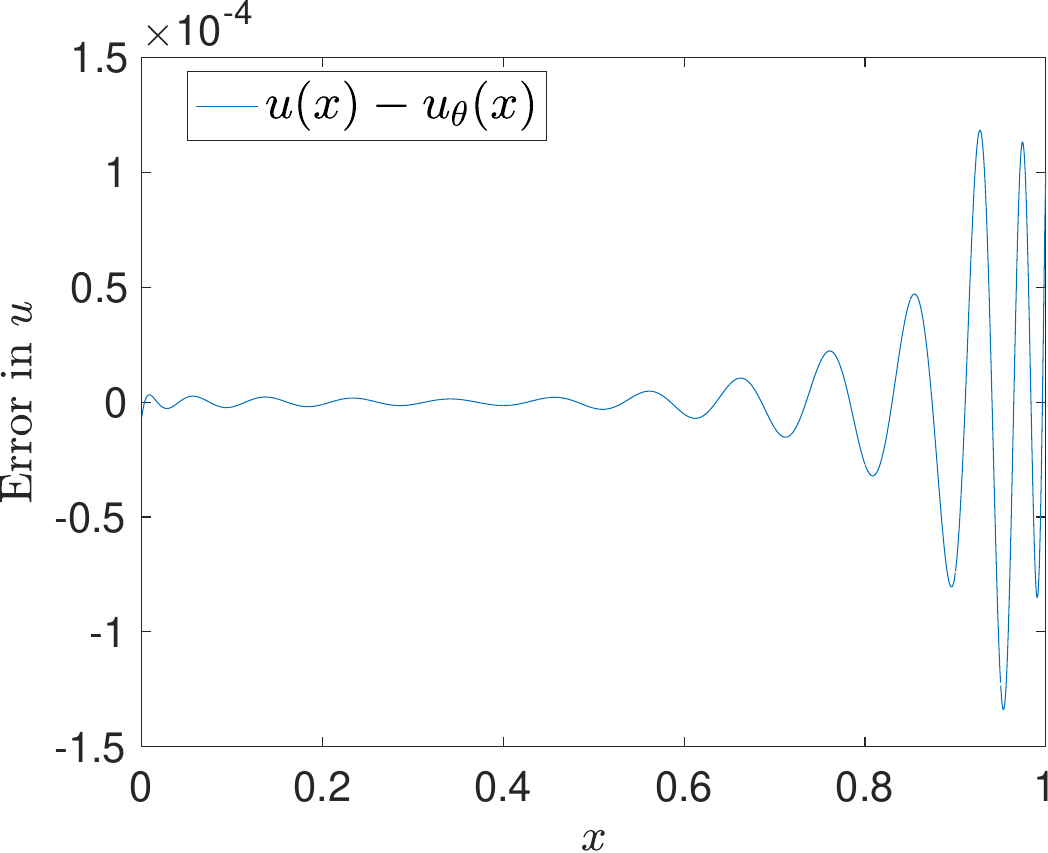}
\subcaption{}\label{fig:alpha50_p7q8_bsp_cd_c}
\end{subfigure} \hfill
\begin{subfigure}{0.48\textwidth}
\includegraphics[width=\textwidth]{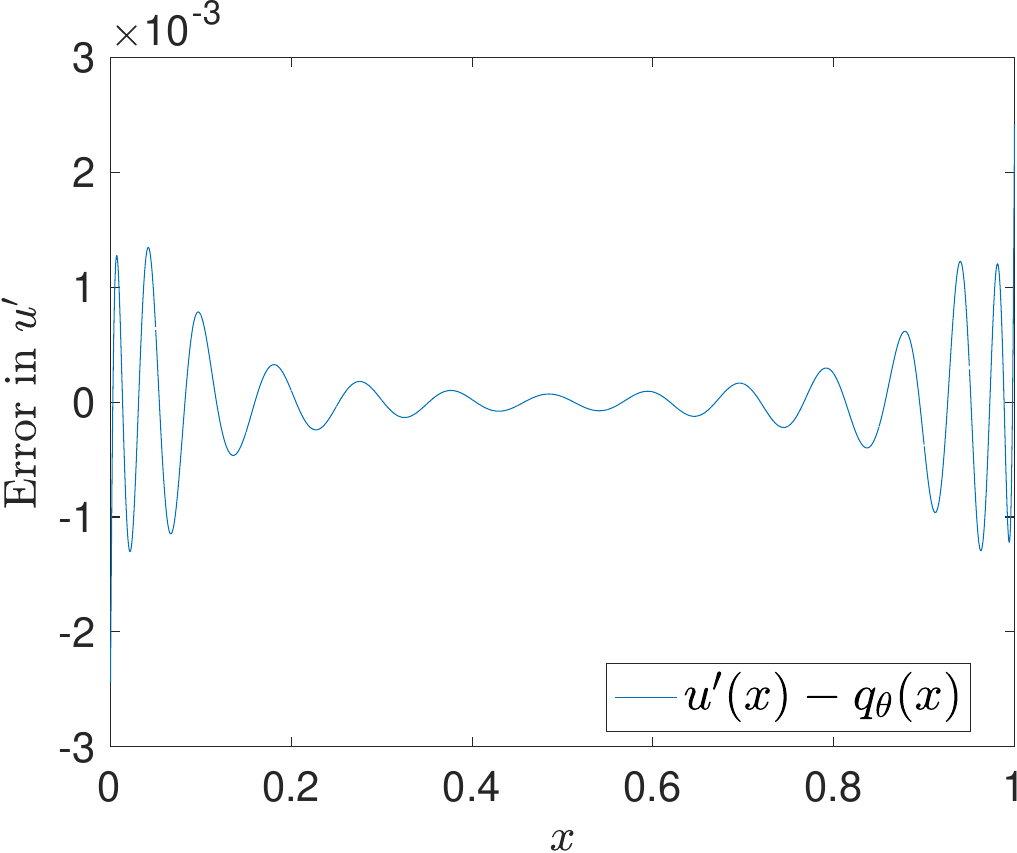}
\subcaption{}\label{fig:alpha50_p7q8_bsp_cd_d}
\end{subfigure}
\caption{B-spline computations to solve the steady-state 
         convection-diffusion problem ($\alpha = 50$). 
         For $n = 20$, B-spline basis functions to 
         form (a) $\mu_\theta(x)$ ($p = 7$) and
         (b) $\lambda_\theta(x)$ ($q = 8$);
         (c) $u - u_\theta$; and (d) $u^\prime - q_\theta$.}
         \label{fig:alpha50_p7q8_bsp_cd}
\end{figure}

We perform a convergence study with B-splines to assess the rate of convergence of the method. The dual fields 
$\mu_\theta(x) \in C^{p-1}(\Omega)$
and $\lambda_\theta(x) \in C^{q-1}(\Omega)$ 
are approximated using B-spline
basis functions of degree $p$ and $q$,
respectively. The primal fields are obtained from the
DtP map in~\eref{eq:dual_steadystate_cd_a}:
\begin{equation*}
 u_\theta(x) = \mu_\theta^\prime (x), \quad q_\theta (x) = 
 \mu_\theta (x) - \alpha \lambda_\theta (x) - 
 \lambda_\theta^\prime (x) .
\end{equation*}
\revised{One uses standard $L^2$ interpolation estimates for Galerkin methods with approximations that possess $p$-th degree completeness to arrive at the rates of convergence of $u$ and $u^\prime$.
From the DtP map and the fact that
B-spline approximations (degree $p$) possess $p$-th 
degree polynomial
completeness, one can deduce that the error in $u$ should
decay as ${\cal O}(h^p)$ in the $L^2$ norm ($h$ is the mesh spacing),
and likewise the error in the $L^2$ norm for the flux $u^\prime$ should decay as 
${\cal O}(h^p)$ if $q = p$ and
${\cal O}(h^{p+1})$ if $q = p+1$.}
A sequence of uniformly refined meshes
with  $n = [4,\, 8, \, 16, \, 32, \, 64]$ is selected.  
For $\alpha = 10$ and $\alpha = 50$, four different choices for
 $p$ and $q$ are considered: $p = q = 1$;
$p = 1$, $q = 2$; $p = 2$, $q = 3$; and $p = 3$, $q = 4$. The
total number of degrees of freedom is $2n + p + q - 2$.
In Figs.~\ref{fig:cd_alpha10_convergence} and~\ref{fig:cd_alpha50_convergence}, the plots
of the relative errors versus the 
number of degrees of freedom
are presented. Note that the convergence plots
in Figs.~\ref{fig:cd_alpha10_convergence_a} 
and~\ref{fig:cd_alpha50_convergence_a} are for 
linear $C^0$ finite elements
($p = q = 1$), where primal fields that are
discontinuous at the knot locations are used in
the error norm computations. For the convergence plots shown
in Figs.~\ref{fig:cd_alpha10_convergence_a}--\ref{fig:cd_alpha10_convergence_d},
the rates of convergence in $(u,u^\prime)$ are found to be
$(1,1)$, $(1,2)$, $(2.1,3)$, $(3.1,4.1)$, which is in agreement with 
a priori error estimates of ${\cal O}(h^p)$  
and ${\cal O}(h^{p+1})$ ($q = p+1$) 
in $u$ and $u^\prime$, respectively. 
For the convergence plots shown
in Figs.~\ref{fig:cd_alpha50_convergence_a}--\ref{fig:cd_alpha50_convergence_d},
the rates of convergence in $(u,u^\prime)$ are found to be
$(1.2,0.9)$, $(0.9,2)$, $(2,3)$ and $(3,3.7)$.
We consistently observe that 
the B-spline solution
for the derivative of $u$ is more accurate than 
the B-spline solution for $u$.
Due to the five-fold increase in $\alpha$, the
relative errors in~\fref{fig:cd_alpha50_convergence} are about 
two orders of magnitude less accurate than the corresponding results
shown in~\fref{fig:cd_alpha10_convergence}.
% Convergence study for alpha = 10
\begin{figure}[!tbh]
\centering
\begin{subfigure}{0.48\textwidth}
\includegraphics[width=\textwidth]{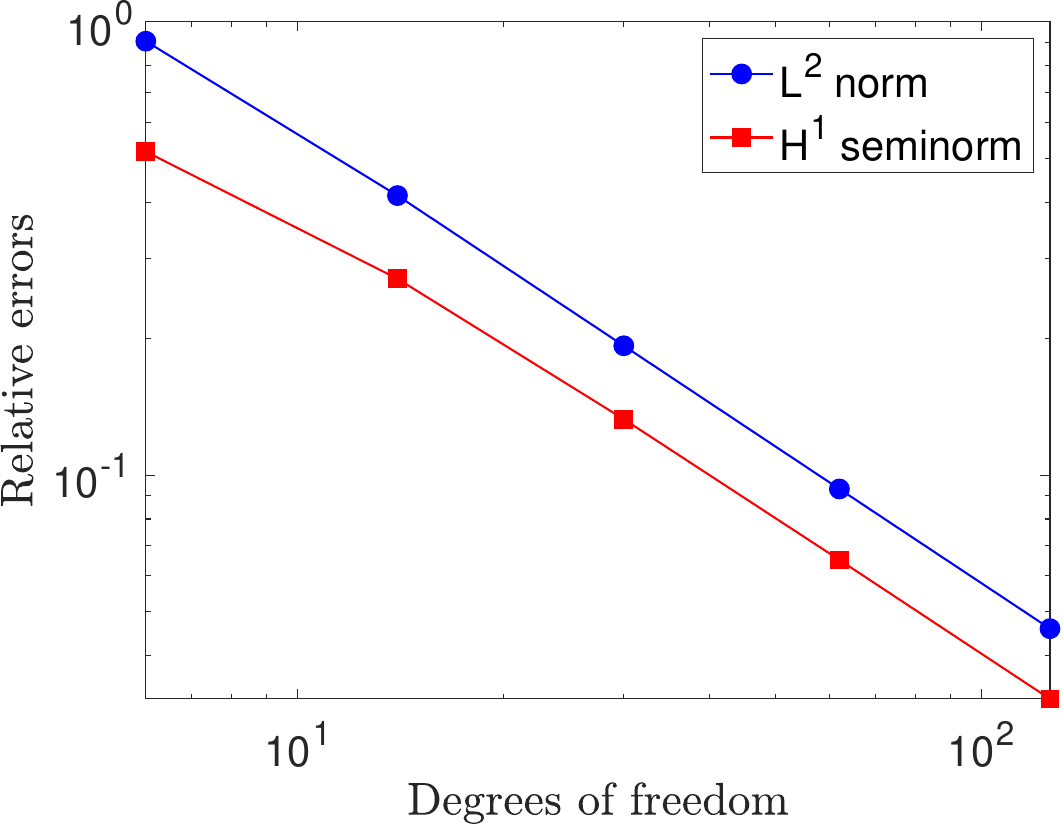}
\subcaption{}\label{fig:cd_alpha10_convergence_a}
\end{subfigure} \hfill
\begin{subfigure}{0.48\textwidth}
\includegraphics[width=\textwidth]{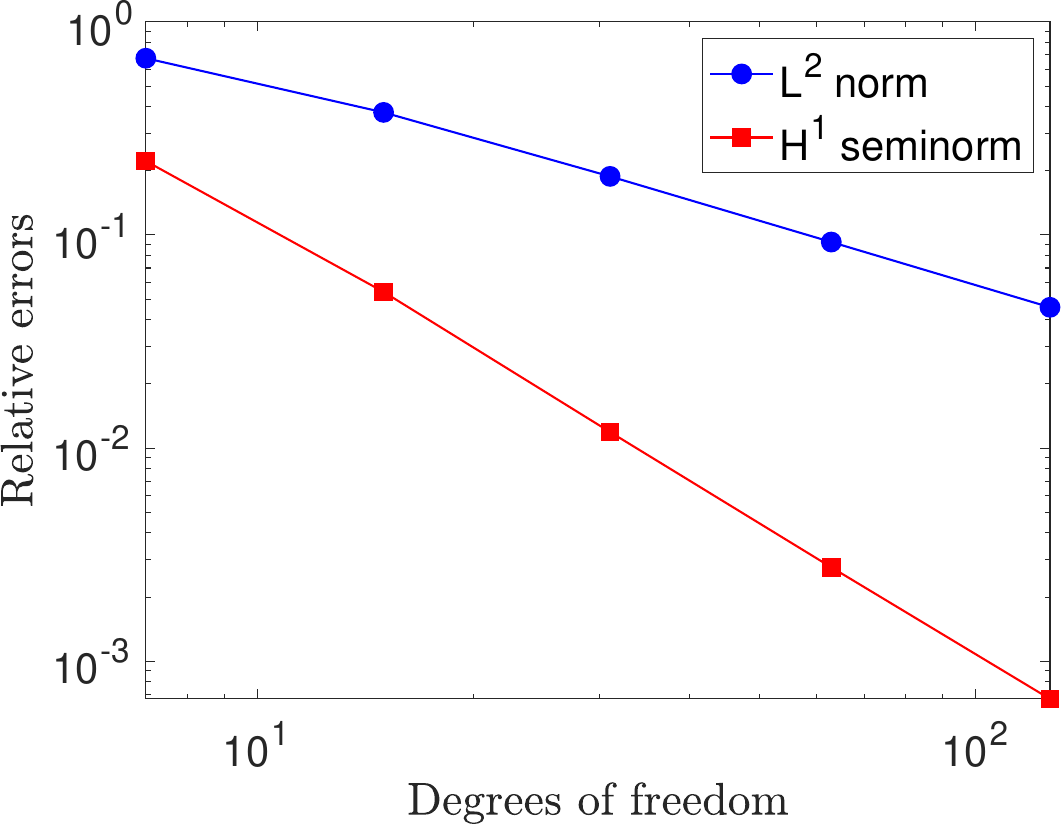}
\subcaption{}\label{fig:cd_alpha10_convergence_b}
\end{subfigure} \hfill
\begin{subfigure}{0.48\textwidth}
\includegraphics[width=\textwidth]{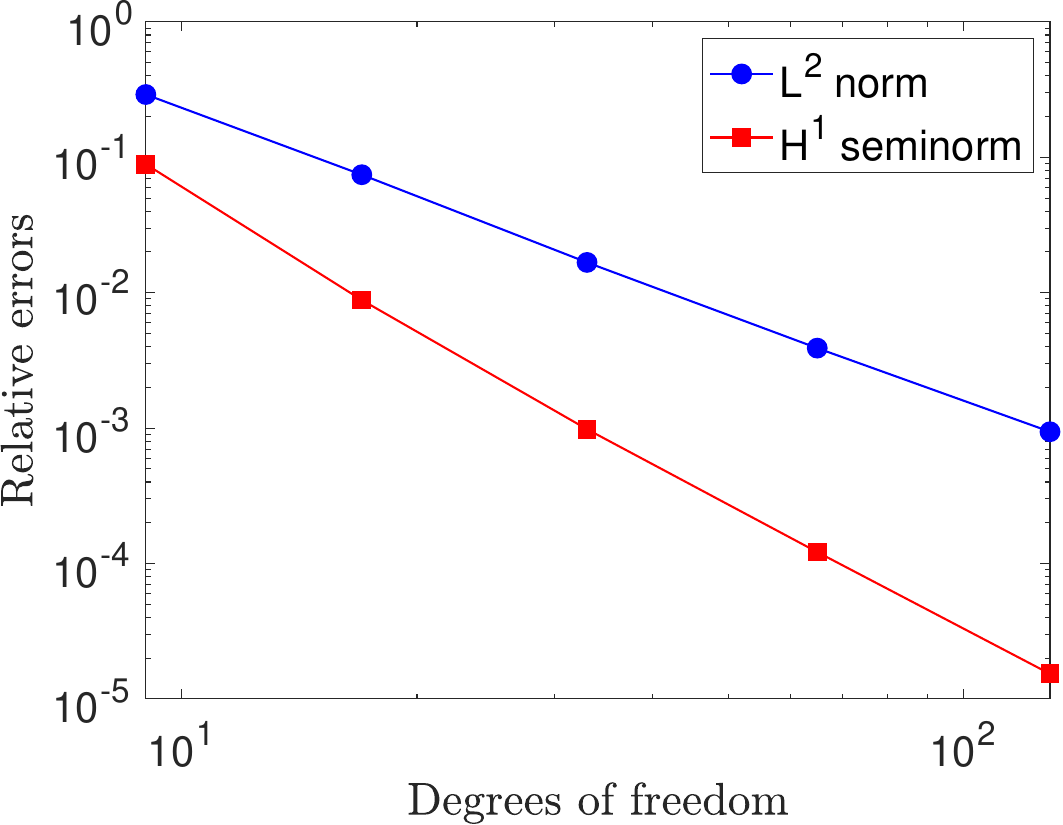}
\subcaption{}\label{fig:cd_alpha10_convergence_c}
\end{subfigure} \hfill
\begin{subfigure}{0.48\textwidth}
\includegraphics[width=\textwidth]{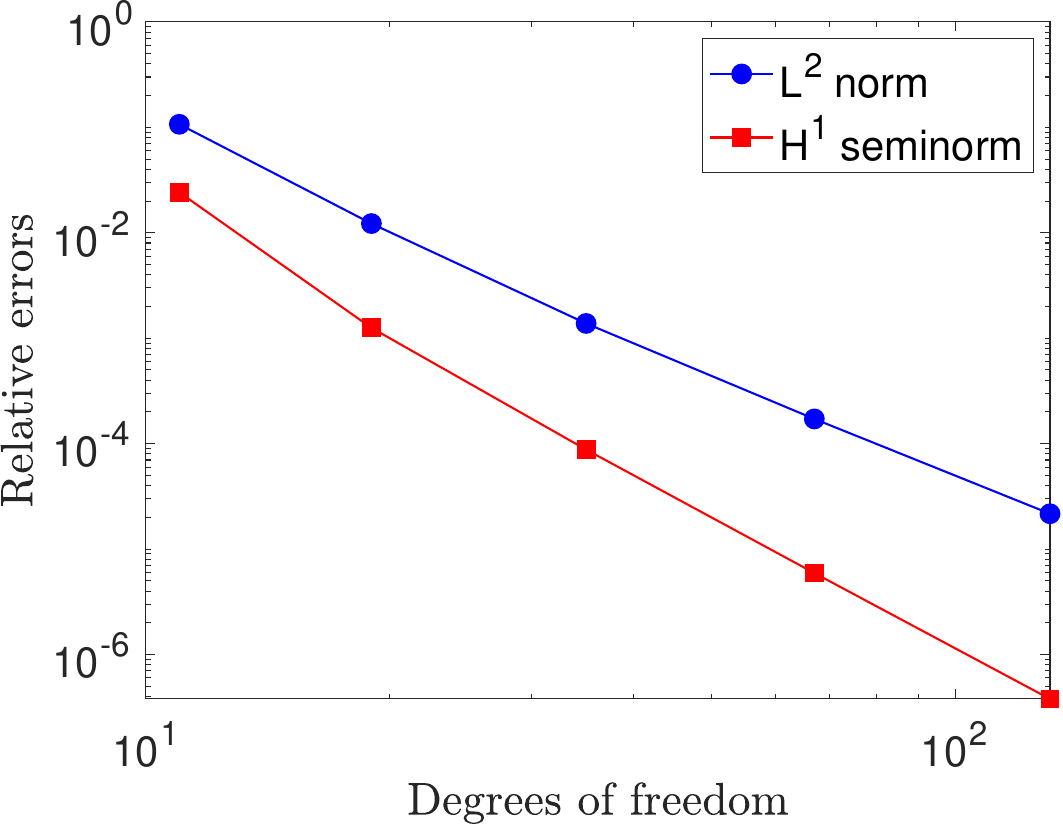}
\subcaption{}\label{fig:cd_alpha10_convergence_d}
\end{subfigure} 
\caption{Convergence study
         with B-splines for the steady-state
         convection-diffusion problem ($\alpha = 10$). 
         The dual fields $\mu_\theta(x)$ and $\lambda_\theta(x)$ are
         approximated using polynomials of degree 
         (a) $p = 1$ and $q = 1$; 
         (b) $p = 1$ and $q = 2$;
         (c) $p = 2$ and $q = 3$; and
         (d) $p = 3$ and $q = 4$. 
         The rate of convergence in $u$ and $u^\prime$
         is $p$ in (a), and are $p$ and $p+1$ in (b)--(d).}
         \label{fig:cd_alpha10_convergence}
\end{figure}
% Convergence study for alpha = 50
\begin{figure}[!tbh]
\centering
\begin{subfigure}{0.48\textwidth}
\includegraphics[width=\textwidth]{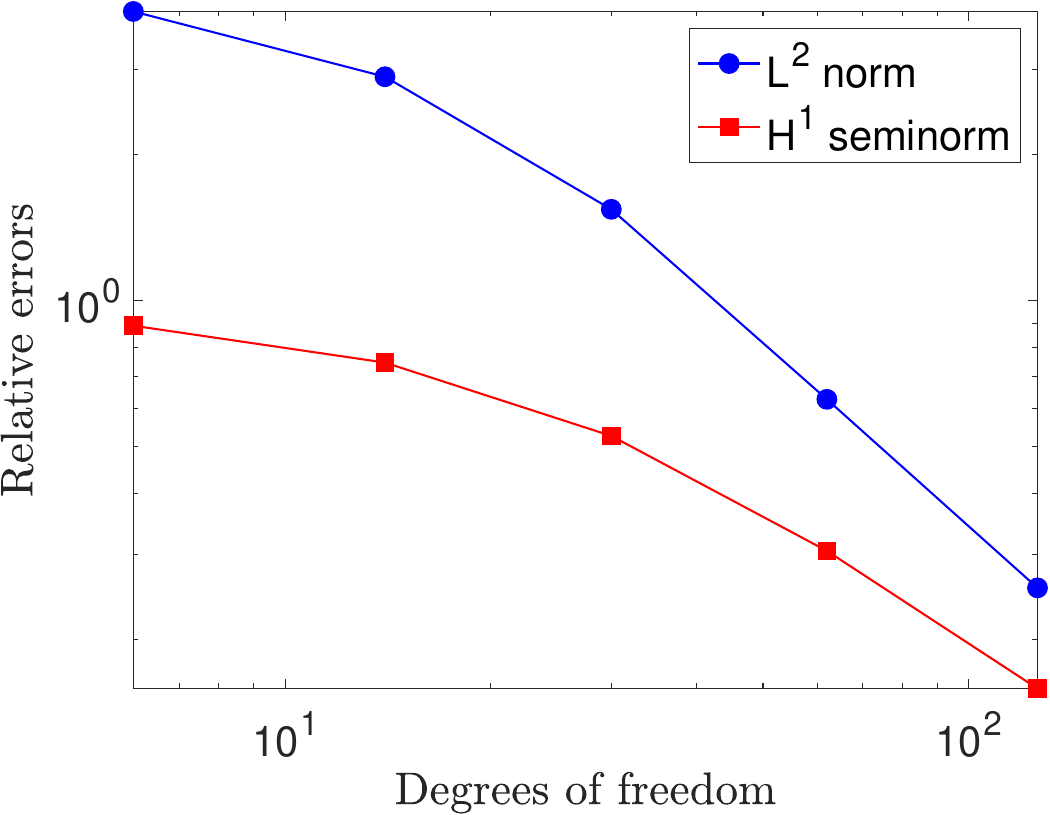}
\subcaption{}\label{fig:cd_alpha50_convergence_a}
\end{subfigure} \hfill
\begin{subfigure}{0.48\textwidth}
\includegraphics[width=\textwidth]{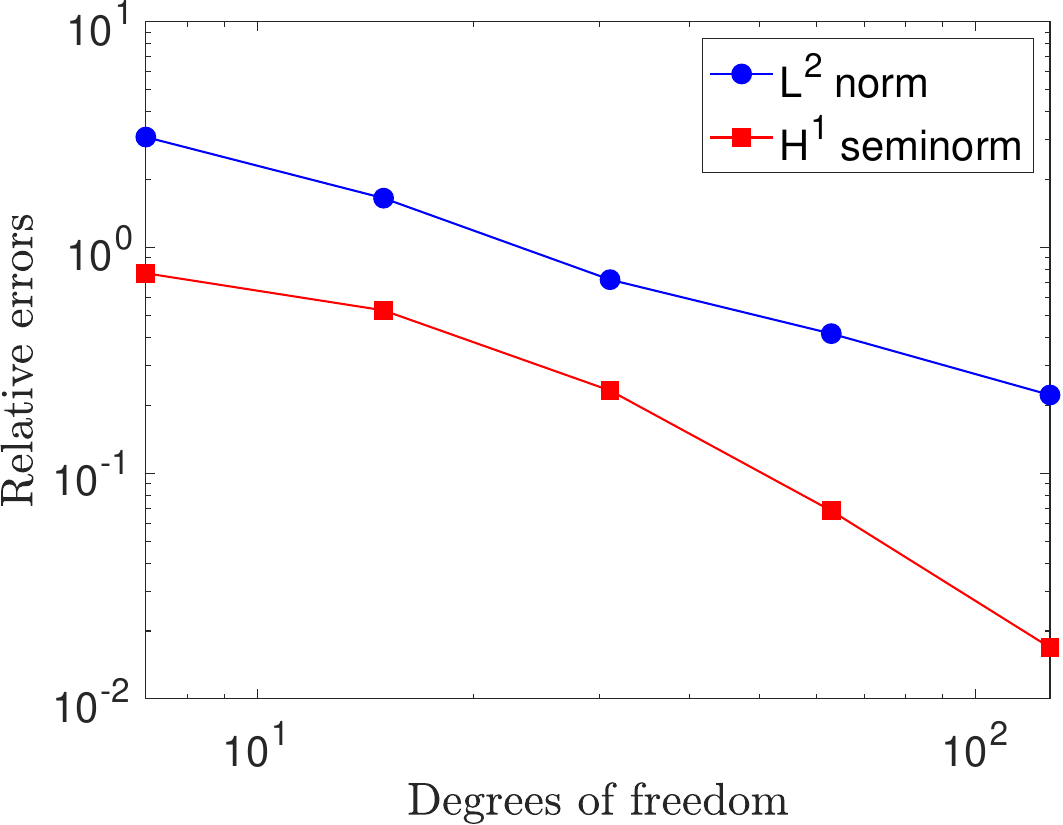}
\subcaption{}\label{fig:cd_alpha50_convergence_b}
\end{subfigure} \hfill
\begin{subfigure}{0.48\textwidth}
\includegraphics[width=\textwidth]{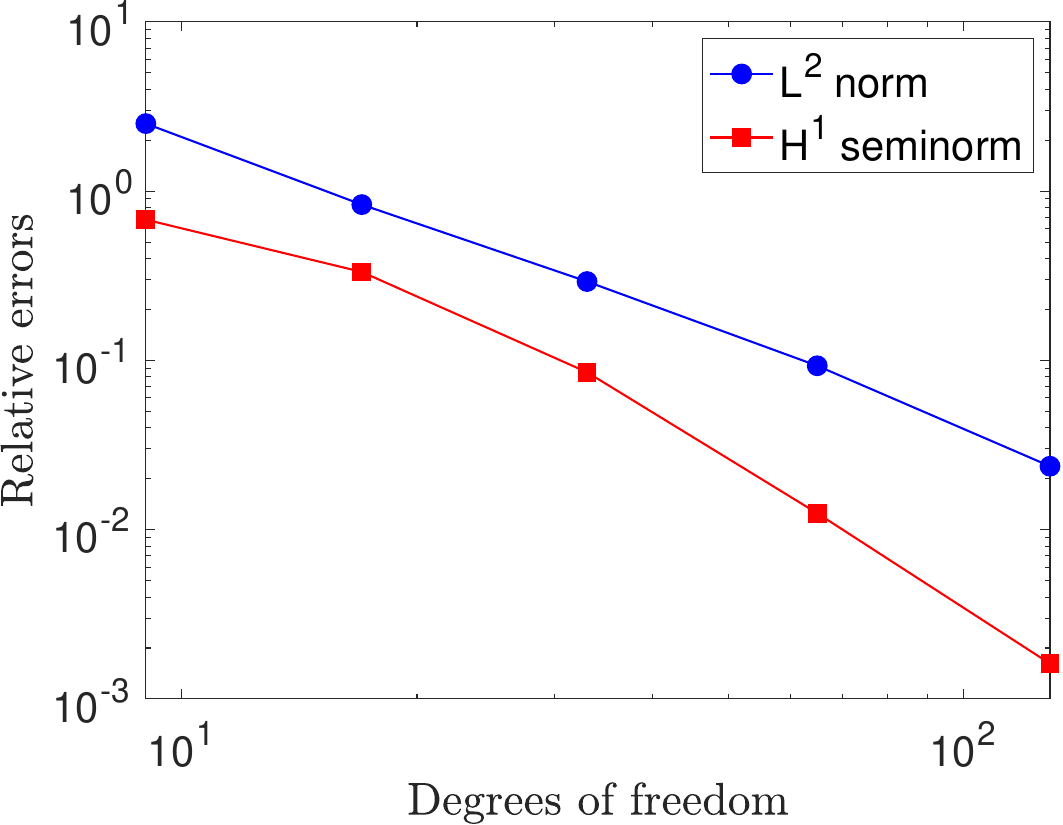}
\subcaption{}\label{fig:cd_alpha50_convergence_c}
\end{subfigure} \hfill
\begin{subfigure}{0.48\textwidth}
\includegraphics[width=\textwidth]{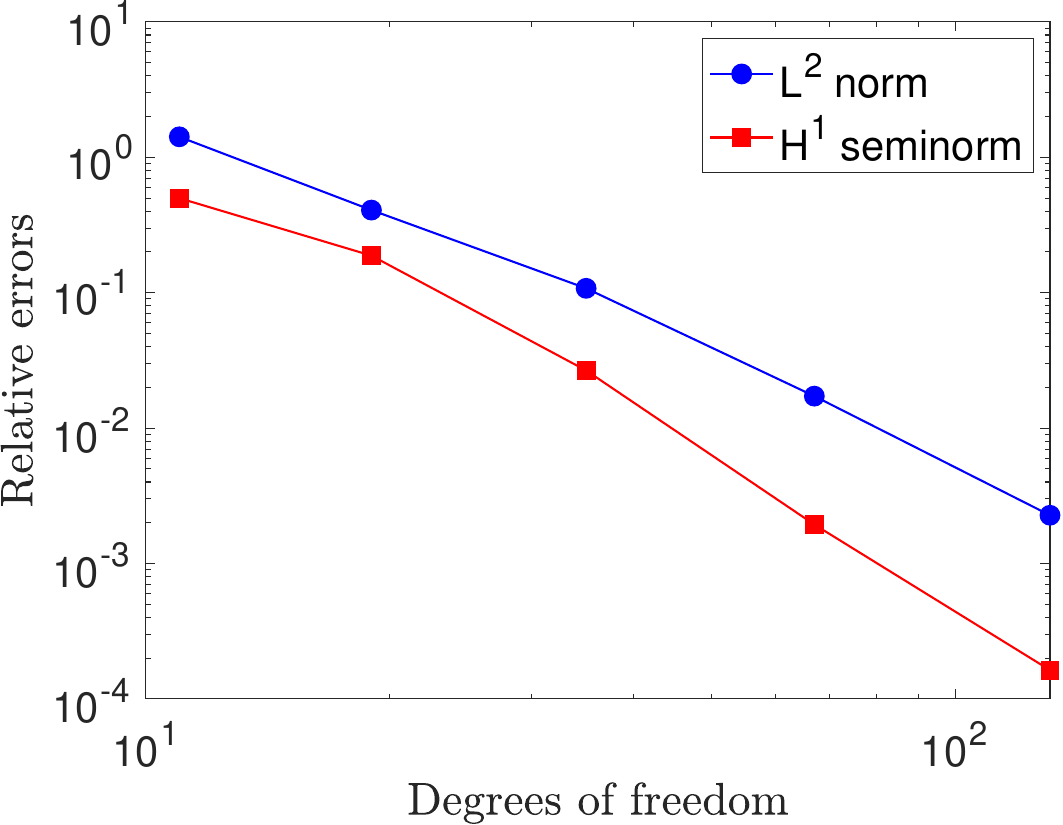}
\subcaption{}\label{fig:cd_alpha50_convergence_d}
\end{subfigure}
\caption{Convergence study
         with B-splines for the steady-state
         convection-diffusion problem ($\alpha = 50$). 
         The dual
         fields $\mu_\theta(x)$ and $\lambda_\theta(x)$ are
         approximated using polynomials of degree 
         (a) $p = 1$ and $q = 1$; 
         (b) $p = 1$ and $q = 2$;
         (c) $p = 2$ and $q = 3$; and
         (d) $p = 3$ and $q = 4$. 
         The rate of convergence in $u$ and $u^\prime$
         is $p$ in (a), and are $p$ and $p+1$ in (b)--(d).}
         \label{fig:cd_alpha50_convergence}
\end{figure}

\subsection{Transient convection-diffusion equation with B-splines}
\label{subsec:transient_cd_solution}
Consider the following transient convection-diffusion model problem:
\begin{subequations}\label{eq:ibvp_cd}
\begin{align}
\kappa \dfrac{\partial^2 u}{\partial x^2} - \alpha
\dfrac{\partial u}{\partial x} &= \dfrac{\partial u}{\partial t}
 \ \ \textrm{in } \Omega = \Omega_0 \times \Omega_t = (0,1) \times (0,1), 
    \label{eq:ibvp_cd_a} \\
u(0,t) &= 0, \quad u(1,t) = 0 ,
\label{eq:ibvp_cd_b} \\
u(x,0) &= u_0(x), 
\label{eq:ibvp_cd_c}
\intertext{with the exact solution~\cite{Cox:2008:ASD}:}
\begin{split}\label{eq:ibvp_cd_d}
u(x,t) &= \exp \left(- \frac{\alpha^2}{4 \kappa } t \right) \,
\exp \left(\frac{\alpha}{2 \kappa } x \right) \,
\sum_{n=1}^\infty
b_n \sin ( n \pi x ) 
\exp \left( - \kappa n^2 \pi^2 t \right) , \\
b_n &= 2 \int_0^1 \exp \left(- \frac{\alpha}{2 \kappa } x \right)
      \, u_0(x) \, \sin ( n \pi x) \, dx .
\end{split}
\end{align}
\end{subequations}
The ansatz for $\mu_\theta(x,t)$ and $\lambda_\theta(x,t)$ are formed by
tensor products of univariate B-splines given 
in~\eref{eq:bsp_mulambda}. The bivariate open knot vector is given by
\begin{equation}\label{eq:knot_2D}
\Xi \times \Xi, \quad \Xi := [\underset{p+1}{\underbrace{0,0,\dots,0}}, x_1, x_2, \dots, x_{n-1},
        \underset{p+1}{\underbrace{1,1,\dots,1}}].
\end{equation}
where $x_0 = 0$ and $x_n = 1$, and the trial functions are
\begin{equation}\label{eq:bsp_2D_mulambda}
  \mu_\theta(x,t) = \sum_{i=1}^{n+p} \sum_{j=1}^{n+p}
                       B_i^p(x) \, B_j^p(t) \, a_{ij}
                              , \quad 
  \lambda_\theta(x,t) = \sum_{i=1}^{n+q} \sum_{j=1}^{n+q}
                       B_i^q(x) \, B_j^q(t) \, b_{ij} .
\end{equation}
Appropriate coefficients are set to zero so that the
Dirichlet boundary conditions on $\lambda$ and $\mu$
are met.  On setting $\bar{u}_1 = 0$ and $\bar{u}_2 = 0$ to
meet the Dirichlet boundary conditions in~\eref{eq:ibvp_cd_b}
and choosing $u_0(x) = \sin(2\pi x)$ for the initial
condition in~\eref{eq:ibvp_cd_c},
the expression for the stiffness matrix and force vector
are obtained from~\eref{eq:Kd=f}. 

In the numerical computations, we choose $\kappa = 0.01$ and
$\alpha = 0.1$, so that advection dominates diffusion. Since
we performed a single solve over the entire 
space-time domain, $\Omega = (0,1)^2$, we judiciously selected values for these constants so that $u$ does not vary sharply in time and can be captured to modest accuracy by a high-order B-spline approximation. The exact solution in~\eref{eq:ibvp_cd_d} is computed with 1000 terms 
$(n = 1,2,\dots,1000)$, which satisfies the
initial condition in~\eref{eq:ibvp_cd_c}
to ${\cal O}(10^{-7})$ accuracy.  

The dual fields are constructed using bivariate
tensor-product
B-spline approximations.
Plots of the two-dimensional basis 
functions for $\mu_\theta(x,t)$ (degree $p = 9$; 100 basis
functions)
and $\lambda_\theta(x,t)$ (degree $q = 10$; 90 basis functions) 
are presented in~\fref{fig:bsp_2D_mulambda_cd}. 
Observe that
$\lambda_\theta(0,t) = \lambda_\theta(x,1) = 
\lambda_\theta(x,1) = 0$ so
that $\lambda_\theta \in S_\lambda$ and $\mu_\theta \in S_\mu$
satisfy~\eref{eq:dual_cd_c}, and are therefore kinematically
admissible.
\begin{figure}
\centering
\begin{subfigure}{0.48\textwidth}
\includegraphics[width=\textwidth]{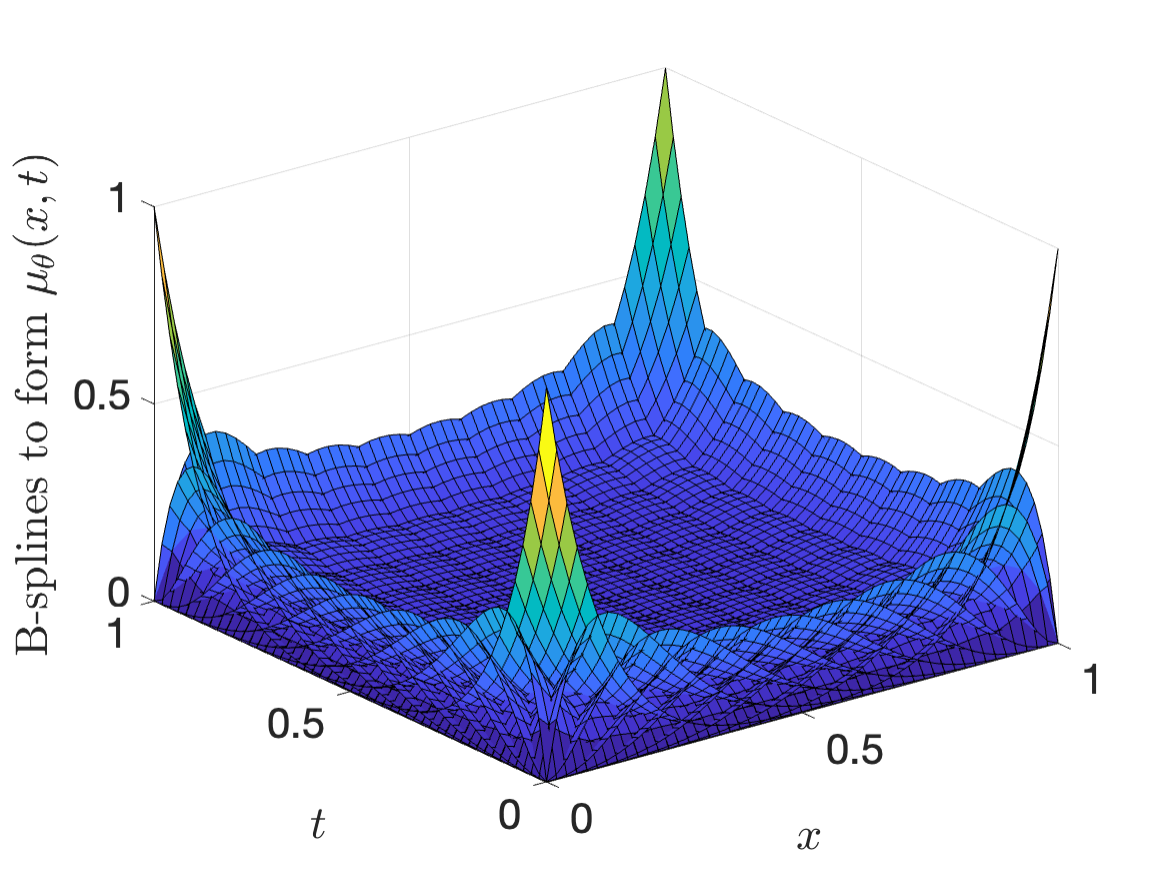}
\subcaption{}\label{fig:bsp_2D_mulambda_cd_a}
\end{subfigure} \hfill
\begin{subfigure}{0.48\textwidth}
\includegraphics[width=\textwidth]{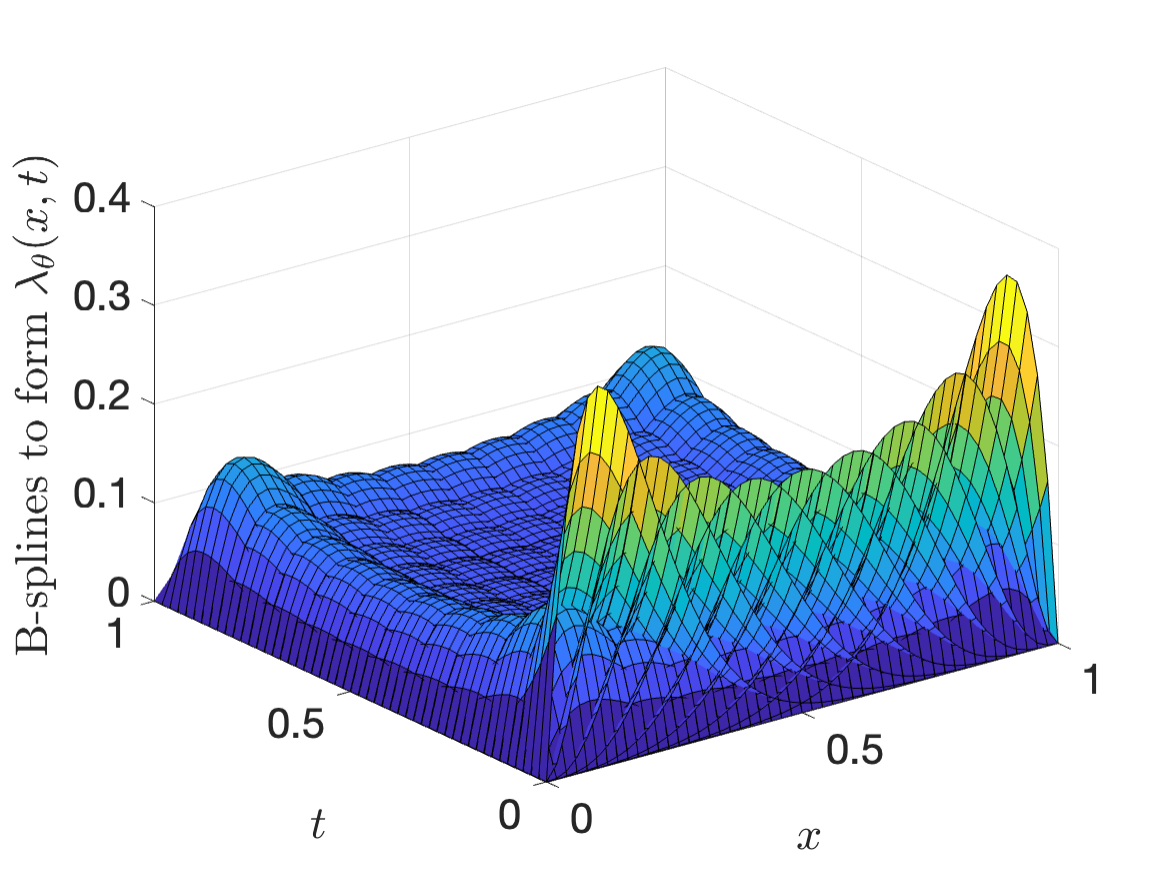}
\subcaption{}\label{fig:bsp_2D_mulambda_cd_b}
\end{subfigure}
\caption{Plots of two-dimensional B-spline basis functions
         to form the dual fields (a) $\mu_\theta(x,t)$ ($p = 9$) and 
         (b) $\lambda_\theta(x,t)$ ($q = 10$) to solve the 
         transient convection-diffusion problem. 
         The bivariate open knot vector
         ($n = 1$) that is given
         in~\eref{eq:knot_2D} is used.}
         \label{fig:bsp_2D_mulambda_cd}
\end{figure}
The B-spline solution for the dual fields, 
$\mu_\theta(x,t)$ and
$\lambda_\theta(x,t)$, are presented 
in~\fref{fig:mulambda_cd}. 
The B-spline solution for $u$ and $q := \partial u / \partial x$ are computed from the dual fields using
the DtP map given in~\eqref{eq:dtp_cd}. These primal
B-spline solutions 
are compared to the exact solutions in~\fref{fig:bsp_transient_cd}. We find that the
B-spline solutions for $u$ and $q$ are in fair agreement with the
exact solutions. Since $||u||_\infty = 1$ and
$|\partial u / \partial x |_\infty \approx 10$, the
relative maximum error in $u$ and $q$ are 0.06 and 0.1, respectively
(see Figs.~\ref{fig:bsp_transient_cd_c} and~\ref{fig:bsp_transient_cd_f}).
In addition, we observe that the errors grow larger
close to the terminal time, $t = 1$. This is
evident from
the time history plots shown in Figs.~\ref{fig:bsp_transient_cd_g} to~\ref{fig:bsp_transient_cd_j}. Notably the maximum error
in $u$ for $t \in [0,1]$ (see~\fref{fig:bsp_transient_cd_g})
is about six times that for
$t \in [0,0.9]$ (see~\fref{fig:bsp_transient_cd_h}).
The same trend is also noticed for the error in $q$ from
Figs.~\ref{fig:bsp_transient_cd_i} and~\ref{fig:bsp_transient_cd_j}. 
We attribute the larger errors that concentrate in the
vicinity of $t = 1$ (terminal time)
to the interactions between the boundary conditions in the primal problem and the terminal boundary condition in the dual problem, 
which may be understood via the explanation that follows.

\subsubsection{Understanding the solution behavior near the terminal time}
\label{subsubsec:T_issues_cd}
Consider the DtP mapping equations \eqref{eq:dtp_cd_b}, which are reproduced below:
\begin{align}\label{eq:dtp_cd_numerics}
\begin{split}
u_H &= \dfrac{\partial \lambda}{\partial t} 
                        + \dfrac{\partial \mu}{\partial x}
                        =: \partial_t \lambda + \partial_x \mu, \\
q_H &= \mu - \alpha \lambda - \kappa \dfrac{\partial \lambda}{\partial x} 
    =: \mu - \alpha \lambda - \kappa \,\partial_x \lambda .
\end{split}
\end{align} 
Let $T$ be the terminal time ($T = 1$ herein) and recall that $q(x,t) = \partial_x u(x,t)$. 
Now, for the sake of this argument, assume 
the functions
$u(x,T), \, q(x,T)$ for $x \in [0,1]$ 
are known from the unique solution to the primal problem with specified initial and boundary conditions, which the discrete solution
aspires to reproduce.
Given the prescribed
Dirichlet boundary condition on the top edge 
of the space-time domain, say $\lambda(x,T) = \lambda_\textrm{top}(x)$,
it is clear from the expression for $q_H$ in~\eref{eq:dtp_cd_numerics} that 
$\mu(x,T)$ for $x \in [0,1]$ on the 
top edge become fixed, which implies that $\partial_x \mu (x,T)$ is also fixed. From the 
expression for $u_H$ in~\eref{eq:dtp_cd_numerics}, this 
also means that $\partial_t \lambda(x,T)$ is fixed on the top edge, and 
in particular $\partial_t\lambda (1,1) =: d_{\textrm{top}}$. 
However, there is a Dirichlet boundary condition that is specified on the 
right edge of the domain, say $\lambda(1,t) = \lambda_{\textrm{right}}(t)$, thus fixing $\partial_t \lambda (1,\cdot)$, and in particular $\partial_t \lambda(1,T) =: d_{\textrm{right}}$. In our computations, 
we chose the functions 
$\lambda_{\textrm{top}} = 0, \lambda_{\textrm{right}} = 0$,
but note that the value $d_{\textrm{top}}$ is not entirely determined from the function $\lambda_{\textrm{top}}$. 
Now, if $u$ is continuous at $(1,1)$ but $d_{\textrm{top}}
\neq d_{\textrm{right}}$, values that depend on the arbitrarily 
specified Dirichlet specifications of $\lambda_{\textrm{top}}$ and $\lambda_{\textrm{right}}$, then approximating such a dual solution (this is not an impediment for a weak formulation) with smooth interpolation can result in 
higher errors at the $(1,1)$ corner of the domain. On
applying similar
arguments to the behavior of the solutions on the left and top edges lead us to also draw the same inference at
the $(0,1)$ corner. 
Another way to interpret this result is that if $(u_H, q_H)$ have to 
equal 
the exact solution, say $(u, q)$, of the primal problem on the top edge, then it is possible that in general $d_{\textrm{top}}$ may not equal $d_{\textrm{right}}$, and similarly for the left top corner. This implies that if such discontinuities are not allowed in the dual solution, then either $u_H \neq u$ or $q_H \neq q$ (or both) on the
top edge. We note that $\mu$ has no boundary conditions specified
for this problem~\eqref{eq:dual_cd_c}, and hence has substantially more freedom to accommodate the demands of continuity in the primal solution at the top corners of the domain.  In general, even if an outright discontinuity like 
$d_{\textrm{top}} \neq d_{\textrm{right}}$ does not occur, the demands of the dual Dirichlet boundary conditions and the primal problem to be solved can set up boundary layers in the dual 
solution (see the dual solutions shown in~\fref{fig:mulambda_cd}). 
These features were demonstrated 
for the heat equation and pure convection case in~\cite[Secs.~5.1.1, 5.2.1]{Kouskiya:2024:HCH} that led to 
some degradation in solution accuracy (as in our calculations here), but which can be alleviated with 
refinement. We point out that while the `corner' issue does not arise in the solution of ODEs by the dual methodology, a boundary layer in the dual solution near the final time can arise. However, these final-time issues can be robustly 
resolved, 
as explained in the following. 

The degradation of accuracy of the solution on the top edge is  
easily dealt with, as demonstrated in~\cite{Kouskiya:2024:HCH,Kouskiya:2024:IBD}. 
Let the primal initial-boundary value problem, posed
in the time interval $[0,T]$, have a unique solution (and even otherwise). Then, one poses and solves the dual problem in an interval $[0,T+ \delta]$ in time, with the boundary conditions of the primal (original) problem imposed in the interval $[0,T]$, and appending an (arbitrary) continuous extension of these primal boundary conditions in the interval $[T, T+ \delta]$. Recovering the primal fields from the dual solution, it can be seen from the logic of the scheme that, in principle, the correct primal problem would have been solved in the time interval $[0,T]$. 
Thus, approaching this augmented primal problem with a  `buffer-zone' in the time-like direction with the dual technique and discarding the solution in the space-time strip $[0,L] \times (T, T+ \delta]$ ($\delta$ preferably small),
eliminates all issues with dual solution accuracy at top corners/edges in the truncated solution. Furthermore, since the dual solution scheme involves solving boundary-value problems in space-time domains, solution costs can be prohibitive for 
long-time simulations.
For such situations, a time-slicing strategy can be utilized where the problem is solved in a finite number of stages in time, the truncated solution at the end of one stage serving as the primal initial condition for the successive stage; these procedures have been successfully demonstrated 
in~\cite{Kouskiya:2024:HCH,Kouskiya:2024:IBD}.
\begin{figure}[!tbh]
\centering
\begin{subfigure}{0.48\textwidth}
\includegraphics[width=\textwidth]{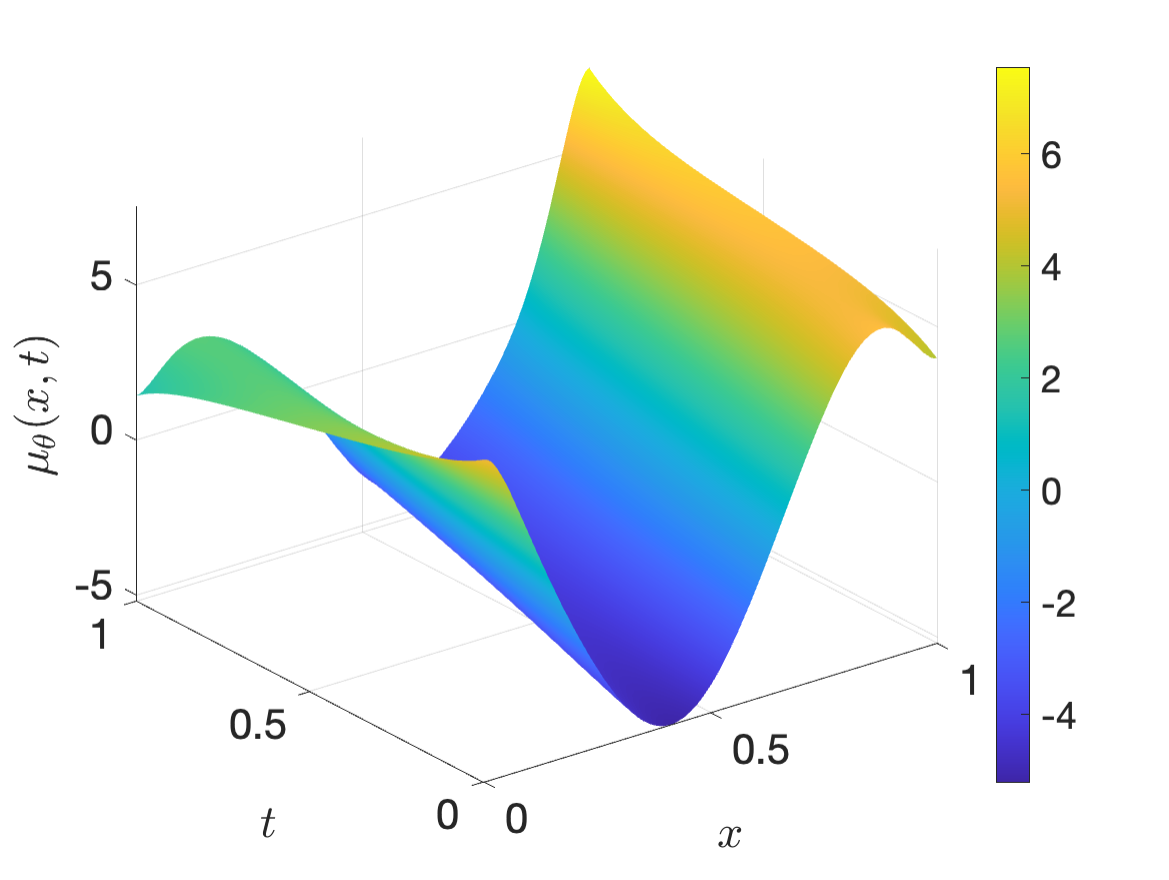}
\subcaption{}\label{fig:mulambda_cd_a}
\end{subfigure} \hfill
\begin{subfigure}{0.48\textwidth}
\includegraphics[width=\textwidth]{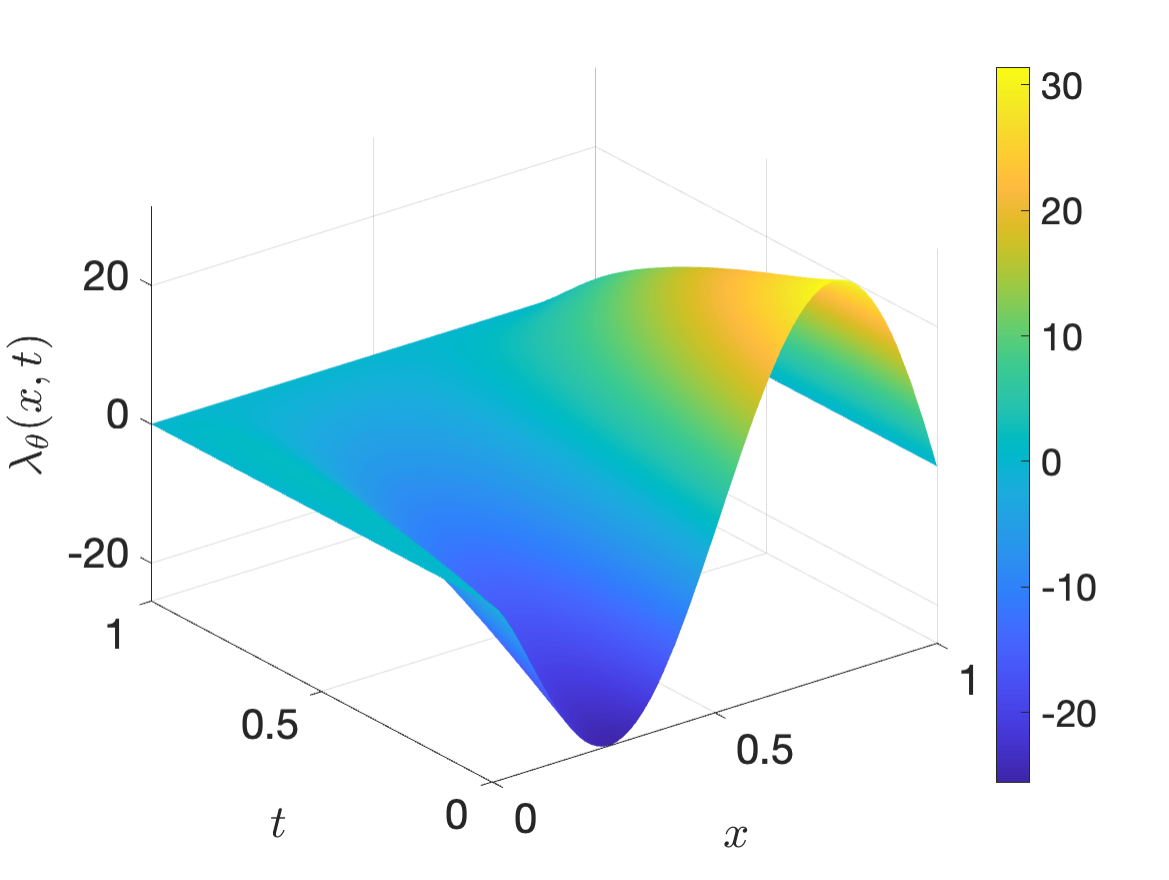}
\subcaption{}\label{fig:mulambda_cd_b}
\end{subfigure}
\caption{B-spline solution (dual fields)
         of the transient convection-diffusion
         equation ($\kappa = 0.01, \alpha = 0.1$).
         (a) $\mu_\theta(x,t)$ ($p = 9$) and 
         (b) $\lambda_\theta(x,t)$ ($q = 10$).
         }
         \label{fig:mulambda_cd}
\end{figure}
\begin{figure}[!tbh]
\centering
\begin{subfigure}{0.33\textwidth}
\includegraphics[width=\textwidth]{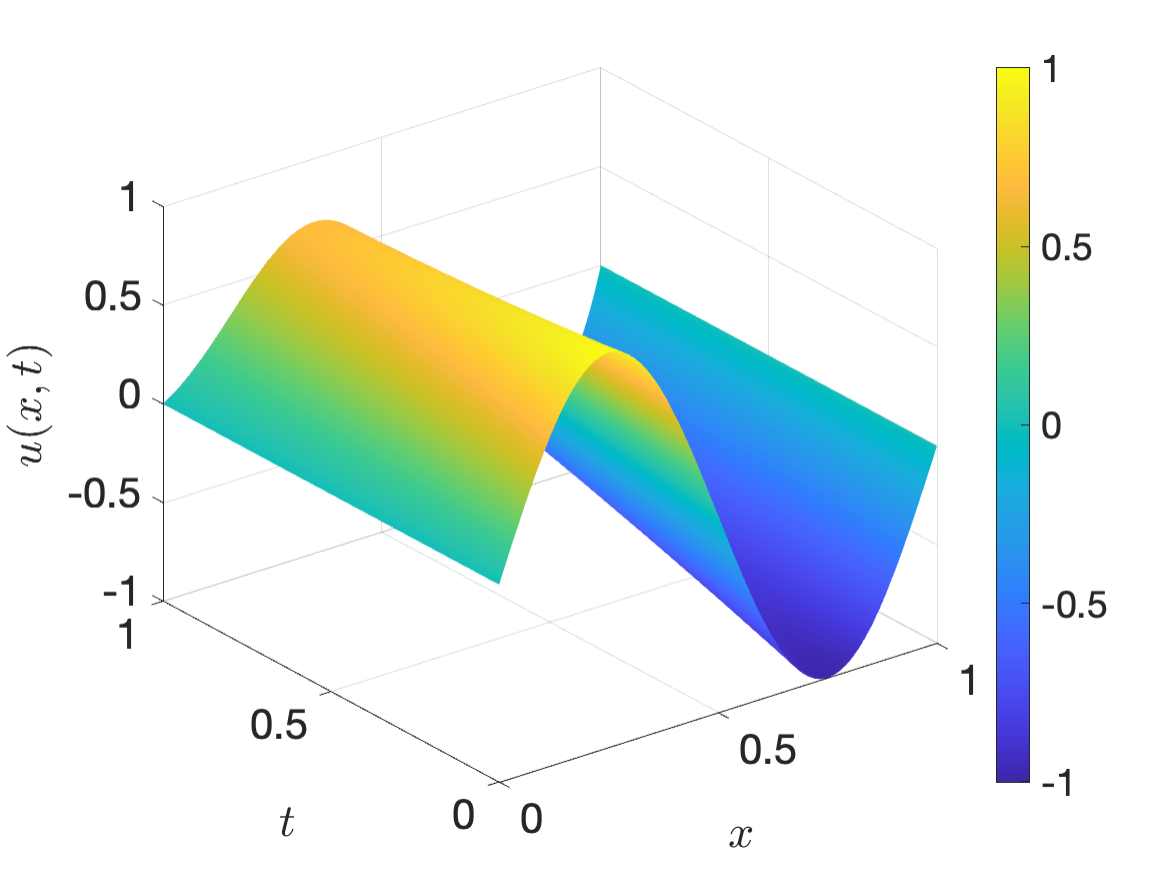}
\subcaption{}\label{fig:bsp_transient_cd_a}
\end{subfigure} \hfill
\begin{subfigure}{0.33\textwidth}
\includegraphics[width=\textwidth]{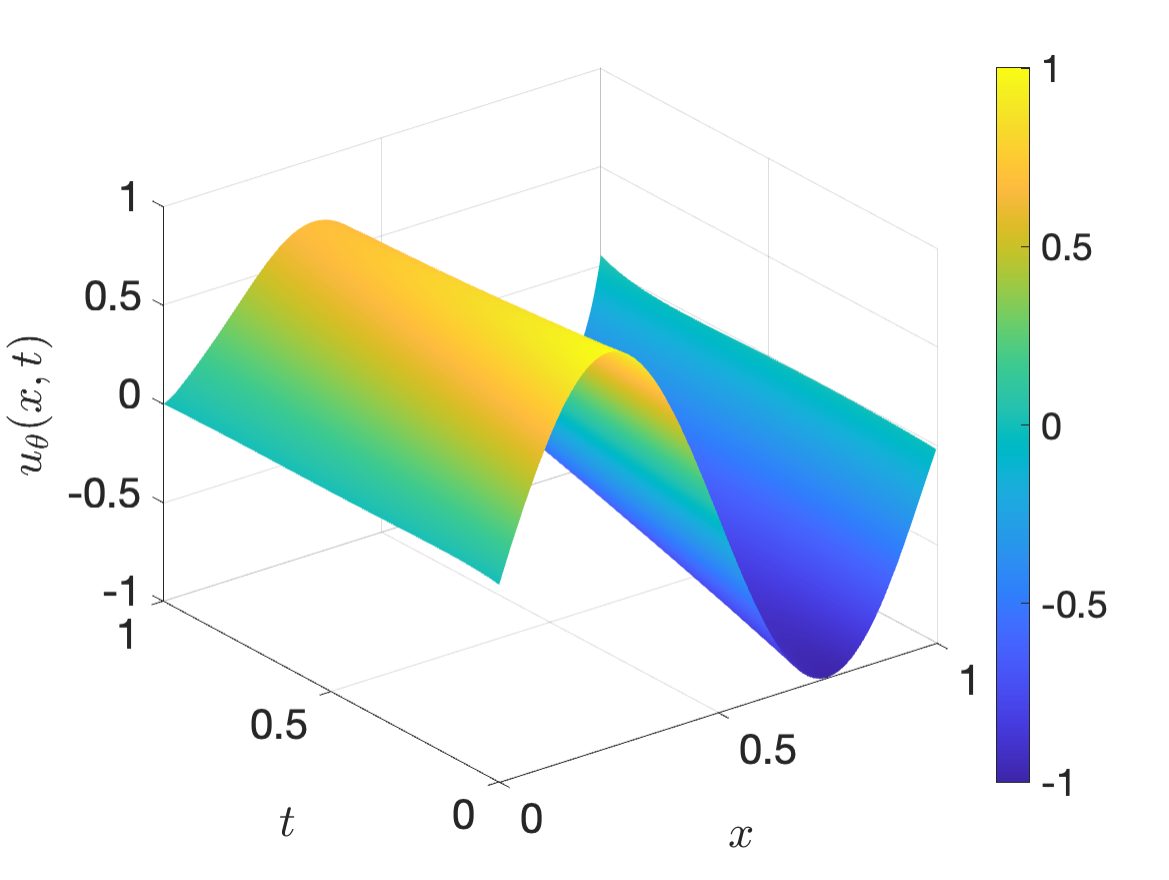}
\subcaption{}\label{fig:bsp_transient_cd_b}
\end{subfigure} \hfill
\begin{subfigure}{0.33\textwidth}
\includegraphics[width=\textwidth]{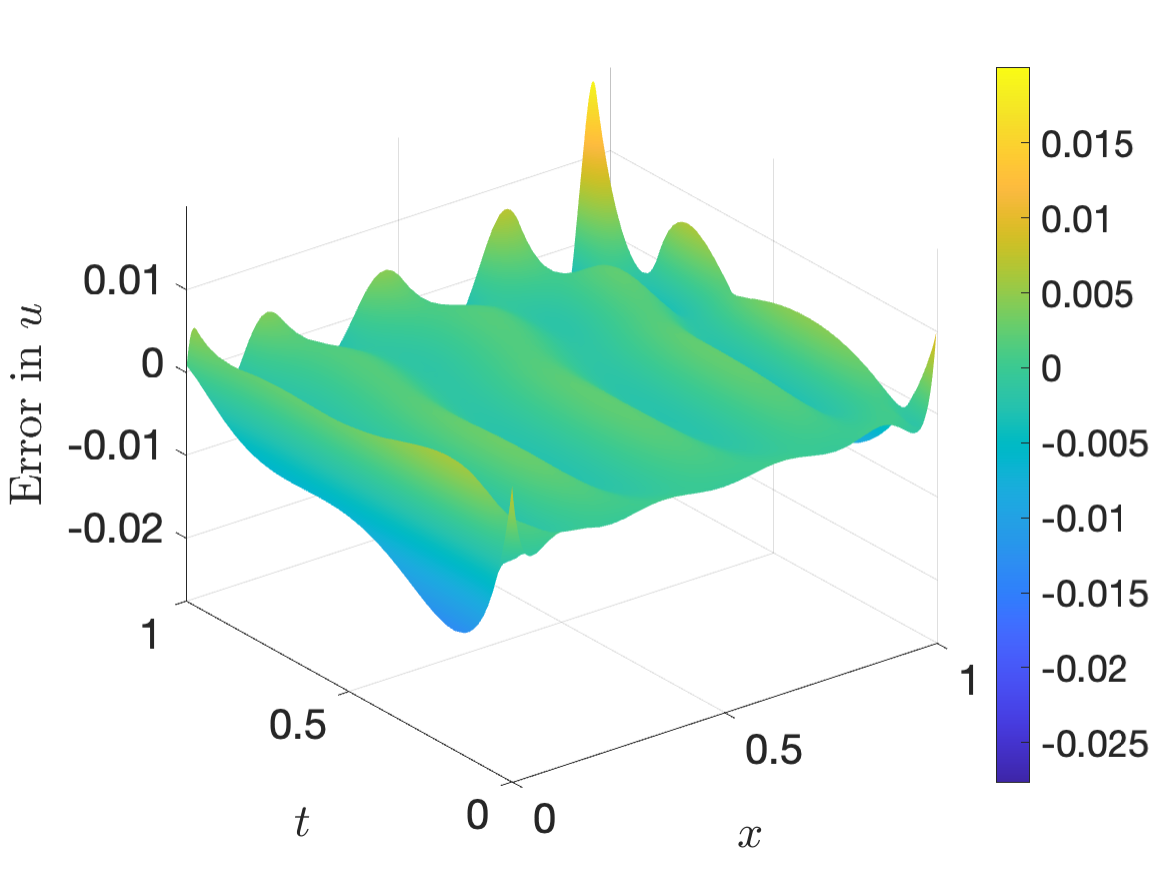}
\subcaption{}\label{fig:bsp_transient_cd_c}
\end{subfigure}
\begin{subfigure}{0.33\textwidth}
\includegraphics[width=\textwidth]{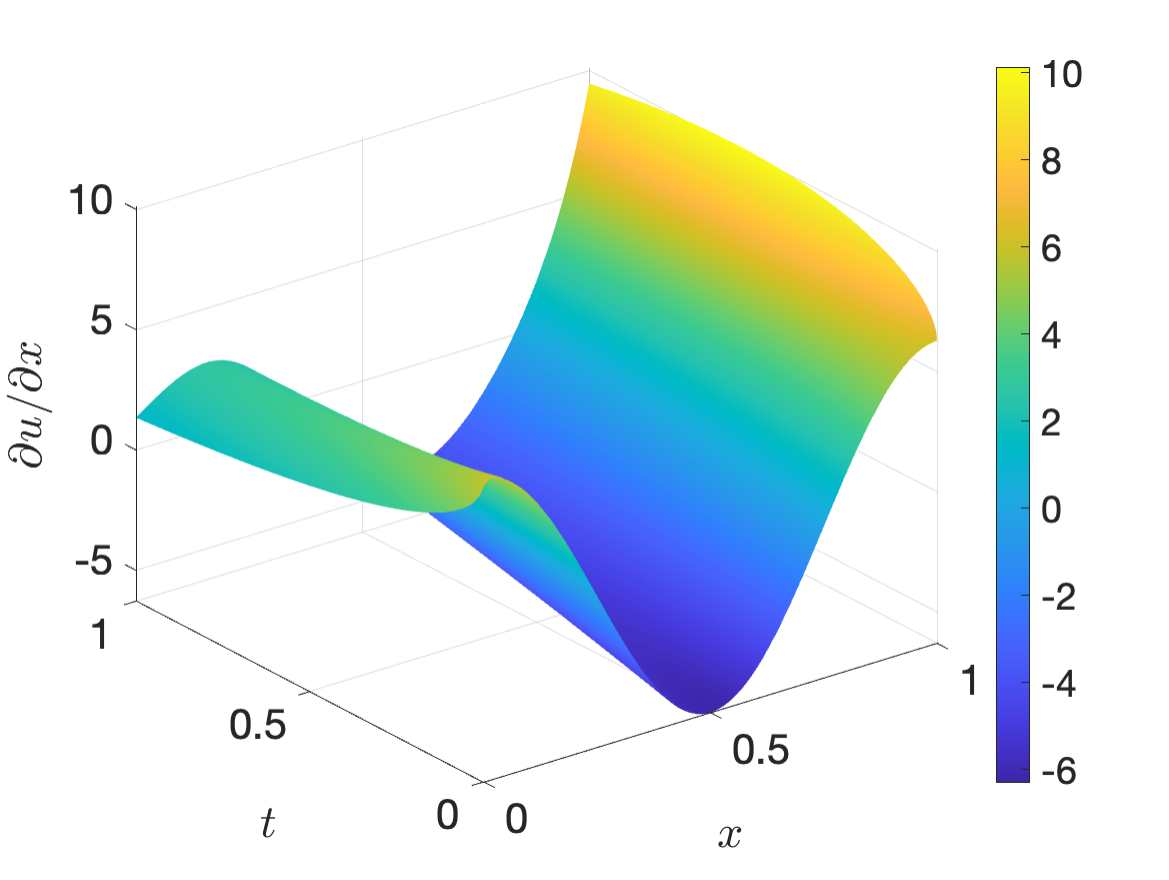}
\subcaption{}\label{fig:bsp_transient_cd_d}
\end{subfigure} \hfill
\begin{subfigure}{0.33\textwidth}
\includegraphics[width=\textwidth]{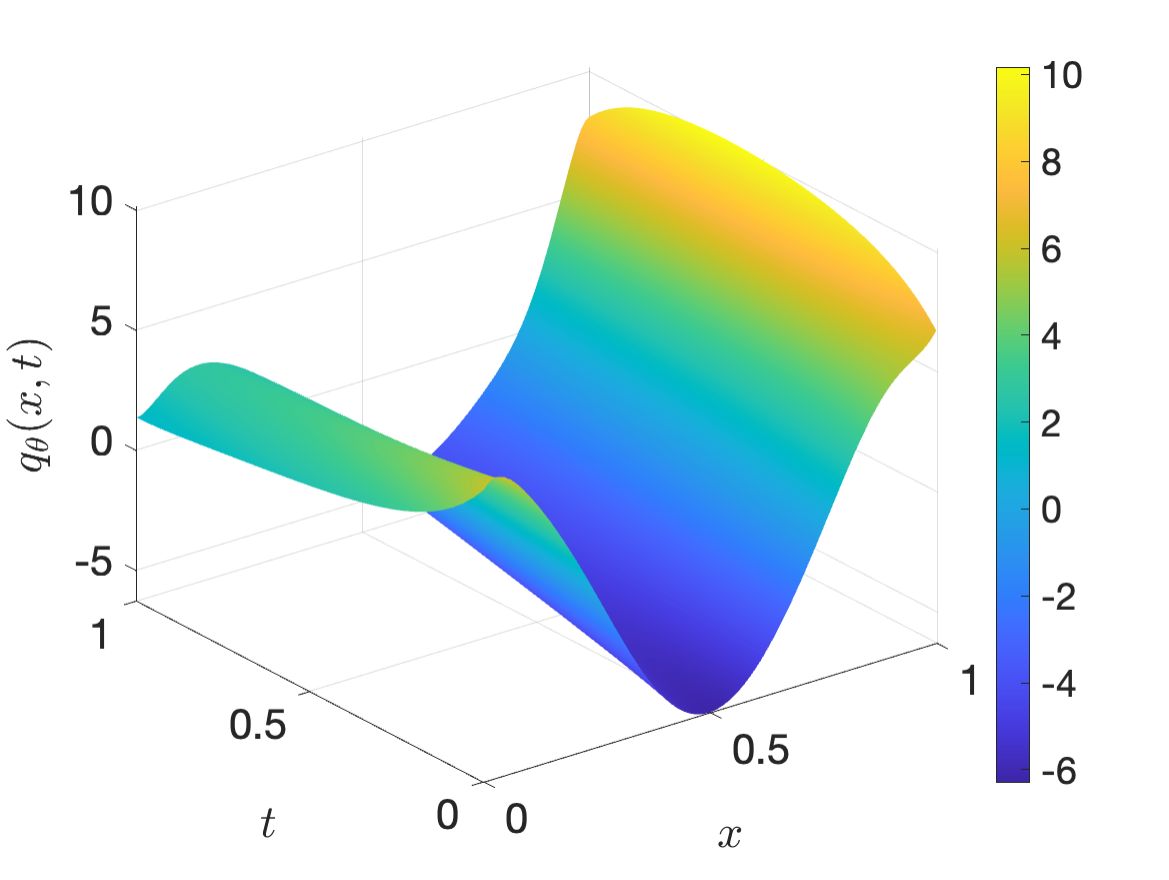}
\subcaption{}\label{fig:bsp_transient_cd_e}
\end{subfigure} \hfill
\begin{subfigure}{0.33\textwidth}
\includegraphics[width=\textwidth]{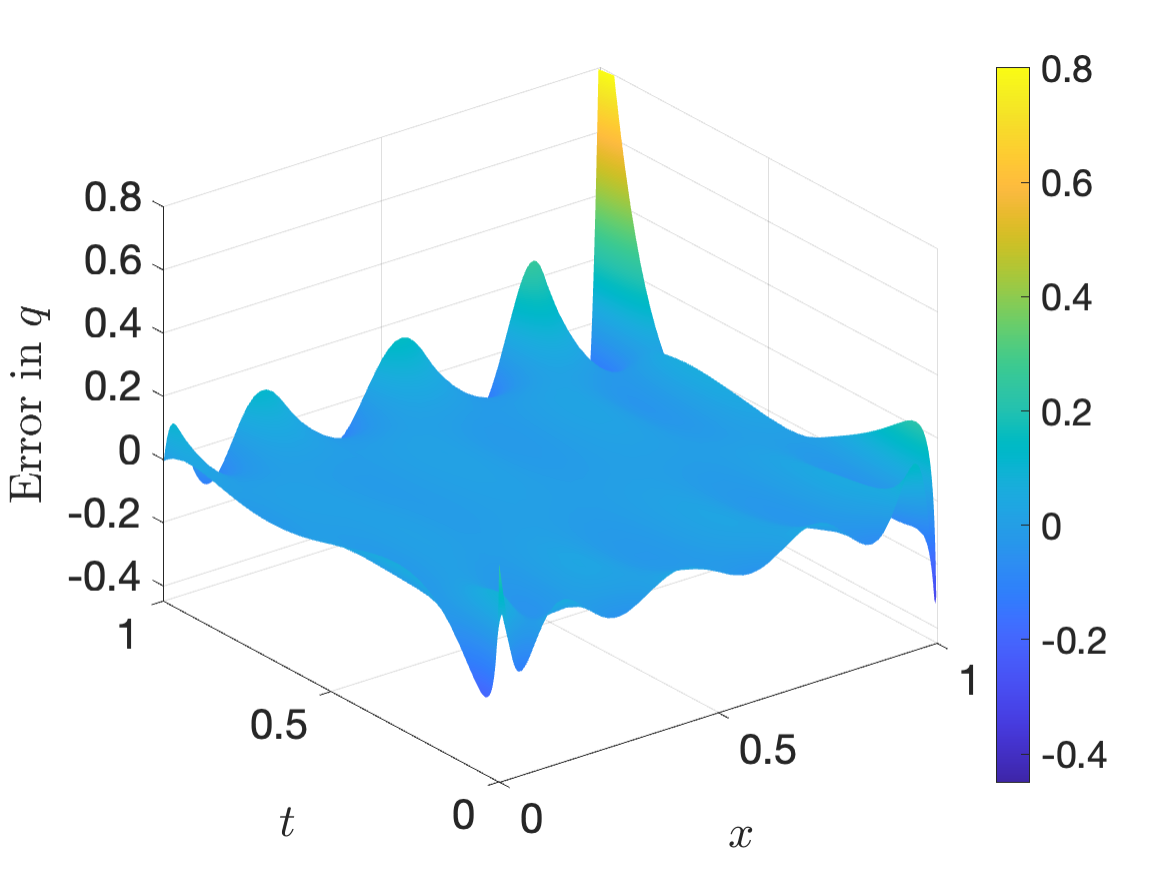}
\subcaption{}\label{fig:bsp_transient_cd_f}
\end{subfigure}
\begin{subfigure}{0.24\textwidth}
\includegraphics[width=\textwidth]{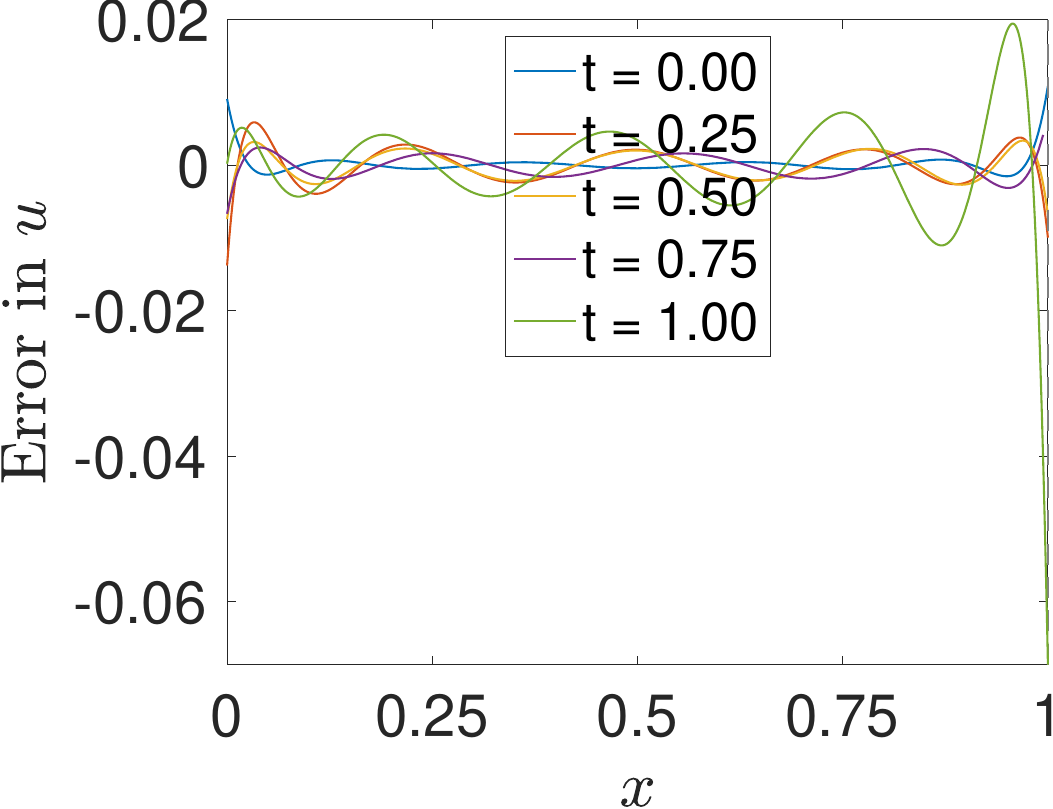}
\subcaption{}\label{fig:bsp_transient_cd_g}
\end{subfigure} \hfill
\begin{subfigure}{0.24\textwidth}
\includegraphics[width=\textwidth]{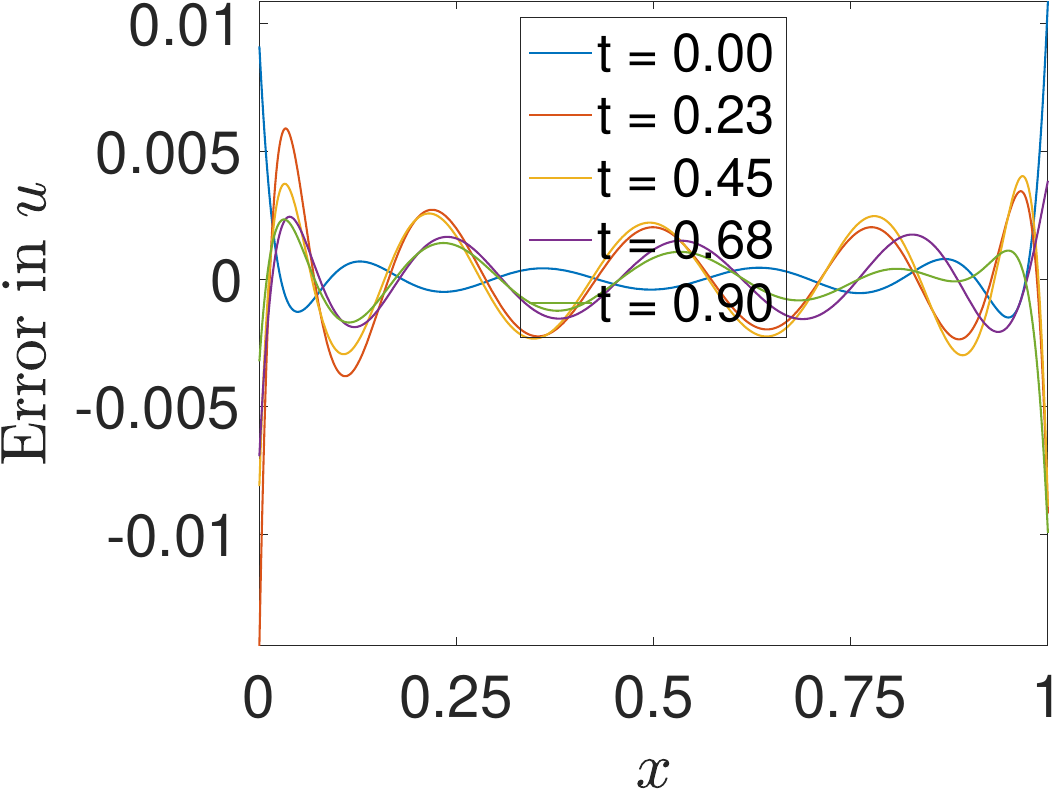}
\subcaption{}\label{fig:bsp_transient_cd_h}
\end{subfigure} \hfill
\begin{subfigure}{0.24\textwidth}
\includegraphics[width=\textwidth]{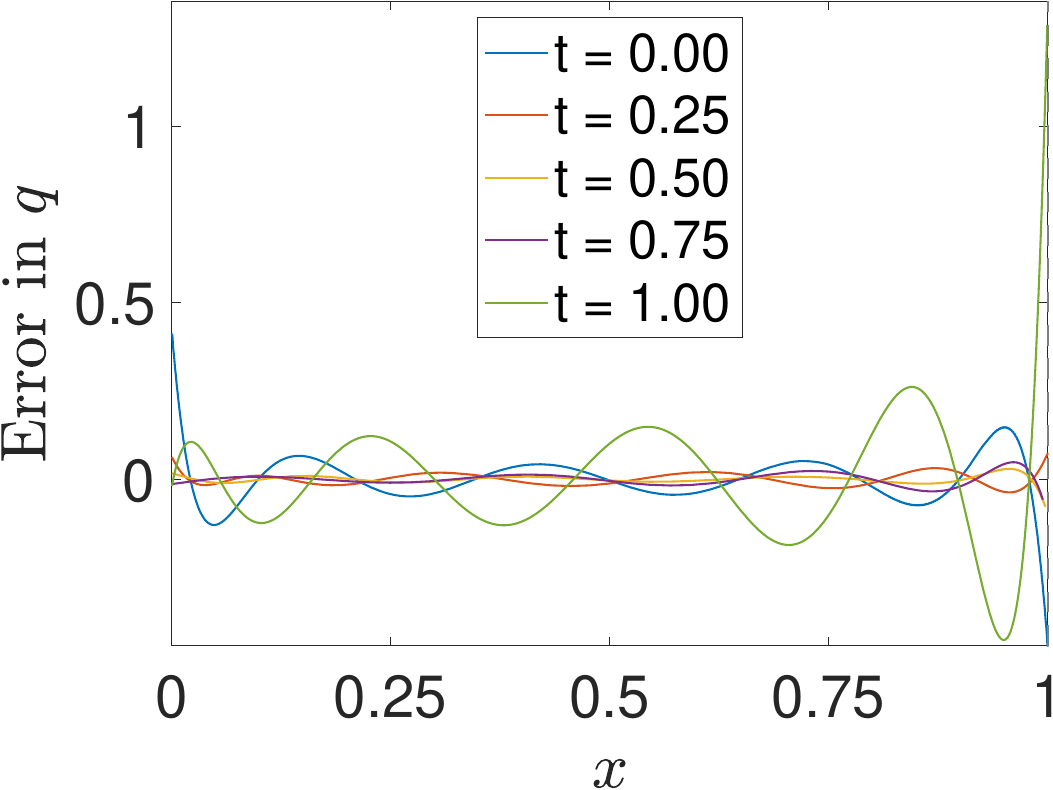}
\subcaption{}\label{fig:bsp_transient_cd_i}
\end{subfigure} \hfill
\begin{subfigure}{0.24\textwidth}
\includegraphics[width=\textwidth]{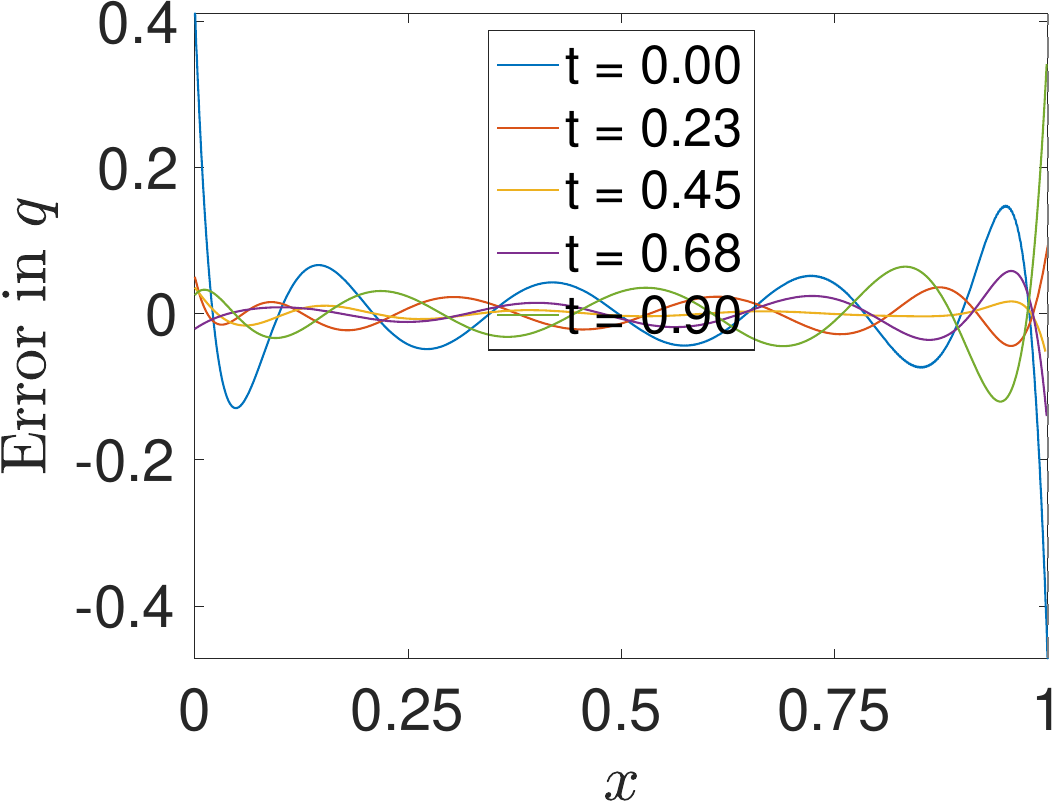}
\subcaption{}\label{fig:bsp_transient_cd_j}
\end{subfigure} 
\caption{Space-time B-splines to solve the
         transient convection-diffusion problem.
         The conductivity coefficient 
         $\kappa = 0.01$ and the convection coefficient
         $\alpha = 0.1$.
         The dual fields $\mu_\theta(x,t)$
         and $\lambda_\theta(x,t)$ are composed of tensor-product
         B-splines of degree $p = 9$ and $q = 10$, respectively.
         (a) Exact solution, $u$; (b) $u_\theta$; 
         (c) Error, $u - u_\theta$;
         (d) Exact $q = \partial u / \partial x$; 
         (e) $q_\theta$; and
         (f) Error, $q_\theta - \partial u 
         / \partial x$.
         Plots of the time history for the error in $u$ 
         at (g) $t \in [0,1]$ and (h) $t \in [0, 0.9]$.
         Plots of the time history for the error in $q$ 
         for (i) $t \in [0,1]$ and (j) $t \in [0, 0.9]$.}
         \label{fig:bsp_transient_cd}
\end{figure}

\subsection{Transient heat equation with 
B-splines}\label{subsec:heat_solution}
Consider the following transient heat conduction model problem:
\begin{subequations}\label{eq:ibvp_heat}
\begin{align}
\kappa \dfrac{\partial^2 u}{\partial x^2} &= \dfrac{\partial u}{\partial t}
 \ \ \textrm{in } \Omega = \Omega_0 \times \Omega_t = (0,1) \times (0,1), 
    \label{eq:ibvp_heat_a} \\
u(0,t) &= 1, \quad \kappa \frac{\partial u}{\partial x}(1,t) = 0 ,
\label{eq:ibvp_heat_b} \\
u(x,0) &= u_0(x) =  1 + \sin \left( \frac{\pi x}{2} \right),
\label{eq:ibvp_heat_c}
\intertext{with the exact solution:}
u(x,t) &= 1 + \sin \left( \frac{\pi x}{2} \right)
\exp \left(- \frac{\pi^2 \kappa }{4} t \right) .
\label{eq:ibvp_heat_d}
\end{align}
\end{subequations}

The dual functional for this IBVP with
$\bar{u}_1 = 1$ is given in~\eref{eq:dual_heat}. We set its first variation to zero and follow the steps presented for the transient 
convection-diffusion problem. The stiffness matrix is obtained by
setting $\alpha = 0$ in~\eref{eq:Kd=f_b} and the force vector is:
\begin{equation}\label{eq:f_heat}
\vm{f}_\lambda = - \int_0^1  u_0(x) \vm{N}_\lambda^\top (x,0)\, dx , \quad
\vm{f}_\mu = - \int_0^1 \vm{N}_\mu^\top (0,t) \, dt.
\end{equation}

The trial functions are chosen to be of the form~\eref{eq:bsp_2D_mulambda} with $n = 1$ in the 
bivariate open knot vector~\eref{eq:knot_2D}. Plots of the two-dimensional basis 
functions for $\mu_\theta(x,t)$
($p = 5$) and for $\lambda_\theta(x,t)$ ($q = 6$) are shown 
in~\fref{fig:bsp_2D_mulambda_heat}. Observe that
$\mu_\theta(1,t) = 0$ and $\lambda_\theta(0,t) 
= \lambda_\theta (x,1) = 0$ so
that $\lambda_\theta \in S_\lambda$ and $\mu_\theta \in S_\mu$
meet~\eref{eq:dual_heat_c}, and are therefore kinematically
admissible.
\begin{figure}[!tbh]
\centering
\begin{subfigure}{0.48\textwidth}
\includegraphics[width=\textwidth]{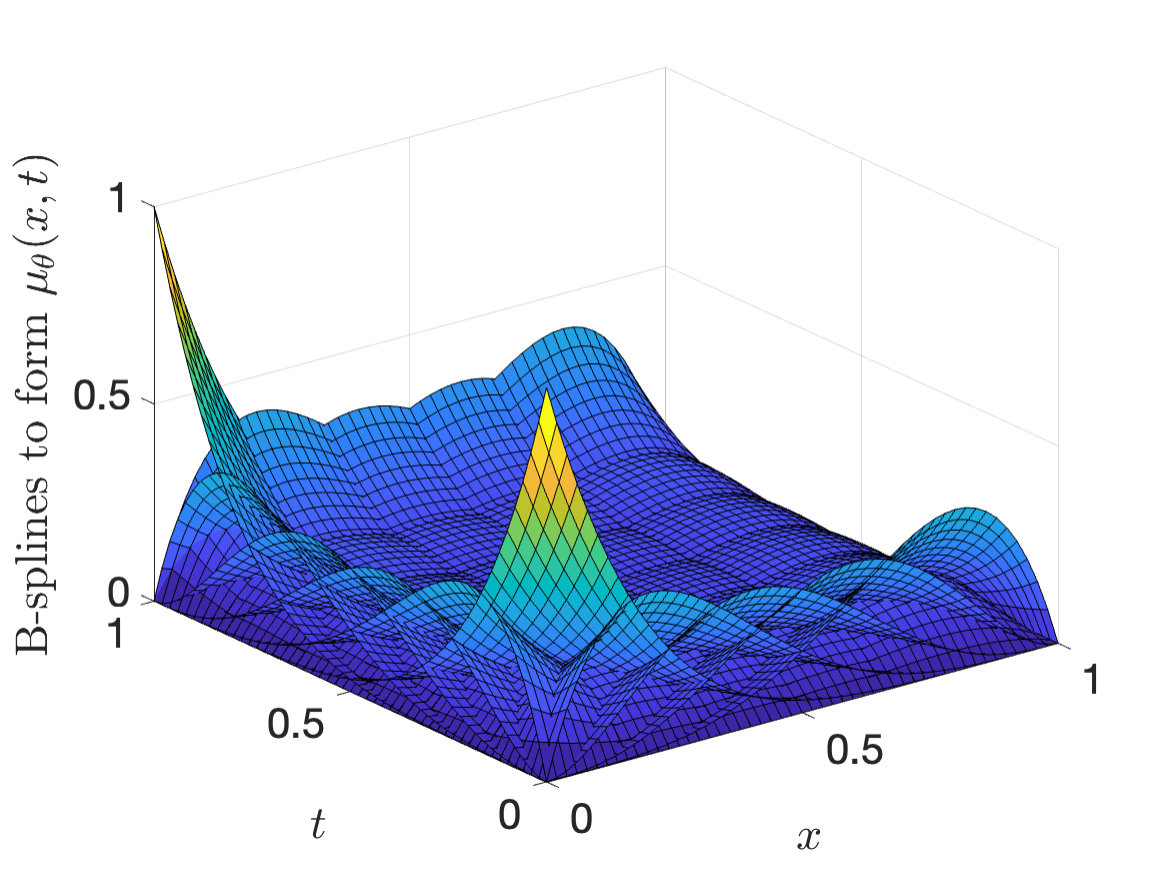}
\subcaption{}\label{fig:bsp_2D_mulambda_heat_a}
\end{subfigure} \hfill
\begin{subfigure}{0.48\textwidth}
\includegraphics[width=\textwidth]{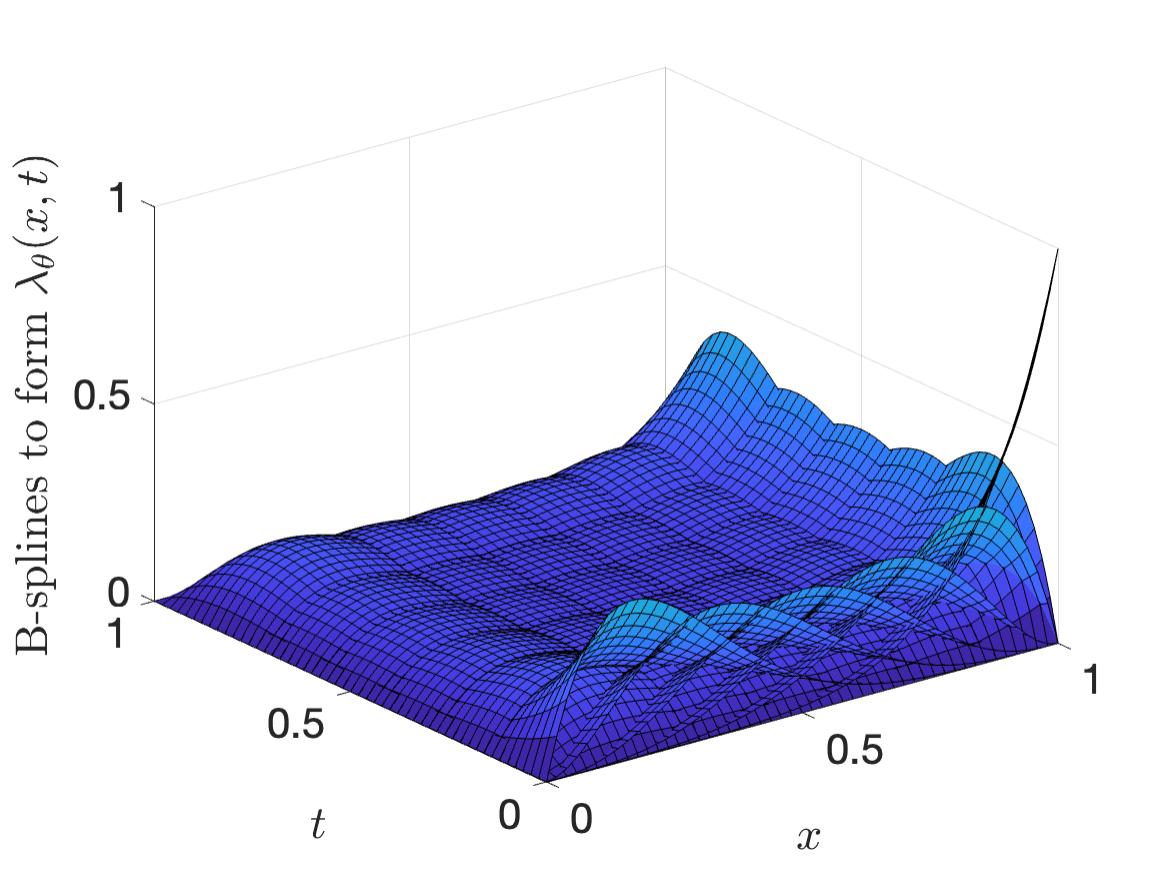}
\subcaption{}\label{fig:bsp_2D_mulambda_heat_b}
\end{subfigure}
\caption{Plots of two-dimensional B-spline basis functions
         to form the dual fields (a) $\mu_\theta(x,t)$ ($p = 5$) and 
         (b) $\lambda_\theta(x,t)$ ($q = 6$) to solve the heat 
         conduction problem. The bivariate open knot vector
         ($n = 1$) that is
         given in~\eref{eq:knot_2D} is used.}
         \label{fig:bsp_2D_mulambda_heat}
\end{figure}
We conduct numerical computations for $\kappa = 1$.  
The B-spline solution for the dual fields, 
$\mu_\theta(x,t)$ and
$\lambda_\theta(x,t)$, are presented 
in~\fref{fig:mulambda_heat}.  The primal fields
are computed from the dual field using the
DtP map in~\eqref{eq:dual_heat_a}.
 Comparisons of the
exact solution to the B-spline 
solution for the temperature and the flux 
are
provided in~\fref{fig:kappa1_p5q6_bsp_heat}. Even on 
this coarse discretization, the space-time
B-spline solutions are fairly accurate: maximum pointwise 
errors in $u$ and $q = \partial u / \partial x$ (flux)
are $4 \times 10^{-3}$ and $9 \times 10^{-2}$, 
respectively. One can observe 
from Figs.~\ref{fig:kappa1_p5q6_bsp_heat_i} and~\ref{fig:kappa1_p5q6_bsp_heat_j} that the accuracy
for notably $q = u^\prime$ worsens 
as one approaches the terminal time $t = 1$.  This behavior
in the vicinity of $t = 1$ is also observed 
(albeit more severe)
in~\fref{fig:bsp_transient_cd} 
for the convection-diffusion problem; refer to the discussion 
in~\sref{subsubsec:T_issues_cd} 
on the behavior of the discrete solution near $t = 1$.
\begin{figure}
\centering
\begin{subfigure}{0.48\textwidth}
\includegraphics[width=\textwidth]{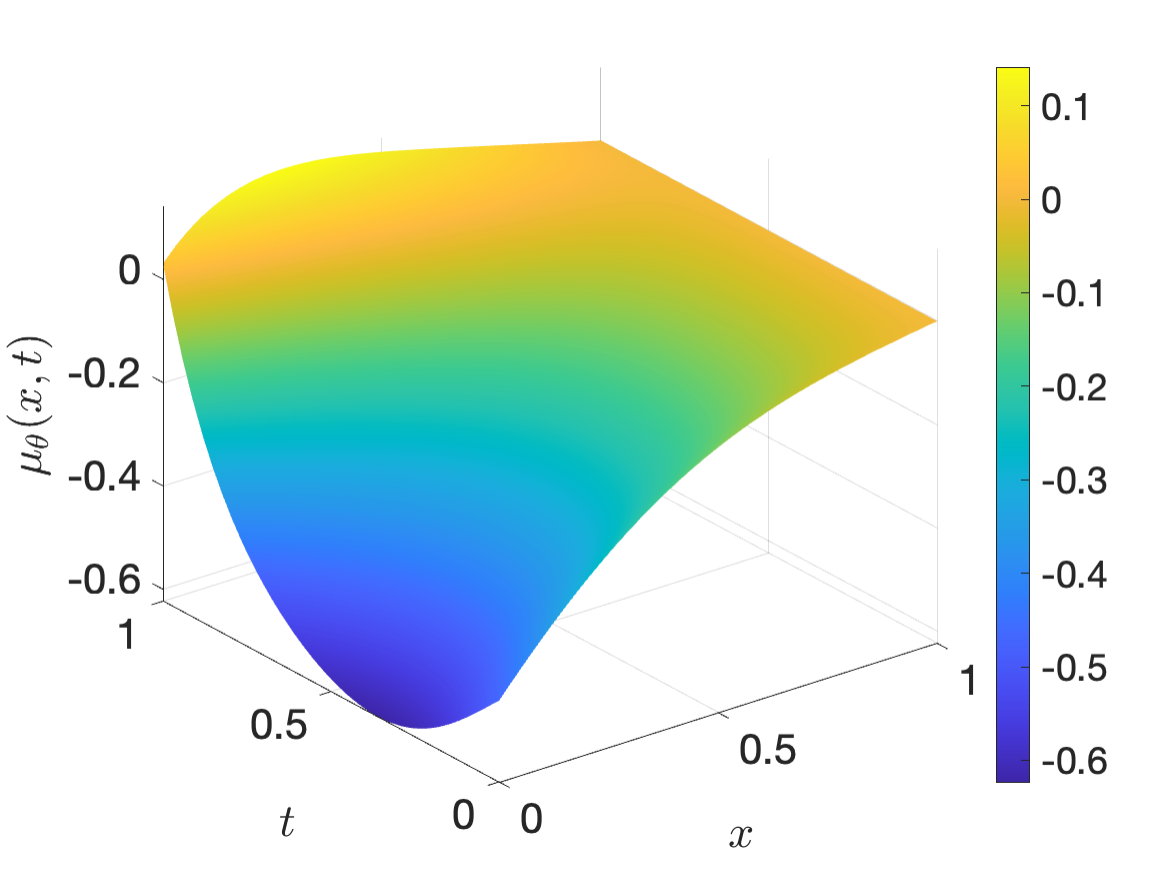}
\subcaption{}\label{fig:mulambda_heat_a}
\end{subfigure} \hfill
\begin{subfigure}{0.48\textwidth}
\includegraphics[width=\textwidth]{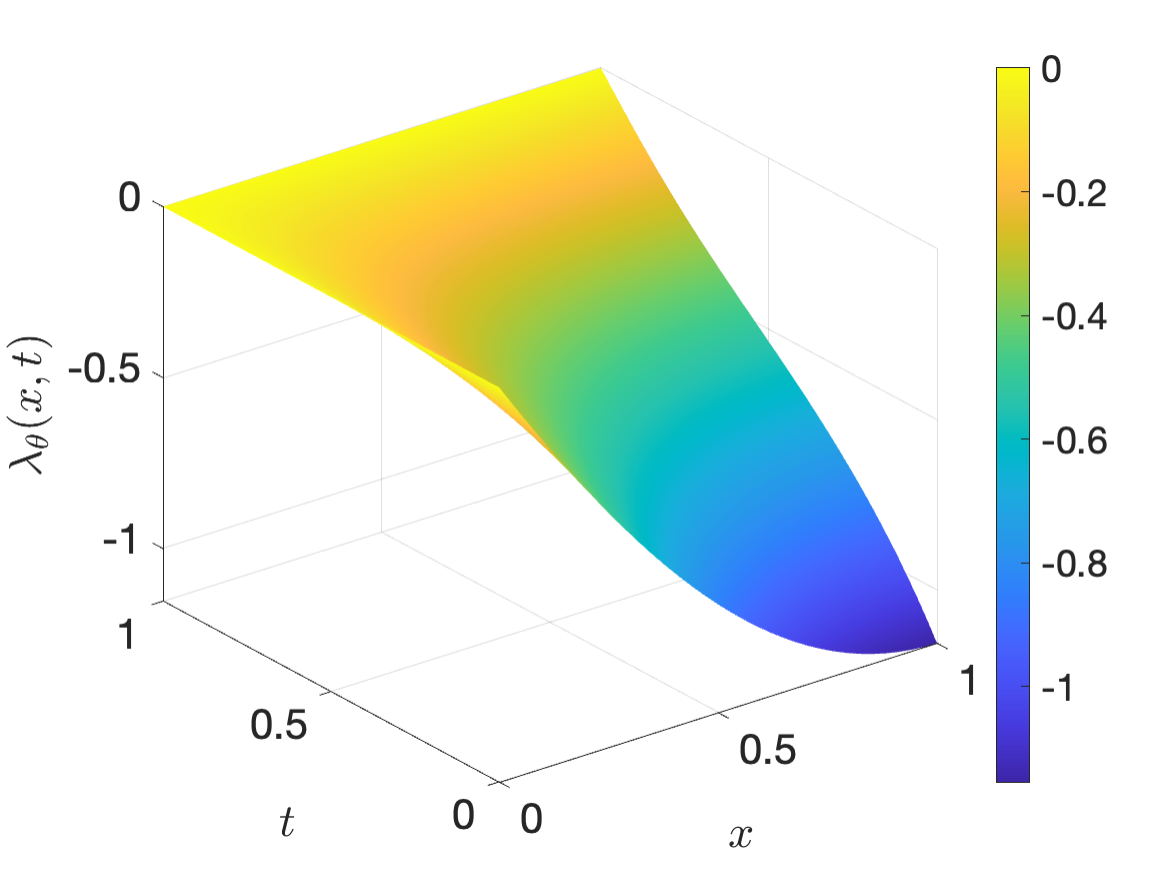}
\subcaption{}\label{fig:mulambda_heat_b}
\end{subfigure}
\caption{B-spline solution (dual fields) of the 
         transient heat equation $(\kappa = 1)$.
         (a) $\mu_\theta(x,t)$ ($p = 5$) and 
         (b) $\lambda_\theta(x,t)$ ($q = 6$).
         }
         \label{fig:mulambda_heat}
\end{figure}
\begin{figure}[!tbh]
\centering
\begin{subfigure}{0.33\textwidth}
\includegraphics[width=\textwidth]{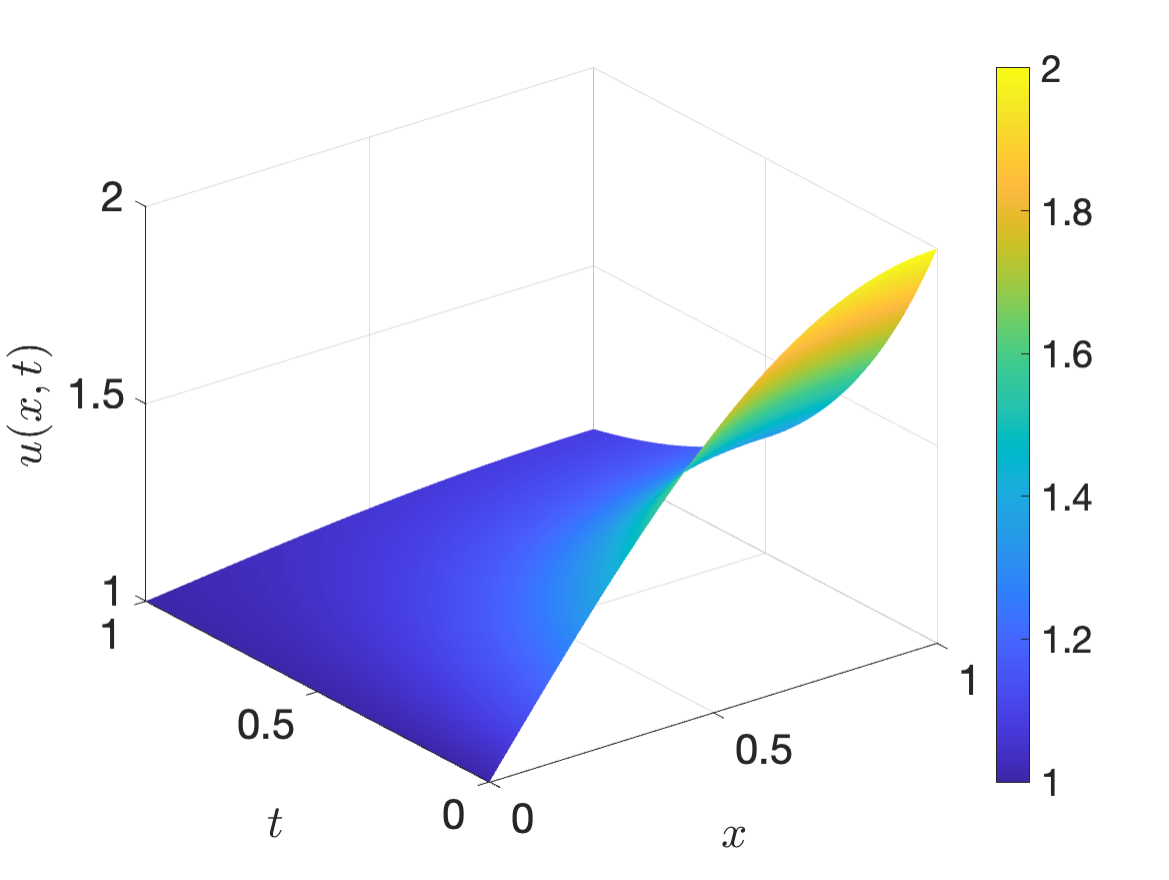}
\subcaption{}\label{fig:kappa1_p5q6_bsp_heat_a}
\end{subfigure} \hfill
\begin{subfigure}{0.33\textwidth}
\includegraphics[width=\textwidth]{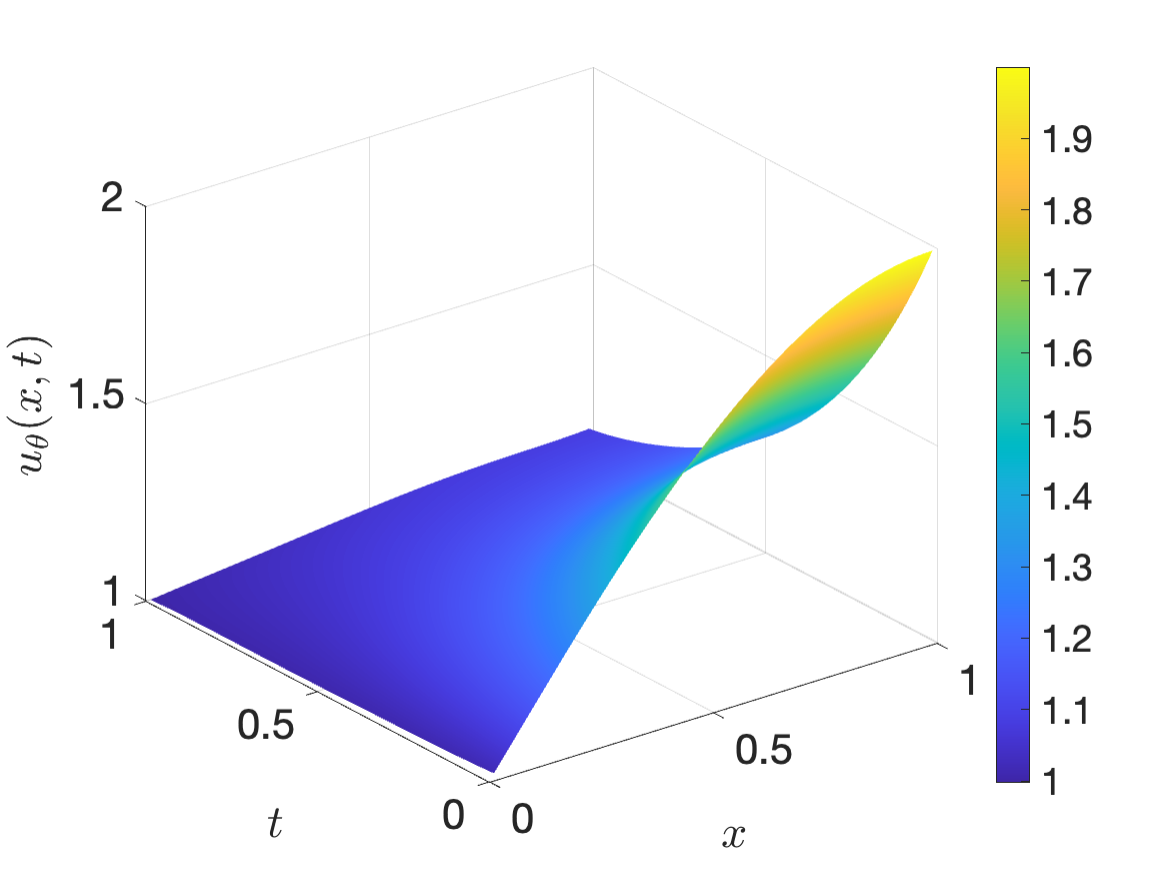}
\subcaption{}\label{fig:kappa1_p5q6_bsp_heat_b}
\end{subfigure} \hfill
\begin{subfigure}{0.33\textwidth}
\includegraphics[width=\textwidth]{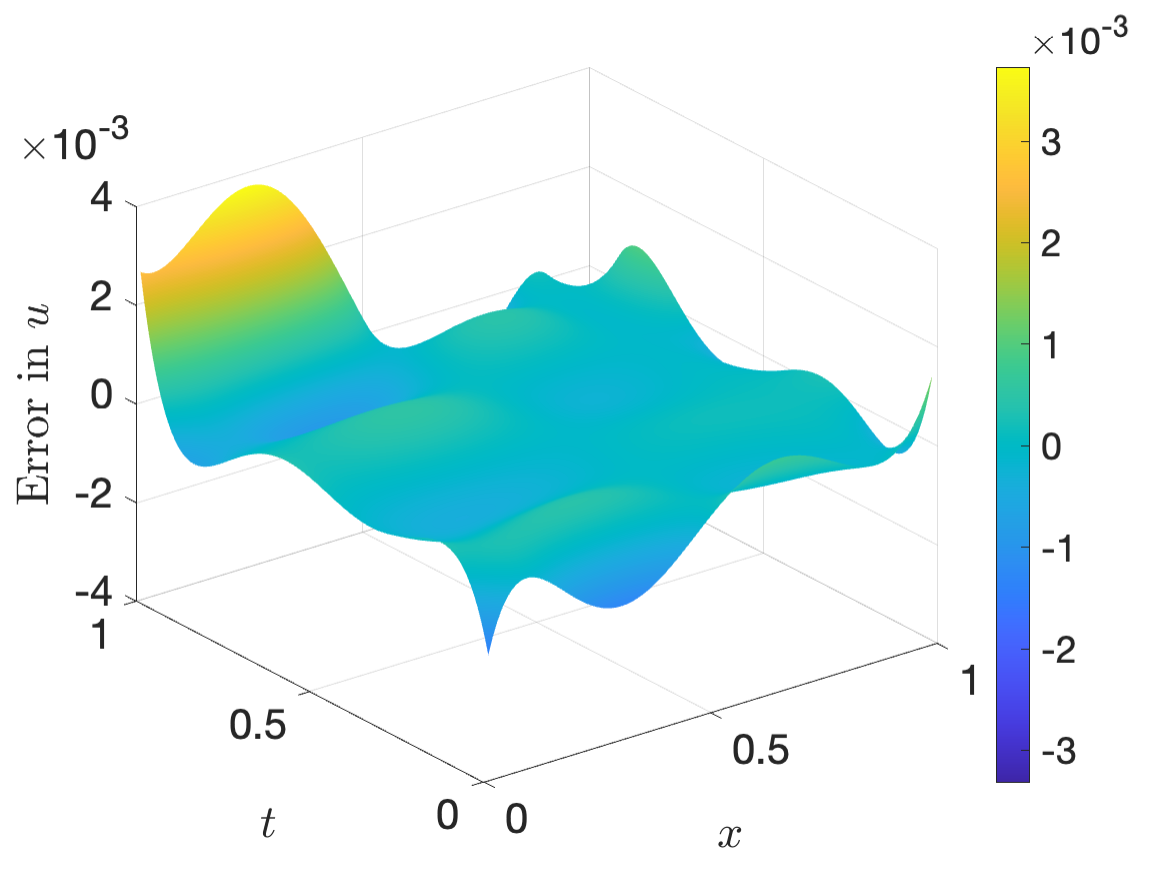}
\subcaption{}\label{fig:kappa1_p5q6_bsp_heat_c}
\end{subfigure}
\begin{subfigure}{0.33\textwidth}
\includegraphics[width=\textwidth]{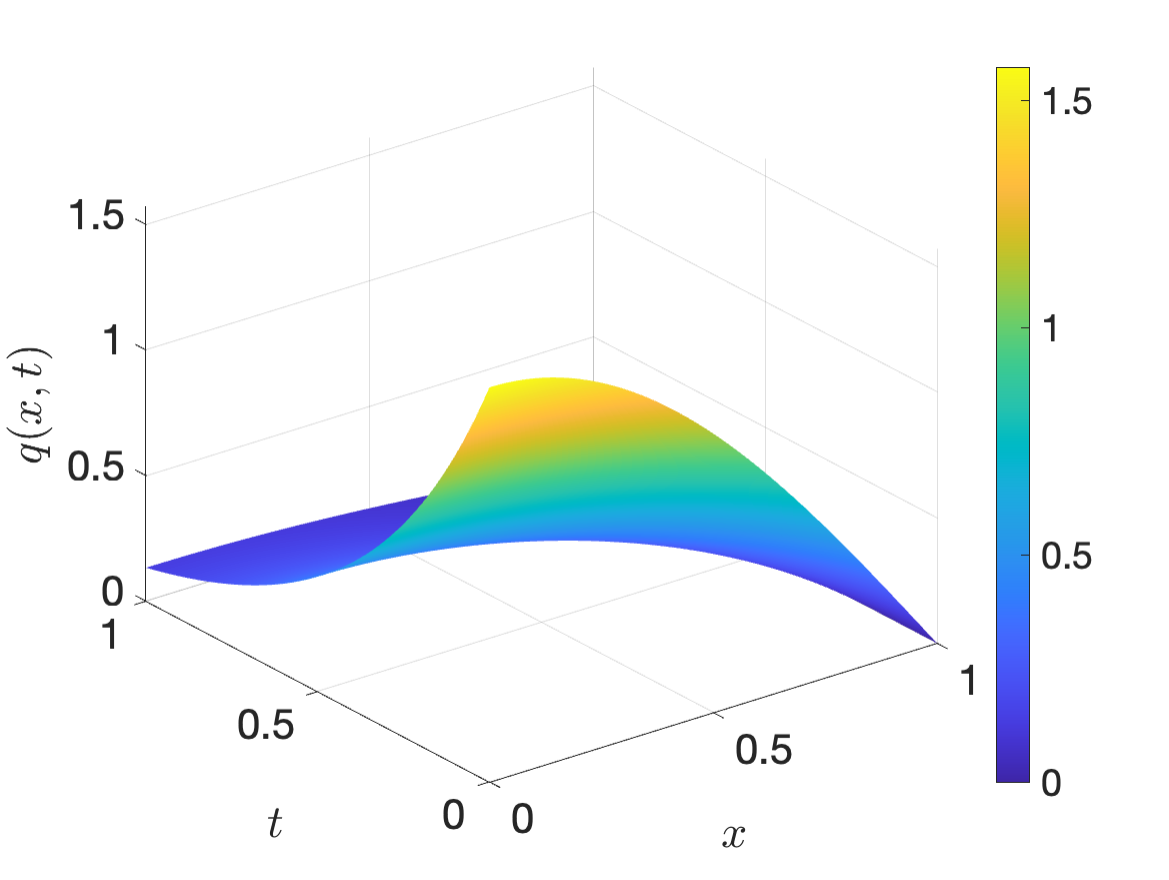}
\subcaption{}\label{fig:kappa1_p5q6_bsp_heat_d}
\end{subfigure} \hfill
\begin{subfigure}{0.33\textwidth}
\includegraphics[width=\textwidth]{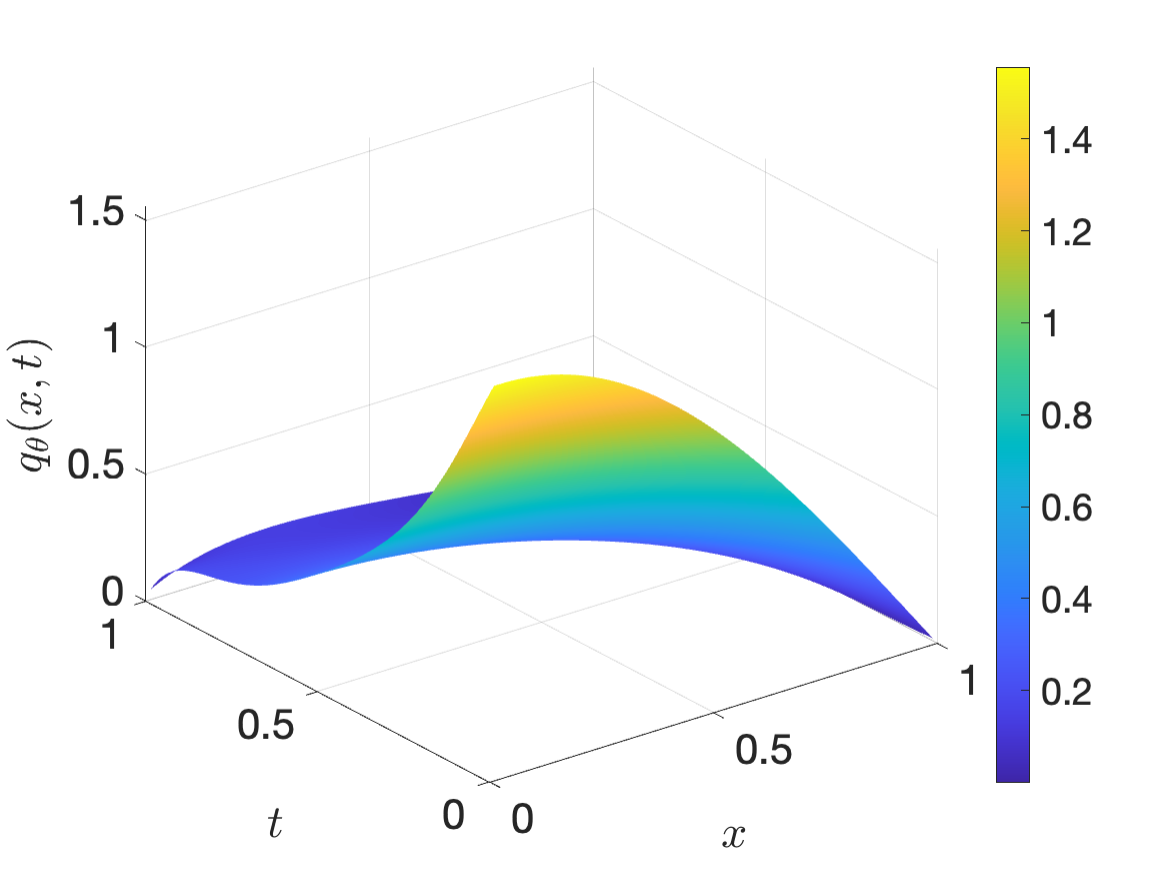}
\subcaption{}\label{fig:kappa1_p5q6_bsp_heat_e}
\end{subfigure} \hfill
\begin{subfigure}{0.33\textwidth}
\includegraphics[width=\textwidth]{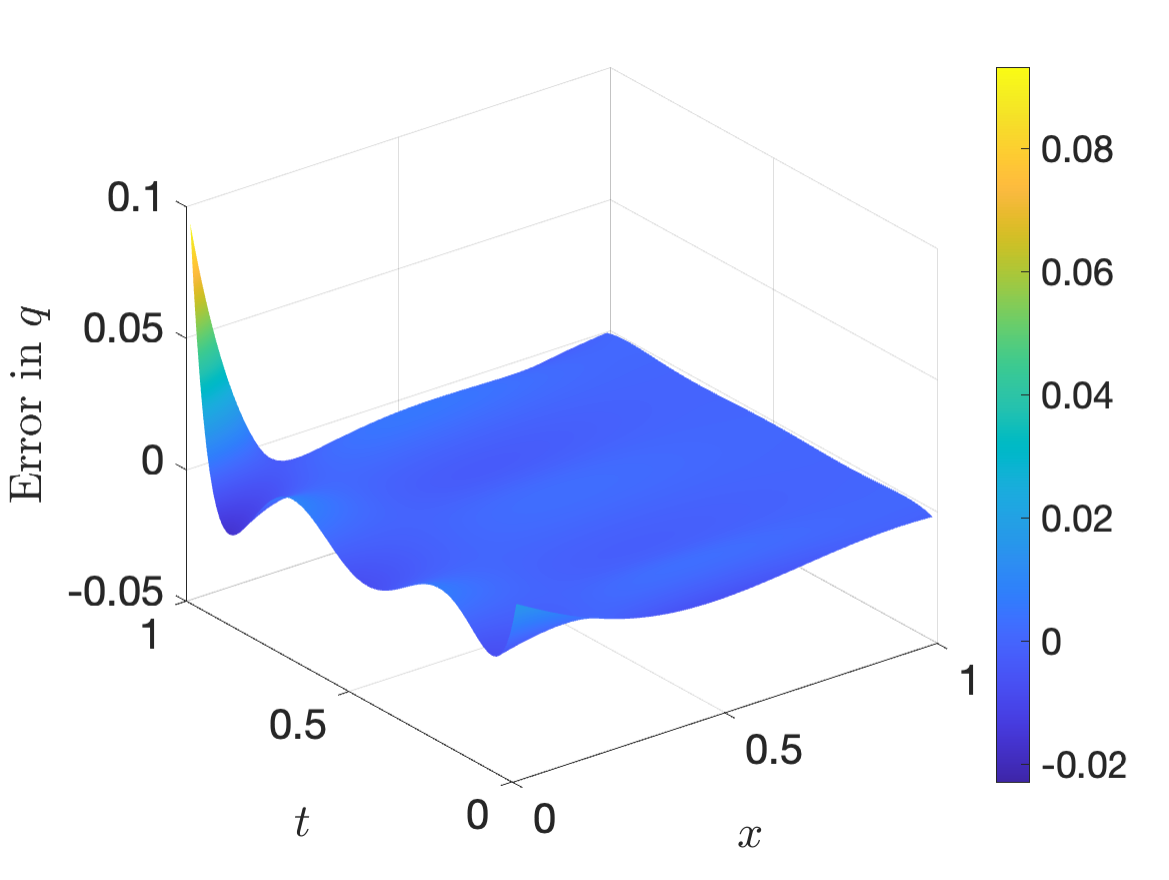}
\subcaption{}\label{fig:kappa1_p5q6_bsp_heat_f}
\end{subfigure}
\begin{subfigure}{0.24\textwidth}
\includegraphics[width=\textwidth]{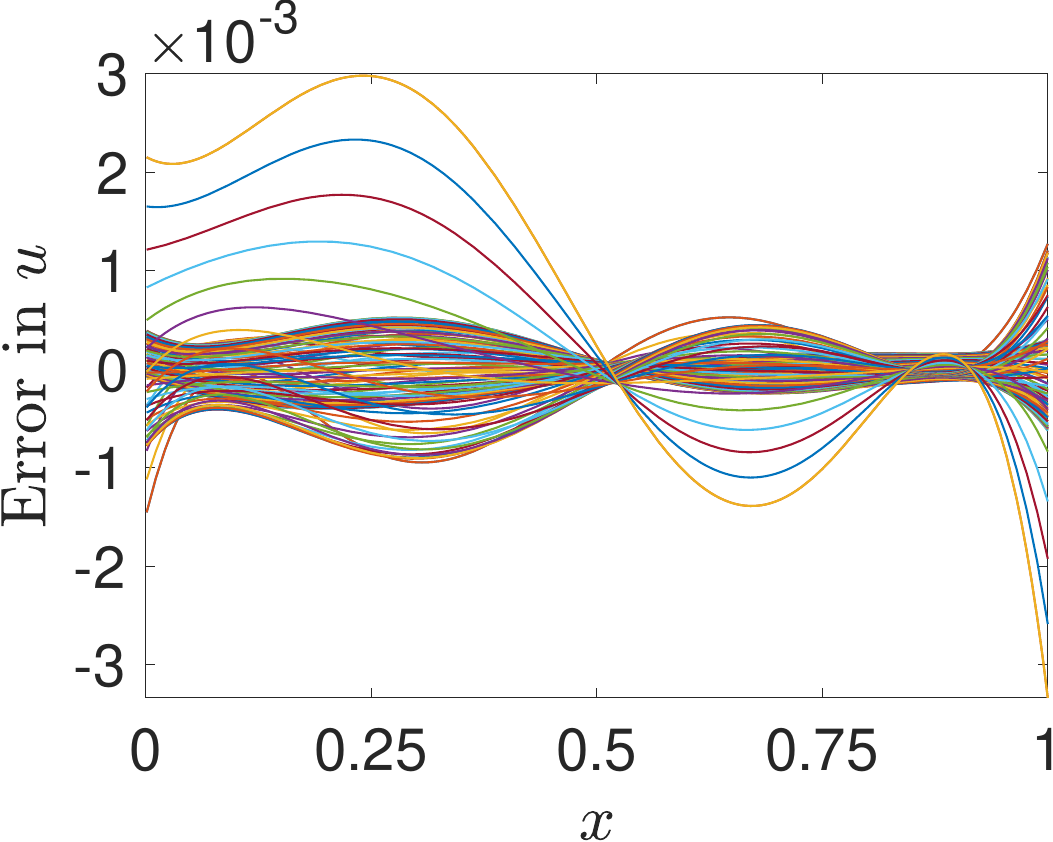}
\subcaption{}\label{fig:kappa1_p5q6_bsp_heat_g}
\end{subfigure} \hfill
\begin{subfigure}{0.24\textwidth}
\includegraphics[width=\textwidth]{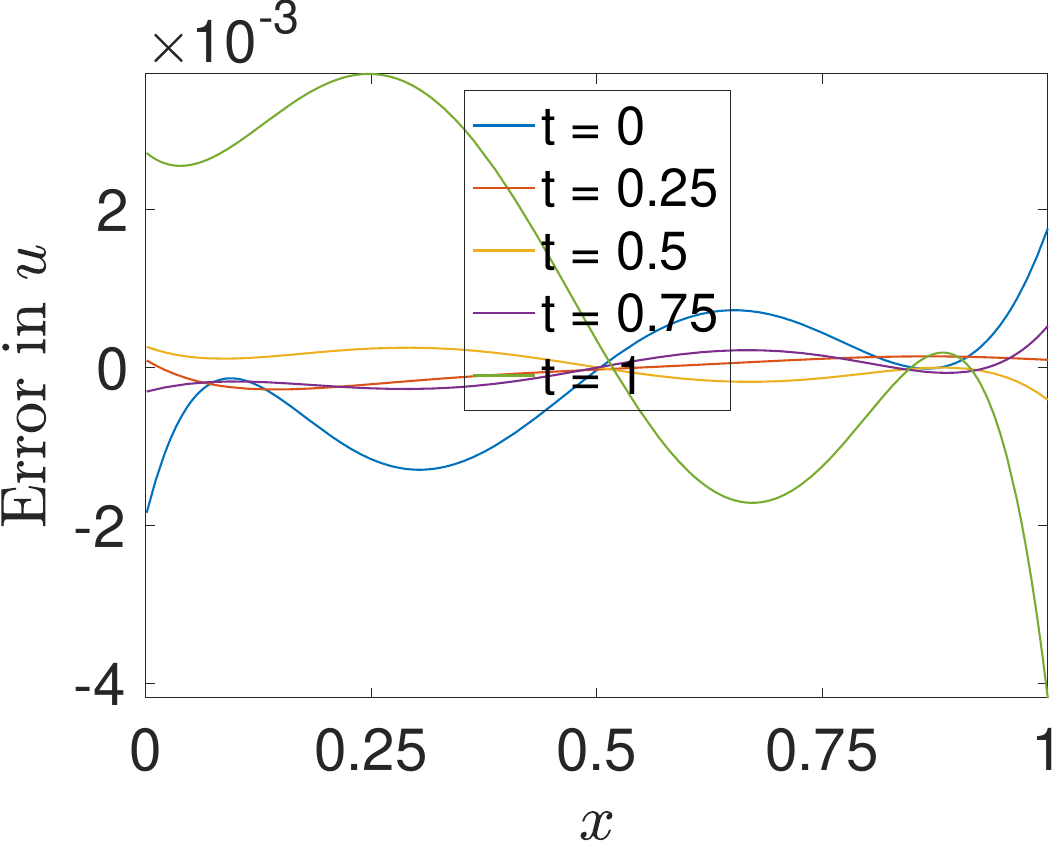}
\subcaption{}\label{fig:kappa1_p5q6_bsp_heat_h}
\end{subfigure} \hfill
\begin{subfigure}{0.24\textwidth}
\includegraphics[width=\textwidth]{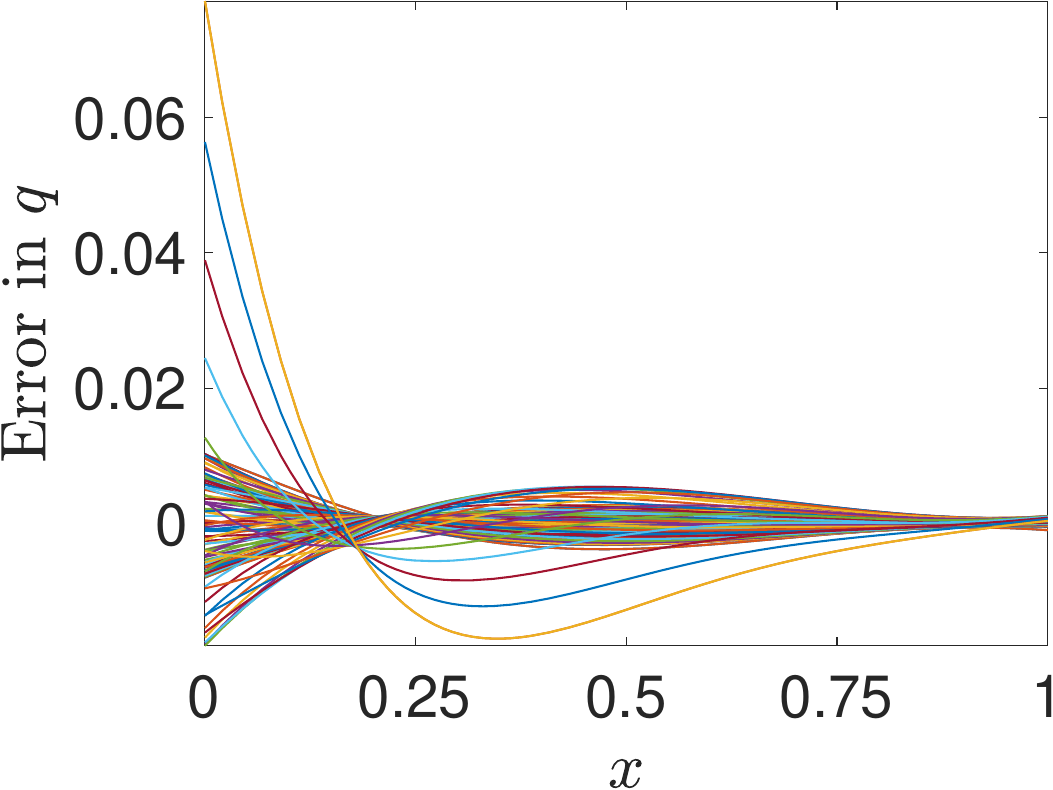}
\subcaption{}\label{fig:kappa1_p5q6_bsp_heat_i}
\end{subfigure} \hfill
\begin{subfigure}{0.24\textwidth}
\includegraphics[width=\textwidth]{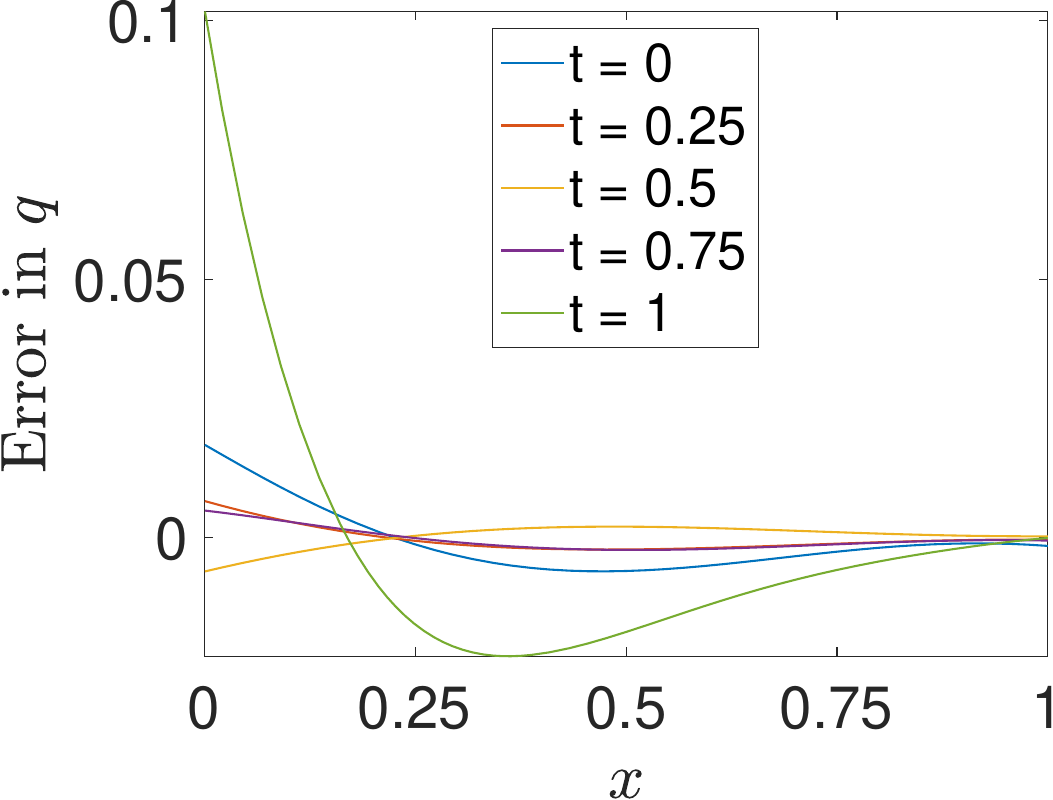}
\subcaption{}\label{fig:kappa1_p5q6_bsp_heat_j}
\end{subfigure}
\caption{Space-time B-splines to solve the
         transient heat conduction problem ($\kappa = 1$).
         The dual fields $\mu_\theta(x,t)$
         and $\lambda_\theta(x,t)$ are composed of tensor-product
         B-splines of degree $p = 5$ and $q = 6$, respectively.
         (a) Exact solution, $u$; (b) $u_\theta$; 
         (c) Error, $u - u_\theta$;
         (d) Exact $q = \partial u / \partial x$; 
         (e) $q_\theta$; and
         (f) Error, $q_\theta - \partial u 
         / \partial x$.
         Plots of the time history for the error in
         $u$ at
         (g) $t = [0,0.01,0.02,\dots,1]$
         and (h) $t = [0, 0.25, 0.5, 0.75, 1]$.
         Plots of the time history for the
         error in $q$ at
         (i) $t = [0,0.01,0.02,\dots,1]$
         and (j) $t = [0, 0.25, 0.5, 0.75, 1]$.}
         \label{fig:kappa1_p5q6_bsp_heat}
\end{figure}

\section{Conclusions}\label{sec:conclusions}
In this paper, we applied a recently proposed dual variational 
principle to solve partial differential equations that do not have a variational structure in primal form.  In this variational approach, the primal 
partial differential equation 
is treated as a constraint and an arbitrarily chosen convex auxiliary potential is minimized. This leads to requiring a concave dual functional to be 
maximized subject to Dirichlet boundary conditions on dual 
variables. Prescribed Dirichlet boundary conditions in the primal problem appear as natural boundary conditions in the dual functional.
The overarching goal of the duality approach is to solve 
ordinary/partial differential equations using a variational strategy.
\revised{We point out that 
nonlinear ODEs/PDEs when treated by
the dual approach become problems of convex optimization---formal theoretical justification is given in Reference~\cite{Acharya:2024:HCC}
and rigorous analysis is provided in References~\cite{Singh:2024:HCN,Acharya:2024:VDS,acharya2024variational}.}
Beyond its use in the present work and previous studies~\cite{Acharya:2022:VPN, Acharya:2023:DVP, Acharya:2024:HCC,Singh:2024:HCN,Acharya:2024:VDS,
Acharya:2024:APD}, one can also solve for critical points 
of primal functionals (for
example, the Chern--Simons functional~\cite{acharya2024variational}) 
that are not bounded below or above. 

We presented the general formalism and illustrated it for
a system of linear equations, a quadratic system of two equations,
and to obtain nonnegative solutions for an underdetermined system of 
linear equations.  It was also shown how a causal 
initial-value problem (ODE) 
can be exactly solved (and accurately approximated) as an acausal 
boundary-value problem 
in time. The dual weak form 
was derived for the transient convection-diffusion equation,  which was
discretized using a Galerkin method with finite-dimensional trial and test 
functions that were formed using machine learning approximants and 
B-splines. Uniqueness of solutions
of the dual variational equations for the transient convection-diffusion equation 
was established.  Shallow neural network with smooth RePU activation function
and univariate B-splines delivered accurate solutions
for the steady-state convection-diffusion equation. 
Smooth B-spline approximations of degrees $p$ and $p+1$ for 
the dual fields delivered more accurate solutions for the flux 
($q = u^\prime$) than $u$; 
\revised{optimal rates
of convergence of orders $p$ and $p+1$
were obtained in the $L^2$ norm and $H^1$ seminorm of the primal field $u$, respectively.}
One-dimensional transient (heat and convection-diffusion) 
problems were successfully solved as a space-time Galerkin method 
with tensor-product B-spline basis functions. 
It was shown that 
the errors were small but they tended to become relatively larger 
near the terminal
time $t = T$, and an explanation 
of this behavior and solution strategies for its remedy were put forth. 
As part of future work, we plan to pursue the 
application of the dual variational scheme
to solve challenging linear and nonlinear systems of ordinary and
partial differential equations.

\appendix

\section{Adjoint method for constrained-ODE}\label{appendix:A}
Consider the ODE in~\eref{eq:ODE_a}, and denote
the initial condition in~\eref{eq:ODE_b}
as $u(0) = u_0 := p$ ($p$ is now considered as a parameter),
so that the exact solution is:
\begin{equation}\label{eq:u_exact_for_adj}
u(t,p) = p e^{at}.
\end{equation}
Consider the objective function~\cite{Bradley:2024:TAM}
\begin{equation}\label{eq:adj_objective}
F(u,p) = \int_0^T u(t,p) \, dt.
\end{equation}
The sensitivity of
$F$ with respect to $p$,
$dF/dp$, is sought. On using the method of Lagrange multipliers, we
form the Lagrangian:
\begin{equation}
L(u,\lambda) = F(u,p) - 
\int_0^T \lambda \, (\dot{u} - a u) \, dt ,
\end{equation}
where $\lambda$ is known as the adjoint variable. 
Since $au - \dot{u} = 0$, then
\begin{equation}\label{eq:dfdp_0}
\begin{aligned}
\dfrac{dF}{dp} = \frac{dL}{dp} 
&= \int_0^T \dfrac{\partial u}{\partial p} \, dt
+
\int_0^T \lambda \left( 
a \dfrac{\partial u}{\partial p}
- \dfrac{\partial \dot{u} }{\partial p} 
\right) \, dt \\
&= \int_0^T \dfrac{\partial u}{\partial p} \, dt
- \left[ \lambda \dfrac{\partial u}{\partial p} \right]_{t = 0}^{t = T}
+ \int_0^T (\dot{\lambda} + a \lambda) \frac{\partial u}{\partial p} 
\, dt , \\
&= - \left[ \lambda \dfrac{\partial u}{\partial p} \right]_{t = 0}^{t = T}
+ \int_0^T (\dot{\lambda} + a \lambda + 1) \frac{\partial u}{\partial p} 
\, dt .
\end{aligned}
\end{equation}
Now we choose $\lambda(T) = 0$ so that the
boundary term at $t = T$ vanishes, and we eliminate the
integral by requiring $\lambda$ to satisfy the ODE 
(adjoint equation) $\dot{\lambda} + a \lambda + 1 = 0$. In addition, since $du(0)/dp = 1$ for the boundary term at 
$t = 0$, \eref{eq:dfdp_0} simplifies to
\begin{equation}\label{eq:dfdp}
\dfrac{dF}{dp} = \lambda(0).
\end{equation}
Now, $\lambda$ must satisfy the adjoint system:
\begin{equation}\label{eq:adj_system}
\dot{\lambda} + a \lambda + 1 = 0, \ \ \lambda(T) = 0,
\end{equation}
which admits the exact solution 
\begin{equation}
\lambda(t) = \dfrac{e^{a(T-t)} - 1}{a}.
\end{equation}
On using this in~\eref{eq:dfdp},
the solution from the adjoint method is:
\begin{equation}
\dfrac{dF}{dp} = \dfrac{e^{aT}- 1}{a},
\end{equation}
which matches the result for $dF/dp$ using the
exact solution in~\eqref{eq:u_exact_for_adj}. 

In the optimization of parametric 
differential-algebraic equation (DAE) systems,
a given objective function is required to be optimized with
respect to parameters $\vm{p}$, subject to
satisfying an initial-boundary value problem (forward or primal problem) that depends on the parameters. It is assumed that solvers
are available for the forward problem. In the dual approach for differential-algebraic systems, a user-specified convex potential
is optimized subject to treating the system as a constraint.
The potential is easily designed using 
the knowledge of the DAE system to be solved.  
Both methods use Lagrange multipliers (referred to as \emph{adjoint variable} and \emph{dual variable} in the
two approaches) to form the Lagrangian of the respective problems. The goal of the adjoint method is to provide an efficient method to compute the
sensitivity of the objective function with respect to $\vm{p}$.
The goal of the duality approach is to solve 
ordinary/partial differential-algebraic  
systems using a variational strategy.
This eases the burden of solving problems with nonstandard structure~\cite{Acharya:2023:DVP} or those
without existence-of-solution guarantees when viewed 
through the lens of current
methods for the primal problem~\cite{Singh:2024:HCN}. The adjoint method in constrained optimization does not concern itself with solution strategies for the forward problem, 
the availability of a robust technique for which is an essential ingredient of executing the method.
The dual variational approach to DAEs concerns itself with a strategy to solve such (primal) systems, but does not address the question of optimization of an objective function with respect to
parameters defining the primal problem. 

There is a common step (albeit superficial) 
in the two approaches. 
For parametric constrained function optimization, the scheme in 
the adjoint method is as follows~\cite{Bradley:2024:TAM}:
\begin{equation}\label{eq:appendix_adj}
\begin{aligned}
 \textrm{Given } f(\vx,p), \  \textrm{subject to} \ 
 \vm{g}(\vx,p) = 0,& \quad \mbox{find } 
\dfrac{df}{dp} .  \\
 \textrm{Crucial step is the definition of the adjoint $\vm{\lambda}$ from }& \quad 
 \dfrac{\partial f}{\partial \vx} + \vm{\lambda}^\top 
 \dfrac{\partial \vm{g} }{\partial \vx} = \vm{0} .   \\
 \textrm{Evaluate }& \quad \dfrac{df}{dp} 
 = \vm{\lambda}^\top \dfrac{\partial 
    \vm{g}}{\partial p} +
    \dfrac{\partial f}{\partial p}. 
\end{aligned}
\end{equation}
With the replacement $f \to H$, $\vm{g} \to \vm{G}$, 
and $\vx \to \vm{U}$
(see~\sref{sec:formalism}), the equation for 
the adjoint in~\eref{eq:appendix_adj} is identical to the
DtP mapping, but with a very different goal 
and interpretation. The adjoint method defines the adjoint
(dual variable $\vm{\lambda}$) using the knowledge of the primal 
$\vx$; the dual approach does exactly the opposite. In constrained optimization using the adjoint method, 
the solution of the primal problem is required to solve for the
adjoint variable. The dual approach defines the dual problem,
which provides the dual solution, and the DtP (dual-to-primal) mapping
yields the primal solution that solves the primal problem.

\section{Weak form of the dual BVP for the IVP}\label{appendix:B}
We consider the IVP posed in~\eqref{eq:ODE}, with the dual functional 
$S[\lambda]$ given in~\eqref{eq:S_IVP} and the strong form of the dual BVP presented in~\eref{eq:dual_ODE}.  On setting $\delta S[\lambda;
\delta \lambda] = 0$ and using~\eref{eq:uH}, we obtain
\begin{align}\label{eq:varform_ODE_BVP_0}
\begin{split}
\delta S[\lambda; \delta \lambda] &= 
- \int_0^T (\dot{\lambda} + a \lambda) (\dot{\delta \lambda} + a \delta \lambda) \,
dt - u_0 \delta \lambda(0) \\
&= - \int_0^T [ \dot{\lambda} \delta \dot{\lambda}   +
a^2 \lambda \delta \lambda ] \, dt  - 
\int_0^T  a [ \dot{\lambda} \delta \lambda + \lambda \delta \dot{\lambda} ] \, dt
 - u_0 \delta \lambda(0) \\
 &= 0 \ \ \forall \delta \lambda \in {\sf S}_\lambda,
\end{split}
\end{align}
where ${\sf S}_\lambda = \{ \delta \lambda \in H^1(0,T), \ \delta \lambda (T) = 0\}$ 
is the space of admissible variations.  We find that
\begin{align*}
\int_0^T  a ( \dot{\lambda} \delta \lambda + \lambda \delta \dot{\lambda} ) \, dt
&= \int_0^T a \delta \left[ \frac{d}{dt} \left( \frac{1}{2} \lambda^2 \right) \right] dt
= \delta \left[ \int_0^T a \frac{d}{dt} \left( \frac{1}{2} \lambda^2 \right) dt
\right]
= \delta \left[ a \frac{\lambda^2}{2} \right]_{t = 0}^{t=T} \\
&= a \lambda(T) \delta \lambda(T) - a \lambda(0) \delta \lambda (0)
= - a \lambda(0) \delta \lambda (0) ,
\end{align*}
since $\delta \lambda (T) = 0$. On using the above in~\eqref{eq:varform_ODE_BVP_0},
we can write the variational form as:
\begin{equation}\label{eq:varform_ODE_BVP}
\int_0^T [ \dot{\lambda} \delta \dot{\lambda} +  
a^2 \lambda \delta \lambda] \, dt - a \lambda(0) \delta \lambda(0) 
= -u_0 \delta \lambda(0) \ \ \forall \delta \lambda \in {\sf S}_\lambda .
\end{equation}
 
To derive the weak form, 
we multiply the dual BVP in~\eqref{eq:dual_ODE}
by the test function (virtual field) $\delta \lambda$ and perform integration by parts to obtain
\begin{equation*}
- \int_0^T [ \dot{\lambda} \delta \dot{\lambda} + a^2 \lambda \delta \lambda ] \, dt
+ \dot{\lambda}(T) \delta \lambda (T) 
- \dot{\lambda} (0) \delta \lambda (0) = 0 \ \
\forall \delta \lambda \in {\sf S}_\lambda .
\end{equation*}
On substituting $ \delta \lambda(T) = 0$ and the Robin boundary condition 
from~\eref{eq:dual_ODE_b} to replace $\dot{\lambda}(0)$ in the above equation leads us
to the weak form, which is identical to the variational form 
given in~\eref{eq:varform_ODE_BVP}.

\end{document}